\documentclass[12pt]{article}

\usepackage{latexsym, 
amscd, 
graphicx, amssymb}

\newtheorem{thm}{Theorem}[section]

\newtheorem{prop}[thm]{Proposition}
\newtheorem{cor}[thm]{Corollary}
\newtheorem{lemma}[thm]{Lemma}
\newtheorem{rema}[thm]{Remark}

\newcommand{\halmos}{\rule{1ex}{1.4ex}}

\newcommand{\binom}[2]{{{#1}\choose {#2}}}
\newcommand{\text}[1]{\mbox{\rm #1}}
\newcommand{\mod}{\;\;\mbox{\rm mod}\;}

\newcommand{\nn}{\nonumber \\}

 \newcommand{\res}{\mbox{\rm Res}}
 
\renewcommand{\hom}{\mbox{\rm Hom}}
\newcommand{\edo}{\mbox{\rm End}\ }
 \newcommand{\pf}{{\it Proof.}\hspace{2ex}}
 
 \newcommand{\epfv}{\hspace*{\fill}\mbox{$\halmos$}\vspace{1em}}

\newcommand{\wt}{\mbox{\rm wt}\ }
\newcommand{\swt}{\mbox{\rm {\scriptsize wt}}\ }

\newcommand{\tr}{\mbox{\rm Tr}}

\title{ {\bf Differential equations,  duality and 
modular invariance} }
\date{}
\author{Yi-Zhi Huang}

\begin{document}

\bibliographystyle{alpha}
\maketitle

\begin{abstract} 
We solve the problem of constructing all chiral genus-one correlation
functions from chiral genus-zero correlation functions associated to a
vertex operator algebra satisfying the following conditions: (i)
$V_{(n)}=0$ for $n<0$ and $V_{(0)}=\mathbb{C}\mathbf{1}$, (ii) every
$\mathbb{N}$-gradable weak $V$-module is completely reducible and (iii)
$V$ is $C_{2}$-cofinite.  We establish the fundamental properties of
these functions, including suitably formulated commutativity,
associativity and modular invariance.  The method we develop and use
here is completely different from the one previously used by Zhu and
other people. In particular, we show that the $q$-traces of products of
certain geometrically-modified intertwining operators satisfy modular
invariant systems of differential equations which, for any fixed modular
parameter, reduce to doubly-periodic systems with only regular singular
points. Together with the results obtained by the author in the
genus-zero case, the results of the present paper solves essentially the
problem of constructing chiral genus-one weakly conformal field
theories from the representations of a vertex operator algebra
satisfying the conditions above. 
\end{abstract}

\tableofcontents

\renewcommand{\theequation}{\thesection.\arabic{equation}}
\renewcommand{\thethm}{\thesection.\arabic{thm}}
\setcounter{equation}{0}
\setcounter{thm}{0}
\setcounter{section}{-1}

\section{Introduction}

The present paper is one of the main steps in the author's program of
constructing conformal field theories in the sense of Segal \cite{S1}
\cite{S2} \cite{S3} and Kontsevich from representations of suitable
vertex operator algebras.  We construct all the chiral genus-one
correlation functions from chiral genus-zero correlation functions (or
equivalently, from intertwining operators) for a vertex operator algebra
satisfying suitable natural conditions, and we establish their
fundamental properties, including suitably formulated commutativity,
associativity and modular invariance.  Together with the results previously
obtained by the author in the genus-zero case, the results of the
present paper solve essentially the problem of constructing chiral
genus-one weakly conformal field theories. The results of the present
paper imply the rigidity of the ribbon tensor category of
modules for a suitable vertex operator algebra.  Together with the
results in the genus-zero case and the results obtained by Moore and
Seiberg in \cite{MS1} and \cite{MS2}, the results of the present paper
also imply the Verlinde conjecture that the fusion rules are
diagonalized by the action of the modular transfermation given by
$\tau\mapsto -1/\tau$. These applications will be discussed in \cite{H12}
and \cite{H13}. We believe that the results of the present paper will be
interesting and useful not only to people working on vertex operator
algebras and conformal field theories, but also to people in other
related fields. 

In \cite{H3.5}, \cite{H4}, \cite{H5}, 
\cite{H8}, \cite{H10} and \cite{H11}, the author has constructed 
a chiral genus-zero weakly conformal field theory in the sense of Segal 
\cite{S2} \cite{S3} 
from representations of a vertex operator algebra satisfying the
$C_{1}$-cofiniteness condition and certain finite reductivity
conditions. The genus-zero part of the program is thus
more or less finished.
The next natural step is the construction of all the chiral genus-one
correlation functions and the proof of their basic properties.  
In \cite{Z}, Zhu constructed chiral
genus-one correlation functions associated to elements of a vertex operator
algebra $V$ from $V$-modules when $V$ satisfies the
$C_{2}$-cofiniteness condition and certain (stronger) finite
reductivity conditions, and he proved their basic properties, 
including modular invariance.  Zhu's work was generalized by Dong-Li-Mason
\cite{DLM} to the case of twisted
representations of vertex operator algebras and by Miyamoto \cite{M1}
to the case of chiral
genus-one correlation functions associated to elements of modules
among which at most one is not isomorphic to the algebra.
Miyamoto \cite{M2} further constructed certain one-point chiral
genus-one correlation functions involving logarithms 
when the vertex operator algebra 
satisfies only the $C_{2}$-cofiniteness condition. 

The method used by Zhu, Dong-Li-Mason and Miyamoto depends on the
commutator formula for (untwisted or twisted) vertex operators or the
commutator formula between vertex operators and intertwining operators.
These formulas were needed to find recurrence relations such that the
study of the formal $q$-traces of products of more than one operator can
be reduced to the study of the formal $q$-traces of vertex or
intertwining operators, not products of such operators.  

However, this method
cannot be generalized to give a construction of general chiral genus-one
correlation functions because there is no commutator formula for general
intertwining operators, instead see \cite{H1}. 
Note that chiral genus-one correlation functions
are the the main objects to construct in the genus-one case. Also, most
of the conjectures in chiral genus-one weakly conformal field theories
can be formulated and proved directly using these correlation functions.
For example, based on the assumptions that all the chiral genus-zero and
genus-one correlation functions have been constructed and the desired
properties (including the duality properties and the modular invariance)
hold, Moore and Seiberg in \cite{MS1} and \cite{MS2} derived the
genus-one coherence relations and proved the Verlinde conjecture \cite{V}
(the diagonalization of the fusion rules by the action of the modular
transformation $\tau\mapsto -1/\tau$). It turns out that the proof of the 
rigidity of the ribbon tensor category of modules for 
a suitable vertex operator algebra also needs such correlation functions
and their properties \cite{H13}.  In the present paper, we develop
a method completely different from the one used in \cite{Z}, \cite{DLM},
\cite{M1} and \cite{M2} and completely solve the problem of constructing
all the chiral genus-zero correlation functions and establishing all the 
desired properties.

Here is a brief description of our main results. For simplicity, here
we shall discuss the results under certain strong conditions, though,
as we show in this paper, many of
these results in fact hold under weaker conditions.  Let $V$ be a vertex
operator algebra satisfying the following conditions: (i) $V_{(n)}=0$
for $n<0$ and $V_{(0)}=\mathbb{C}\mathbf{1}$, (ii) every
$\mathbb{N}$-gradable weak $V$-module is completely reducible (see
Section 6 for a definition of $\mathbb{N}$-gradable weak $V$-module)
and (iii)
$V$ is $C_{2}$-cofinite (that is, $\dim V/C_{2}(V)<\infty$, where
$C_{2}(V)$ is the subspace of $V$ spanned by the elements of the form
$u_{-2}v$ for $u, v\in V$).  We prove that the $q$-traces of products
of certain geometrically-modified intertwining operators satisfy modular
invariant 
systems of different equations which reduce to doubly-periodic
systems with only regular
singular points when the modular parameter $\tau$ is fixed.  Using these
systems of differential equations, we prove the convergence of the
$q$-traces. The analytic extensions of these convergent $q$-traces
give chiral genus-one correlation functions.  We prove that these chiral
genus-one correlation functions satisfy suitable 
``genus-one commutativity'' and
``genus-one associativity'' properties.  
We also introduce a new product in a vertex
operator algebra and, using this new product, we construct and study an
associative algebra which looks very different from,
but is in fact isomorphic to, Zhu's algebra (see \cite{Z}). Using
these results, we prove the
modular invariance of the space of chiral genus-one correlation
functions. Geometrically,
this modular invariance and the double-periodicity of the reduced 
systems mentioned above give vector bundles with flat connections
over the moduli spaces of genus-one Riemann surfaces with punctures
and standard local coordinates vanishing at the punctures.

Together with the results obtained by the author in the genus-zero case,
the results of the present paper solves essentially the problem of
constructing chiral genus-one weakly conformal field theories.  The
results of the present paper also imply the rigidity of
the ribbon tensor category of modules for a vertex operator algebra. The
results of the present paper are the genus-one parts of the assumptions
on rational conformal field theories in Moore-Seiberg's important work
\cite{MS1} \cite{MS2}. Thus, together with the author's results in
\cite{H1}, \cite{H3} and \cite{H11} in the genus-zero case, the results
of the present paper imply all the results of Moore and Seiberg in
\cite{MS1} and \cite{MS2} obtained using only the genus-zero and
genus-one parts of their fundamental assumptions on rational conformal
field theories.  In particular, the the results of the present paper
together with the author's results in \cite{H1}, \cite{H3} and
\cite{H11} in the genus-zero case imply the Verlinde conjecture for
modules for a vertex operator algebra satisfying the conditions
mentioned above. The details of the proof of the rigidity 
of the ribbon tensor category and the proof of the Verlinde conjecture
will be given in the papers \cite{H12} and \cite{H13}, repectively. 

Our construction of the chiral genus-one correlation functions actually
verifies the sewing property which states that chiral genus-one
correlation functions can be obtained from chiral genus-zero correlation
functions by taking $q$-traces. In particular, we construct the spaces
of chiral genus-one correlation functions and the factorization property
for these spaces is an easy consequence of the sewing property.  We
would like to emphasize that the converse is not true, that is, even if
the spaces of chiral genus-one correlation functions are identified
abstractly and the factorization property is proved, additional
structures such as those constructed in the present paper are still
needed to construct the chiral genus-one correlation functions and to
prove the sewing property. 

Our modular invariance result gives, in particular, 
new proofs of the modular invariance 
result of Zhu in \cite{Z} and its direct, straightforward generalization
by Miyamoto in \cite{M1} since these results 
are (very) special cases. The differential equations and 
the duality properties obtained in the present paper also
have straightforward generalizations to the case of twisted 
modules and intertwining operators among them. Together with a
straightforward generalization of some of the results on
the associative algebras discussed in the present paper, 
they also give a modular invariance result in the twisted case. 
In particular, as a special case, this modular invariance 
gives a new proof of the 
modular invariance result of Dong, Li and Mason in \cite{DLM}. 
We shall discuss the twisted case elsewhere.

We assume that the reader is familiar with basic notions and results in
the algebraic and geometric theory of vertex operator algebras as
presented in \cite{FLM}, \cite{FHL} and \cite{H3.5}. We also assume that
the reader is familiar with the theory of intertwining operator algebras
as developed and presented by the author in \cite{H1}, \cite{H3},
\cite{H4}, \cite{H5}, \cite{H9}, \cite{H11}, based on the tensor product
theory developed by Lepowsky and the author in \cite{HL1}--\cite{HL4}
and \cite{H1}.  We do not, however, assume that the reader is familiar
with the modular invariance result proved by Zhu in \cite{Z} and the
method used there and in other papers, since our method is completely
different.

The present paper is organized as follows: In Section 1, we discuss 
geom\-etrically-modified intertwining operators, which play an 
important role in our construction. 
Then we prove some fundamental 
identities for the $q$-traces of products of 
these geometrically-modified intertwining operators in Section 2.
We establish the existence of the systems of differential equations 
in Section 3. In Section 4, we prove the convergence of the 
$q$-traces of products of 
geometrically-modified intertwining operators and construct
chiral genus-one correlation functions. In the same section, we also
establish genus-one commutativity and associativity for 
these chiral genus-one correlation functions. 
In Section 5, we prove that for fixed modular parameter $\tau$, 
the systems consisting of those equations in the systems mentioned above
not involving derivatives of $\tau$ have only regular
singular points. We introduce our new product on a vertex operator algebra and
construct and study the corresponding associative algebra in Section 6.
The modular invariance of the space of 
chiral genus-one correlation functions is proved in Section 7. 

\paragraph{Notations and conventions} In this paper, $i$ is either 
$\sqrt{-1}$ or an index, and it should be easy to tell which is 
which. The symbols $\mathbb{N}$, 
$\mathbb{Z}_{+}$, $\mathbb{Z}$, $\mathbb{Q}$, $\mathbb{R}$, $\mathbb{C}$
and $\mathbb{H}$ 
denote the nonnegative integers, positive integers, 
integers, rational numbers, real numbers, complex numbers and
the upper half plane, respectively. We shall 
use $x, , y, ...$ to denote commuting formal variables and 
$z, z_{1}, z_{2}, ...$ to denote complex numbers or complex variables.
For any nonzero $z\in \mathbb{C}$ and any $n\in \mathbb{Z}$, 
$\log z=\log |z|+i\arg z$ with $0\le \arg z <2\pi$ and 
$z^{n}=e^{n\log z}$. For any $z\in \mathbb{C}$, $q_{z}=e^{2\pi iz}$. 
As in \cite{FLM} and \cite{FHL}, for
a vertex operator algebra $(V, Y, \mathbf{1}, \omega)$ or a module
$(W, Y)$,
we use the $u_{n}$ for $u\in V$ and $n\in \mathbb{Z}$ 
to denote the coefficients of the formal series $Y(u, x)$. 
We shall use $w_{1}, w_{2}, ...$ to denote elements of $V$-modules. 
Note that when the $V$-module is the adjoint module, 
the notations $u_{n}$ and $w_{1}, w_{2}, ...$
clearly have different meanings. 

\paragraph{Acknowledgment}  This research 
is supported in part 
by NSF grant DMS-0070800. I am grateful to Jim Lepowsky and Antun Milas
for helpful comments.

\renewcommand{\theequation}{\thesection.\arabic{equation}}
\renewcommand{\thethm}{\thesection.\arabic{thm}}
\setcounter{equation}{0}
\setcounter{thm}{0}

\section{Geometrically-modified intertwining operators}

We need the theory of composition-invertible power series
and their actions on modules for the Virasoro algebra developed 
in \cite{H3.5}. Let $A_{j}$, $j\in \mathbb{Z}_{+}$, be the 
complex numbers defined 
by 
\begin{eqnarray*}
\frac{1}{2\pi i}\log(1+2\pi i y)=\left(\exp\left(\sum_{j\in \mathbb{Z}_{+}}
A_{j}y^{j+1}\frac{\partial}{\partial y}\right)\right)y.
\end{eqnarray*}
In particular, a simple calculation gives
\begin{eqnarray}
A_{1}&=&-\pi i,\label{a1}\\
A_{2}&=&-\frac{\pi^{2}}{3}.\label{a2}
\end{eqnarray}
Note that the composition inverse of $\frac{1}{2\pi i}\log(1+2\pi i y)$
is $\frac{1}{2\pi i}(e^{2\pi iy}-1)$ and thus we have
\begin{eqnarray*}
\frac{1}{2\pi i}(e^{2\pi iy}-1)=\left(\exp\left(-\sum_{j\in \mathbb{Z}_{+}}
A_{j}y^{j+1}\frac{\partial}{\partial y}\right)\right)y.
\end{eqnarray*}

In this paper, we fix a vertex operator algebra $V$.
We shall use $Y$ to denote the vertex operator maps for the algebra
$V$ and for $V$-modules. For any $V$-module 
$W$, we shall denote the operator 
$\sum_{j\in \mathbb{Z}_{+}}
A_{j}L(j)$
on $W$ by $L_{+}(A)$. Then 
$$e^{\displaystyle -\sum_{j\in \mathbb{Z}_{+}}
A_{j}L(j)}=e^{-L_{+}(A)}.$$
For convenience, we introduce the operator
\begin{equation}\label{ux}
\mathcal{U}(x)=(2\pi ix)^{L(0)}e^{-L^{+}(A)}\in (\mbox{\rm End}\;W)\{x\}
\end{equation}
where $(2\pi i)^{L(0)}=e^{(\log 2\pi +i \frac{\pi}{2})L(0)}$
(and we shall use this convention throughout this paper). 
In fact we can substitute for $x$ in the operator 
$\mathcal{U}(x)$ a number or any formal expression
which makes sense. For example,
$$\mathcal{U}(1)=(2\pi i)^{L(0)}e^{-L_{+}(A)}$$
and we have
\begin{equation}
\mathcal{U}(xy)=x^{L(0)}\mathcal{U}(y).\label{ux-u1}
\end{equation}
For any $z\in \mathbb{C}$, we shall use 
$q_{z}$ to denote $e^{2\pi iz}$.

\begin{lemma}\label{u1-omega}
For the conformal element $\omega=L(-2)\mathbf{1}\in V$,
$$\mathcal{U}(1)\omega=(2\pi i)^{2}
\left(\omega-\frac{c}{24}\mathbf{1}\right),$$
where $c$ is the central charge of $V$. 
\end{lemma}
\pf
By definition, we have
\begin{eqnarray*}
\mathcal{U}(1)\omega&=&(2\pi i)^{L(0)}e^{-L_{+}(A)}L(-2)\mathbf{1}\nn
&=&(2\pi i)^{L(0)}L(-2)\mathbf{1}-(2\pi i)^{L(0)}
A_{2}L(2)L(-2)\mathbf{1}\nn
&=&(2\pi i)^{2}L(-2)\mathbf{1}+\frac{\pi^{2}c}{6}\mathbf{1}\nn
&=&(2\pi i)^{2}\left(\omega-\frac{c}{24}\mathbf{1}\right), 
\end{eqnarray*}
where we have used (\ref{a1}) and (\ref{a2}). \epfv

We have the following conjugation and commutator formulas:

\begin{prop} 
Let $W_{1}, W_{2}$ and $W_{3}$ be $V$-modules and  $\mathcal{Y}$ 
an intertwining operator of type 
$\binom{W_{3}}{W_{1}W_{2}}$. 
Then for $u\in V$ and $w\in W_{1}$,
\begin{equation}\label{chg-var0}
\mathcal{U}(1)\mathcal{Y}(w, x)
(\mathcal{U}(1))^{-1}
=\mathcal{Y}(\mathcal{U}(e^{2\pi ix})w, e^{2\pi ix}-1)
\end{equation}
and
\begin{eqnarray}\label{x-comm0}
\lefteqn{[Y(\mathcal{U}(x_{1})u, x_{1}), \mathcal{Y}(\mathcal{U}(x_{2})w,
x_{2})]}\nn
&&=2\pi i\res_{y}\delta\left(\frac{x_{1}}{e^{2\pi iy}x_{2}}\right)
\mathcal{Y}(\mathcal{U}(x_{2})Y(u, y)w, 
x_{2}). 
\end{eqnarray} 
\end{prop} 
\pf 
Let $B_{j}$, $j\in \mathbb{Z}_{+}$, be rational numbers defined 
by 
\begin{eqnarray*}
\log(1+ y)=\left(\exp\left(\sum_{j\in \mathbb{Z}_{+}}
B_{j}y^{j+1}\frac{\partial}{\partial y}\right)\right)y.
\end{eqnarray*}
Since the composition inverse of $\log(1+y)$
is $e^{y}-1$, we have
\begin{eqnarray*}
e^{y}-1=\left(\exp\left(-\sum_{j\in \mathbb{Z}_{+}}
B_{j}y^{j+1}\frac{\partial}{\partial y}\right)\right)y.
\end{eqnarray*}
By the change-of-variable formula  in \cite{H3.5}, we have
\begin{equation}\label{chg-var}
Y(e^{yL(0)}u, e^{y}-1)=e^{-L_{+}(B)}
Y(e^{L_{+}(B)}u, y)e^{L_{+}(B)}.
\end{equation}
By the definition of the sequences $\{A_{j}\}_{j\in \mathbb{Z}_{+}}$
and $\{B_{j}\}_{j\in \mathbb{Z}_{+}}$,
we have
\begin{eqnarray*}
\lefteqn{\left(\exp\left(\sum_{j\in \mathbb{Z}_{+}}
B_{j}y^{j+1}\frac{\partial}{\partial y}\right)\right)y}\nn
&&=(2\pi i)^{\displaystyle -y\frac{\partial}{\partial y}}
\left(\exp\left(\sum_{j\in \mathbb{Z}_{+}}
A_{j}y^{j+1}\frac{\partial}{\partial y}\right)\right)
(2\pi i)^{\displaystyle y\frac{\partial}{\partial y}}y.
\end{eqnarray*}

Thus
\begin{eqnarray}
e^{-L_{+}(B)}&=&(2\pi i)^{L(0)}e^{-L_{+}(A)}(2\pi i)^{-L(0)}\nn
&=&\mathcal{U}(1)(2\pi i)^{-L(0)},\label{B-A-1}\\
e^{L_{+}(B)}&=&(2\pi i)^{L(0)}e^{L_{+}(A)}(2\pi i)^{-L(0)}\nn
&=&(2\pi i)^{L(0)}(\mathcal{U}(1))^{-1}.\label{B-A-2}
\end{eqnarray}
Substituting $\mathcal{U}(1)w$ for $w$ in (\ref{chg-var})
and using (\ref{ux-u1}), (\ref{B-A-1}), (\ref{B-A-2})
and the $L(0)$-conjugation formula
\begin{equation}\label{l-0-conj}
(2\pi i)^{-L(0)}\mathcal{Y}((2\pi i)^{L(0)}w, x)(2\pi i)^{L(0)}
=\mathcal{Y}\left(w, \frac{x}{2\pi i}\right),
\end{equation}
we obtain 
$$\mathcal{U}(1)\mathcal{Y}\left(w, \frac{x}{2\pi i}\right)
(\mathcal{U}(1))^{-1}
=\mathcal{Y}(\mathcal{U}(e^{x})w, e^{x}-1),$$
which is equivalent to (\ref{chg-var0}).

As a special case of the Jacobi identity 
\cite{FHL} defining intertwining operators, we
have the commutator formula between vertex operators and  intertwining
operators
$$[Y(u, x_{1}), \mathcal{Y}(w, x_{2})]=
\res_{x_{0}}x_{2}^{-1}\delta\left(\frac{x_{1}-x_{0}}{x_{2}}\right)
\mathcal{Y}(Y(u, x_{0})w, x_{2}).$$
Thus
\begin{eqnarray}\label{x-comm1}
\lefteqn{[Y(x_{1}^{L(0)}u, x_{1}), \mathcal{Y}(x_{2}^{L(0)}w,
x_{2})]}\nn
&&=\res_{x_{0}}x_{2}^{-1}\delta\left(\frac{x_{1}-x_{0}}{x_{2}}\right)
\mathcal{Y}(Y(x_{1}^{L(0)}u, x_{0})x_{2}^{L(0)}w, x_{2})\nn
&&=\res_{x_{0}}x_{2}^{-1}\delta\left(\frac{x_{1}-x_{0}}{x_{2}}\right)
\mathcal{Y}\left(x_{2}^{L(0)}Y\left(\left(\frac{x_{1}}{x_{2}}\right)^{L(0)}u, 
\frac{x_{0}}{x_{2}}\right))w, x_{2}\right).\nn
&&
\end{eqnarray}
Now we change the variable $x_{0}$ to $y$ as follows:
$$x_{0}=\sum_{k\ge 1}\frac{x_{2}(2\pi iy)^{k}}{k!}=(e^{2\pi iy}-1)x_{2}.$$
(Note that by definition, for any $r\in \mathbb{C}$,
$$(x_{2}+(e^{2\pi iy}-1)x_{2})^{r}=e^{2\pi iyr}x_{2}^{r}.)$$
Then the right-hand side of (\ref{x-comm1}) becomes
\begin{eqnarray}\label{x-comm2}
\lefteqn{\res_{y}\frac{\partial x_{0}}{\partial y}
x_{2}^{-1}\delta\left(\frac{x_{1}-(e^{2\pi iy}-1)x_{2}}{x_{2}}\right)
\mathcal{Y}(x_{2}^{L(0)}Y((x_{1}/x_{2})^{L(0)}u, e^{2\pi iy}-1)w, x_{2})}\nn
&&=2\pi i\res_{y}e^{2\pi iy}x_{2}
(x_{2}+(e^{2\pi iy}-1)x_{2})^{-1}
\delta\left(\frac{x_{1}}{x_{2}+(e^{2\pi iy}-1)x_{2}}\right)\cdot\nn
&&\quad\quad\quad\quad\quad\quad\quad\quad\quad\quad\quad\quad\quad\cdot
\mathcal{Y}(x_{2}^{L(0)}Y((x_{1}/x_{2})^{L(0)}u, e^{2\pi iy}-1)w, x_{2})\nn
&&=2\pi i\res_{y}
\delta\left(\frac{x_{1}}{e^{2\pi iy}x_{2}}\right)
\mathcal{Y}(x_{2}^{L(0)}Y((x_{1}/x_{2})^{L(0)}u, e^{2\pi iy}-1)w, x_{2})\nn
&&=2\pi i\res_{y}\delta\left(\frac{x_{1}}{e^{2\pi iy}x_{2}}\right)
\mathcal{Y}(x_{2}^{L(0)}Y(e^{2\pi iyL(0)}u, e^{2\pi iy}-1)w, x_{2})\nn
&&=2\pi i\res_{y}\delta\left(\frac{x_{1}}{e^{2\pi iy}x_{2}}\right)
\mathcal{Y}(x_{2}^{L(0)}
\mathcal{U}(1)Y((\mathcal{U}(1))^{-1}u, y)
(\mathcal{U}(1))^{-1}w, x_{2}),\nn
&&
\end{eqnarray}
where in the last step, we have used (\ref{chg-var0}).
From (\ref{x-comm1}) and
(\ref{x-comm2}), we obtain 
\begin{eqnarray}\label{x-comm3}
\lefteqn{[Y(x_{1}^{L(0)}u, x_{1}), \mathcal{Y}(x_{2}^{L(0)}w,
x_{2})]}\nn
&&=2\pi i\res_{y}\delta\left(\frac{x_{1}}{e^{2\pi iy}x_{2}}\right)
\mathcal{Y}\left(\mathcal{U}(x_{2})Y((\mathcal{U}(1))^{-1}u, y)
(\mathcal{U}(1))^{-1}w, x_{2}\right). \nn
&&
\end{eqnarray} 
Substituting $\mathcal{U}(1)u$ for $u$ and 
$\mathcal{U}(1)w$ for $w$ in (\ref{x-comm3})
and using (\ref{ux-u1}) and (\ref{l-0-conj}), we obtain
(\ref{x-comm0}).
\epfv

The operators $Y(\mathcal{U}(x_{1})u, x_{1})$ and
$\mathcal{Y}(\mathcal{U}(x_{2})w, x_{2})$ are called {\it 
geometrically-modified vertex operator} and {\it 
geometrically-modified intertwining operator}, respectively. 
These operators
play an important role in this paper. See Remark \ref{geometry} in Section 3
for the reason why we need these modified operators. 

\begin{prop}\label{l-1-e(y)}
Let $W_{1}, W_{2}$ and $W_{3}$ be $V$-modules and  $\mathcal{Y}$ 
an intertwining operator of type 
$\binom{W_{3}}{W_{1}W_{2}}$.  Then for any $w_{1}\in W_{1}$,
$$2\pi ix\frac{d}{dx}\mathcal{Y}(\mathcal{U}(x)w_{1}, 
x)=\mathcal{Y}(\mathcal{U}(x)L(-1)w_{1}, x).$$
\end{prop}
\pf
Using the $L(-1)$-derivative property for $\mathcal{Y}$,
we obtain
\begin{eqnarray}\label{der-e(y)}
\lefteqn{2\pi ix\frac{d}{dx}\mathcal{Y}(\mathcal{U}(x)w_{1}, 
x)}\nn
&&=2\pi ix\frac{d}{dx}\mathcal{Y}((2\pi ix)^{L(0)}
e^{-L_{+}(A)}w_{1}, 
x)\nn
&&=\mathcal{Y}((2\pi ixL(-1)+2\pi iL(0))
(2\pi ix)^{L(0)}e^{-L_{+}(A)}w_{1}, 
x)\nn
&&=\mathcal{Y}(x^{L(0)}
(2\pi iL(-1)+2\pi iL(0))\mathcal{U}(1)w_{1}, 
x).
\end{eqnarray}

Using (\ref{chg-var0}) for the vertex operator $Y$
defining the module $W_{1}$ and the conformal element $\omega$, we obtain
$$\mathcal{U}(1)Y(\omega, x)
=Y(\mathcal{U}(e^{2\pi ix})\omega, e^{2\pi ix}-1)\mathcal{U}(1).$$
This formula together with Lemma \ref{u1-omega}, the identity property
and by changing the variable from $x$ to $y=e^{2\pi ix}-1$ gives
\begin{eqnarray}\label{u-1-l-1-comm}
\mathcal{U}(1)L(-1)
&=&\res_{x}Y(\mathcal{U}(e^{2\pi ix})\omega, e^{2\pi ix}-1)\mathcal{U}(1)\nn
&=&\res_{x}Y(e^{2\pi ixL(0)}\mathcal{U}(1)\omega, 
e^{2\pi ix}-1)\mathcal{U}(1)\nn
&=&\res_{x}(2\pi i)^{2}
Y\left(e^{2\pi ixL(0)}\left(\omega-\frac{c}{24}\mathbf{1}\right), 
e^{2\pi ix}-1\right)\mathcal{U}(1)\nn
&=&\res_{x}(2\pi i)^{2}Y(e^{2\pi ixL(0)}\omega, 
e^{2\pi ix}-1)\mathcal{U}(1)\nn
&=&\res_{x}(2\pi i)^{2}e^{4\pi ix}Y(\omega, 
e^{2\pi ix}-1)\mathcal{U}(1)\nn
&=&\res_{y}(2\pi i+2\pi iy)Y(\omega, 
y)\mathcal{U}(1)\nn
&=&(2\pi iL(-1)+2\pi iL(0))\mathcal{U}(1).
\end{eqnarray}
Using (\ref{u-1-l-1-comm}), we see that the right-hand side of 
(\ref{der-e(y)}) is equal to 
$$\mathcal{Y}(\mathcal{U}(x)L(-1)w_{1}, 
x),$$
proving the proposition. 
\epfv

\begin{prop}\label{e-assoc}
Let $W_{i}$, $i=0, \dots, 5$, be $V$-modules and $\mathcal{Y}_{1}$, 
$\mathcal{Y}_{2}$, $\mathcal{Y}_{3}$ and 
$\mathcal{Y}_{4}$ intertwining operators of types 
$\binom{W_{0}}{W_{1}W_{4}}$, 
$\binom{W_{4}}{W_{2}W_{3}}$, $\binom{W_{5}}{W_{1}W_{2}}$
and $\binom{W_{0}}{W_{5}W_{3}}$, respectively, such that 
when $|z_{1}|>|z_{2}|>|z_{1}-z_{2}|>0$, we have the associativity
$$\langle w'_{0}, \mathcal{Y}_{1}(w_{1}, z_{1})
\mathcal{Y}_{2}(w_{2}, z_{2})w_{3}\rangle
=\langle w'_{0}, \mathcal{Y}_{4}(\mathcal{Y}_{3}(w_{1}, 
z_{1})w_{2}, z_{2})w_{3}\rangle$$
for $w_{i}\in W_{i}$, $i=1, 2, 3$ and $w'_{0}\in W'_{0}$.
Then when $|q_{(z_{1}-z_{2})}|>1>|q_{(z_{1}-z_{2})}-
1|>0$,
\begin{eqnarray*}
\lefteqn{\langle w'_{0}, 
\mathcal{Y}_{1}(\mathcal{U}(q_{z_{1}})w_{1}, q_{z_{1}})
\mathcal{Y}_{2}(\mathcal{U}(q_{z_{2}})w_{2}, q_{z_{2}})w_{3}\rangle}\nn
&&=
\langle w'_{0}, \mathcal{Y}_{4}(\mathcal{U}(q_{z_{2}})
\mathcal{Y}_{3}(w_{1}, z_{1}-z_{2})w_{2}, q_{z_{2}})w_{3}\rangle.
\end{eqnarray*}
for $w_{i}\in W_{i}$, $i=1, 2, 3$ and $w'_{0}\in W'_{0}$.
\end{prop}
\pf
By associativity, (\ref{ux-u1}), the $L(0)$-conjugation 
property and (\ref{chg-var0}), 
\begin{eqnarray}\label{e-to-y}
\lefteqn{\langle w'_{0}, 
\mathcal{Y}_{1}(\mathcal{U}(q_{z_{1}})w_{1}, q_{z_{1}})
\mathcal{Y}_{2}(\mathcal{U}(q_{z_{2}})w_{2}, q_{z_{2}})w_{3}\rangle}\nn
&&=\langle w'_{0}, \mathcal{Y}_{4}(\mathcal{Y}_{3}(
\mathcal{U}(q_{z_{1}})w_{1}, 
q_{z_{1}}-q_{z_{2}})
\mathcal{U}(q_{z_{2}})w_{2}, 
q_{z_{2}})w_{3}\rangle\nn
&&=\langle w'_{0}, \mathcal{Y}_{4}(\mathcal{Y}_{3}(
\mathcal{U}(q_{z_{1}})w_{1}, 
q_{z_{1}}-q_{z_{2}})
q_{z_{2}}^{L(0)}\mathcal{U}(1)w_{2}, 
q_{z_{2}})w_{3}\rangle\nn
&&=\langle w'_{0}, \mathcal{Y}_{4}(q_{z_{2}}^{L(0)}\mathcal{Y}_{3}(
\mathcal{U}(q_{z_{1}-z_{2}})w_{1}, 
q_{z_{1}-z_{2}}-1)
\mathcal{U}(1)w_{2}, 
q_{z_{2}})w_{3}\rangle\nn
&&=\langle w'_{0}, \mathcal{Y}_{4}(q_{z_{2}}^{L(0)}\mathcal{U}(1)
\mathcal{Y}_{3}(w_{1}, 
z_{1}-z_{2})
w_{2}, 
q_{z_{2}})w_{3}\rangle\nn
&&=\langle w'_{0}, \mathcal{Y}_{4}(\mathcal{U}(q_{z_{2}})
\mathcal{Y}_{3}(w_{1}, 
z_{1}-z_{2})
w_{2}, 
q_{z_{2}})w_{3}\rangle.
\end{eqnarray}
Since the associativity used above holds in the region
$|q_{z_{1}}|>|q_{z_{2}}|>|q_{z_{1}}-q_{z_{2}}|>0$
or equivalently $|q_{(z_{1}-z_{2})}|>1>|q_{(z_{1}-z_{2})}-
1|>0$, we see that (\ref{e-to-y}) also holds 
in the region $|q_{(z_{1}-z_{2})}|>1>|q_{(z_{1}-z_{2})}-
1|>0$.
\epfv

\renewcommand{\theequation}{\thesection.\arabic{equation}}
\renewcommand{\thethm}{\thesection.\arabic{thm}}
\setcounter{equation}{0}
\setcounter{thm}{0}

\section{Identities for formal $q$-traces}

In this section, we prove several identities on which the main results
in the present paper are based. We shall use both the formal variable 
and the complex variable approaches. The formal variable approach 
is needed because 
the $q$-traces are still formal in this section 
and the complex variable approach is needed because 
we need to use analytic extensions to obtained our results.

For $m\ge 0$, let 
$$P_{m+1}(x; q)=(2\pi i)^{m+1}
\sum_{l> 0}\left(\frac{l^{m}}{m!}\frac{x^{l}}{1-q^{l}}
-\frac{(-1)^{m}l^{m}}{m!}\frac{q^{l}x^{-l}}{1-q^{l}}\right)$$
where $(1-q^{l})^{-1}$ for $l>0$ is understood as the 
series $\sum_{k\ge 0}q^{kl}$. 
For $\tau, z\in \mathbb{C}$ satisfying 
$|q_{\tau}|<|q_{z}|<1$, the series $P_{m+1}(q_{z}; q_{\tau})$ for 
$m\ge 0$ are absolutely convergent and, for $z\in \mathbb{C}$
satisfying $|q_{z}|<1$, the $q$-coefficients of
$P_{m+1}(q_{z}; q)$ is absolutely convergent. 

As we mentioned in Section 1, we shall fix a vertex operator algebra
$V$ throughout this paper. Let $W$ be a $V$-module.
and let $o(u): W\to W$ be the linear map defined by 
$o(u)=u_{\swt u-1}$ for homogeneous $u\in V$. We have:

\begin{prop}\label{identity0-pr}
Let  $W_{i}$ and $\tilde{W}_{i}$, $i=1, \dots, n$, $V$-modules, and
$\mathcal{Y}_{i}$, $i=1, \dots, n$, intertwining operators of 
types $\binom{\tilde{W}_{i-1}}{ W_{i}\tilde{W}_{i}}$, respectively,
where we use the convention $\tilde{W}_{0}=
\tilde{W}_{n}$. 
Then for any $u\in V$, $w_{i}\in W_{i}$, $i=1, \dots, n$, we have 
\begin{eqnarray}\label{identity0}
\lefteqn{\tr_{\tilde{W}_{n}}Y(\mathcal{U}(x)u, x)
\mathcal{Y}_{1}(\mathcal{U}(x_{1})w_{1}, x_{1})
\cdots \mathcal{Y}_{n}(\mathcal{U}(x_{n})w_{n}, x_{n})q^{L(0)}}\nn
&&=\sum_{i=1}^{n}\sum_{m\ge 0}
P_{m+1}\left(\frac{x_{i}}{x}; q\right)\tr_{\tilde{W}_{n}}
\mathcal{Y}_{1}(\mathcal{U}(x_{1})w_{1}, x_{1})\cdot\nn 
&&\quad\quad\quad\quad\quad
\cdots\mathcal{Y}_{i-1}(\mathcal{U}(x_{i-1})w_{i-1}, x_{i-1})
\mathcal{Y}_{i}(\mathcal{U}(x_{i})
u_{m}w_{i}, x_{i})\cdot\nn
&&\quad\quad\quad\quad\quad\cdot
\mathcal{Y}_{i+1}(\mathcal{U}(x_{i+1})w_{i+1}, x_{i+1})
\cdots \mathcal{Y}_{n}(\mathcal{U}(x_{n})w_{n}, x_{n})q^{L(0)}\nn
&&+\tr_{\tilde{W}_{n}}o(\mathcal{U}(1)u)
\mathcal{Y}_{1}(\mathcal{U}(x_{1})w_{1}, x_{1})
\cdots \mathcal{Y}_{n}(\mathcal{U}(x_{n})w_{n}, x_{n})q^{L(0)}
\end{eqnarray}
and 
\begin{eqnarray}\label{identity0.5}
\lefteqn{\sum_{i=1}^{n}
\tr_{\tilde{W}_{n}}
\mathcal{Y}_{1}(\mathcal{U}(x_{1})w_{1}, x_{1})
\cdots \mathcal{Y}_{i-1}(\mathcal{U}(x_{i-1})w_{i-1}, x_{i-1})\cdot}\nn
&&\cdot 
\mathcal{Y}_{i}(\mathcal{U}(x_{i})u_{0}w_{i}, x_{i})
\mathcal{Y}_{i+1}(\mathcal{U}(x_{i+1})w_{i+1}, x_{i+1})
\cdots \mathcal{Y}_{n}(\mathcal{U}(x_{n})w_{n}, x_{n})q^{L(0)}\nn
&&=0.
\end{eqnarray}
\end{prop}
\pf
The left-hand side of (\ref{identity0}) is equal to
\begin{eqnarray}\label{identity0-1}
\lefteqn{\sum_{i=1}^{n}\tr_{\tilde{W}_{n}}
\mathcal{Y}_{1}(\mathcal{U}(x_{1})w_{1}, x_{1})
\cdots \mathcal{Y}_{i-1}(\mathcal{U}(x_{i-1})w_{i-1}, x_{i-1})\cdot}\nn
&&\quad\quad\quad\cdot [Y(\mathcal{U}(x)u, x), 
\mathcal{Y}_{i}(\mathcal{U}(x_{i})w_{i}, x_{i})]
\mathcal{Y}_{i+1}(\mathcal{U}(x_{i+1})w_{i+1}, x_{i+1})\cdot\nn
&&\quad\quad\quad\cdots \mathcal{Y}_{n}(\mathcal{U}(x_{n})w_{n}, x_{n})
q^{L(0)}\nn
&&\quad +\tr_{\tilde{W}_{n}}\mathcal{Y}_{1}(\mathcal{U}(x_{1}w_{1}, x_{1})
\cdots \mathcal{Y}_{n}(\mathcal{U}(x_{n})w_{n}, x_{n})Y(\mathcal{U}(x)u, x)
q^{L(0)}\nn
&&=\sum_{i=1}^{n}2\pi i\res_{y}
\delta\left(\frac{x}{e^{2\pi iy}x_{i}}\right)\tr_{\tilde{W}_{n}}
\mathcal{Y}_{1}(\mathcal{U}(x_{1})w_{1}, x_{1})
\cdot\nn
&&\quad\quad\quad\cdot
\mathcal{Y}_{i-1}(\mathcal{U}(x_{i-1})w_{i-1}, x_{i-1})
\mathcal{Y}_{i}(\mathcal{U}(x_{i})Y(u, y)
w_{i}, x_{i})\cdot\nn
&&\quad\quad\quad\cdot
\mathcal{Y}_{i+1}(\mathcal{U}(x_{i+1})w_{i+1}, x_{i+1})
 \mathcal{Y}_{n}(\mathcal{U}(x_{n})w_{n}, x_{n})
q^{L(0)}\nn
&&\quad +\tr_{\tilde{W}_{n}}
Y\left(\mathcal{U}\left(\frac{x}{q}\right)u, \frac{x}{q}\right)
\mathcal{Y}_{1}(\mathcal{U}(x_{1})w_{1}, x_{1})
\cdots \mathcal{Y}_{n}(\mathcal{U}(x_{n})w_{n}, x_{n})
q^{L(0)}\nn
&&=\res_{y}\sum_{i=1}^{n}\sum_{m\ge 0}\frac{(2\pi i)^{m+1}y^{m}}{m!}
\left(\left(x_{i}\frac{\partial}{\partial x_{i}}\right)^{m}
\delta\left(\frac{x}{x_{i}}\right)\right)\cdot\nn
&&\quad\quad\quad\cdot\tr_{\tilde{W}_{n}}
\mathcal{Y}_{1}(\mathcal{U}(x_{1})w_{1}, x_{1})
\cdots \mathcal{Y}_{i-1}(\mathcal{U}(x_{i-1})w_{i-1}, x_{i-1})\cdot\nn
&&\quad\quad\quad\cdot 
\mathcal{Y}_{i}(\mathcal{U}(x_{i})Y(u, y)
w_{i}, x_{i})
\mathcal{Y}_{i+1}(\mathcal{U}(x_{i+1})w_{i+1}, x_{i+1})\cdot\nn
&&\quad\quad\quad\cdots \mathcal{Y}_{n}(\mathcal{U}(x_{n})w_{n}, x_{n})
q^{L(0)}\nn
&&\quad +q^{-x\frac{\partial}{\partial x}}\tr_{\tilde{W}_{n}}
Y(\mathcal{U}(x)u, x)\mathcal{Y}_{1}(\mathcal{U}(x_{1})w_{1}, x_{1})
\cdots \mathcal{Y}_{n}(\mathcal{U}(x_{n})w_{n}, x_{n})
q^{L(0)},\nn
&&
\end{eqnarray}
where we have used (\ref{x-comm0}), the $L(0)$-conjugation property
for vertex operators and a property of traces ($\tr AB=\tr BA$). 

From (\ref{identity0-1}), we obtain
\begin{eqnarray}\label{identity0-2}
\lefteqn{(1-q^{-x\frac{\partial}{\partial x}})
\tr_{\tilde{W}_{n}}Y(\mathcal{U}(x)u, x)
\mathcal{Y}_{1}(\mathcal{U}(x_{1})w_{1}, x_{1})
\cdots \mathcal{Y}_{n}(\mathcal{U}(x_{n})w_{n}, x_{n})
q^{L(0)}}\nn
&&=\res_{y}\sum_{i=1}^{n}\sum_{m\ge 0}\frac{(2\pi i)^{m+1}y^{m}}{m!}
\left(\left(x_{i}\frac{\partial}{\partial x_{i}}\right)^{m}
\delta\left(\frac{x}{x_{i}}\right)\right)\cdot\nn
&&\quad\quad\quad\cdot\tr_{\tilde{W}_{n}}
\mathcal{Y}_{1}(\mathcal{U}(x_{1})w_{1}, x_{1})
\cdots \mathcal{Y}_{i-1}(\mathcal{U}(x_{i-1})w_{i-1}, x_{i-1})\cdot\nn
&&\quad\quad\quad\cdot 
\mathcal{Y}_{i}(\mathcal{U}(x_{i})Y(u, y)
w_{i}, x_{i})
\mathcal{Y}_{i+1}(\mathcal{U}(x_{i+1})w_{i+1}, x_{i+1})\cdot\nn
&&\quad\quad\quad\cdots \mathcal{Y}_{n}(\mathcal{U}(x_{n})w_{n}, x_{n})q^{L(0)}\nn
&&=\sum_{i=1}^{n}\sum_{m\ge 0}\sum_{l>0}\frac{(2\pi i)^{m+1}}{m!}
\left(\left(x_{i}\frac{\partial}{\partial x_{i}}\right)^{m}
\left(\frac{x^{l}}{x_{i}^{l}}\right)\right)\cdot\nn
&&\quad\quad\quad\cdot\tr_{\tilde{W}_{n}}
\mathcal{Y}_{1}(\mathcal{U}(x_{1})w_{1}, x_{1})
\cdots \mathcal{Y}_{i-1}(\mathcal{U}(x_{i-1})w_{i-1}, x_{i-1})\cdot\nn
&&\quad\quad\quad\cdot 
\mathcal{Y}_{i}(\mathcal{U}(x_{i})u_{m}w_{i}, x_{i})
\mathcal{Y}_{i+1}(\mathcal{U}(x_{i+1})w_{i+1}, x_{i+1})\cdot\nn
&&\quad\quad\quad\cdots \mathcal{Y}_{n}(\mathcal{U}(x_{n})w_{n}, x_{n})q^{L(0)}\nn
&&\quad +\sum_{i=1}^{n}\sum_{m\ge 0}\sum_{l>0}\frac{(2\pi i)^{m+1}}{m!}
\left(\left(x_{i}\frac{\partial}{\partial x_{i}}\right)^{m}
\left(\frac{x^{-l}}{x_{i}^{-l}}\right)\right)\cdot\nn
&&\quad\quad\quad\cdot\tr_{\tilde{W}_{n}}
\mathcal{Y}_{1}(\mathcal{U}(x_{1})w_{1}, x_{1})
\cdots \mathcal{Y}_{i-1}(\mathcal{U}(x_{i-1})w_{i-1}, x_{i-1})\cdot\nn
&&\quad\quad\quad\cdot 
\mathcal{Y}_{i}(\mathcal{U}(x_{i})u_{m}w_{i}, x_{i})
\mathcal{Y}_{i+1}(\mathcal{U}(x_{i+1})w_{i+1}, x_{i+1})\cdot\nn
&&\quad\quad\quad\cdots \mathcal{Y}_{n}(\mathcal{U}(x_{n})w_{n}, x_{n})q^{L(0)}\nn
&&\quad +2\pi i\sum_{i=1}^{n}
\tr_{\tilde{W}_{n}}
\mathcal{Y}_{1}(\mathcal{U}(x_{1})w_{1}, x_{1})
\cdots \mathcal{Y}_{i-1}(\mathcal{U}(x_{i-1})w_{i-1}, x_{i-1})\cdot\nn
&&\quad\quad\quad\cdot 
\mathcal{Y}_{i}(\mathcal{U}(x_{i})u_{0}w_{i}, x_{i})
\mathcal{Y}_{i+1}(\mathcal{U}(x_{i+1})w_{i+1}, x_{i+1})\cdot\nn
&&\quad\quad\quad\cdots \mathcal{Y}_{n}(\mathcal{U}(x_{n})w_{n}, x_{n})q^{L(0)}.
\end{eqnarray}

Note that the left-hand side of (\ref{identity0-2}) has no constant term
as a series in $x$. Thus the constant term (the last term) 
of the right-hand side as
a series in $x$ must be $0$. So we obtain (\ref{identity0.5}) 
and consequently the left-hand side of (\ref{identity0-2}) is equal to the 
sum of the first two terms in the right-hand side of (\ref{identity0-2}). 
Since $1-q^{-x\frac{\partial}{\partial x}}$ acting on any term
independent of $x$ is $0$, we obtain
\begin{eqnarray}\label{identity0-3}
\lefteqn{(1-q^{-x\frac{\partial}{\partial x}})
\left(\tr_{\tilde{W}_{n}}Y(\mathcal{U}(x)u, x)
\mathcal{Y}_{1}(\mathcal{U}(x_{1})w_{1}, x_{1})
\cdots \mathcal{Y}_{n}(\mathcal{U}(x_{n})w_{n}, x_{n})
q^{L(0)}\right.}\nn
&&\quad\quad\quad\quad\quad\quad\left.-\tr_{\tilde{W}_{n}}o(\mathcal{U}(1)u)
\mathcal{Y}_{1}(\mathcal{U}(x_{1})w_{1}, x_{1})
\cdots \mathcal{Y}_{n}(\mathcal{U}(x_{n})w_{n}, x_{n})
q^{L(0)}\right)\nn
&&=\sum_{i=1}^{n}\sum_{m\ge 0}\sum_{l>0}\frac{(2\pi i)^{m+1}}{m!}
\left(\left(x_{i}\frac{\partial}{\partial x_{i}}\right)^{m}
\left(\frac{x^{l}}{x_{i}^{l}}\right)\right)\cdot\nn
&&\quad\quad\quad\cdot\tr_{\tilde{W}_{n}}
\mathcal{Y}_{1}(\mathcal{U}(x_{1})w_{1}, x_{1})
\cdots \mathcal{Y}_{i-1}(\mathcal{U}(x_{i-1})w_{i-1}, x_{i-1})\cdot\nn
&&\quad\quad\quad\cdot 
\mathcal{Y}_{i}(\mathcal{U}(x_{i})u_{m}w_{i}, x_{i})
\mathcal{Y}_{i+1}(\mathcal{U}(x_{i+1})w_{i+1}, x_{i+1})\cdot\nn
&&\quad\quad\quad\cdots \mathcal{Y}_{n}(\mathcal{U}(x_{n})w_{n}, x_{n})q^{L(0)}\nn
&&\quad +\sum_{i=1}^{n}\sum_{m\ge 0}\sum_{l>0}\frac{(2\pi i)^{m+1}}{m!}
\left(\left(x_{i}\frac{\partial}{\partial x_{i}}\right)^{m}
\left(\frac{x^{-l}}{x_{i}^{-l}}\right)\right)\cdot\nn
&&\quad\quad\quad\cdot\tr_{\tilde{W}_{n}}
\mathcal{Y}_{1}(\mathcal{U}(x_{1})w_{1}, x_{1})
\cdots \mathcal{Y}_{i-1}(\mathcal{U}(x_{i-1})w_{i-1}, x_{i-1})\cdot\nn
&&\quad\quad\quad\cdot \mathcal{Y}_{i}(\mathcal{U}(x_{i})u_{m}w_{i}, x_{i})
\mathcal{Y}_{i+1}(\mathcal{U}(x_{i+1})w_{i+1}, x_{i+1})\cdot\nn
&&\quad\quad\quad\cdots \mathcal{Y}_{n}(\mathcal{U}(x_{n})w_{n}, x_{n})q^{L(0)}\nn
&&=\sum_{i=1}^{n}\sum_{m\ge 0}\sum_{l>0}\left(\frac{(2\pi i)^{m+1}}{m!}
\left(\left(x_{i}\frac{\partial}{\partial x_{i}}\right)^{m}
\left(\frac{x^{l}}{x_{i}^{l}}\right)\right)\right.\nn
&&\left.\quad\quad\quad\quad\quad\quad\quad\quad\quad\quad\quad\quad
+\frac{(2\pi i)^{m+1}}{m!}
\left(\left(x_{i}\frac{\partial}{\partial x_{i}}\right)^{m}
\left(\frac{x^{-l}}{x_{i}^{-l}}\right)\right)\right)\cdot\nn
&&\quad\quad\quad\cdot\tr_{\tilde{W}_{n}}
\mathcal{Y}_{1}(\mathcal{U}(x_{1})w_{1}, x_{1})
\cdots \mathcal{Y}_{i-1}(\mathcal{U}(x_{i-1})w_{i-1}, x_{i-1})\cdot\nn
&&\quad\quad\quad\cdot \mathcal{Y}_{i}(\mathcal{U}(x_{i})u_{m}w_{i}, x_{i})
\mathcal{Y}_{i+1}(\mathcal{U}(x_{i+1})w_{i+1}, x_{i+1})\cdot\nn
&&\quad\quad\quad\cdots \mathcal{Y}_{n}(\mathcal{U}(x_{n})w_{n}, x_{n})q^{L(0)}.
\end{eqnarray}

Note that when acting on nonzero power of $x$, 
$1-q^{-x\frac{\partial}{\partial x}}$ is invertible. But on different 
rings, it has different inverses. 
Since all the terms in our formulas as series in $q$ are lower 
truncated, $(1-q^{-x\frac{\partial}{\partial x}})^{-1}$ acting 
on any nonzero power of $x$ should also be lower truncated 
as a series in $q$. 
Thus we have 
\begin{equation}\label{identity0-4}
(1-q^{-x\frac{\partial}{\partial x}})^{-1}x^{l}
=\left\{\begin{array}{ll}
{\displaystyle \frac{-q^{l}x^{l}}{1-q^{l}}}&l>0\\
{\displaystyle \frac{x^{l}}{1-q^{l}}}&l<0.
\end{array}\right.
\end{equation}
By (\ref{identity0-3}) and (\ref{identity0-4}), we obtain
\begin{eqnarray}\label{identity0-5}
\lefteqn{\tr_{\tilde{W}_{n}}Y(\mathcal{U}(x)u, x)
\mathcal{Y}_{1}(\mathcal{U}(x_{1})w_{1}, x_{1})
\cdots \mathcal{Y}_{n}(\mathcal{U}(x_{n})w_{n}, x_{n})
q^{L(0)}}\nn
&&\quad-\tr_{\tilde{W}_{n}}o(\mathcal{U}(1)u)
\mathcal{Y}_{1}(\mathcal{U}(x_{1})w_{1}, x_{1})
\cdots \mathcal{Y}_{n}(\mathcal{U}(x_{n})w_{n}, x_{n})
q^{L(0)}\nn
&&=(1-q^{-x\frac{\partial}{\partial x}})^{-1}
\sum_{i=1}^{n}\sum_{m\ge 0}\sum_{l>0}\left(\frac{(2\pi i)^{m+1}}{m!}
\left(\left(x_{i}\frac{\partial}{\partial x_{i}}\right)^{m}
\left(\frac{x^{l}}{x_{i}^{l}}\right)\right)\right.\nn
&&\left.\quad\quad\quad\quad\quad\quad\quad\quad\quad\quad\quad
\quad\quad\quad\quad
+\frac{(2\pi i)^{m+1}}{m!}
\left(\left(x_{i}\frac{\partial}{\partial x_{i}}\right)^{m}
\left(\frac{x^{-l}}{x_{i}^{-l}}\right)\right)\right)\cdot\nn
&&\quad\quad\quad\cdot\tr_{\tilde{W}_{n}}
\mathcal{Y}_{1}(\mathcal{U}(x_{1})w_{1}, x_{1})
\cdots \mathcal{Y}_{i-1}(\mathcal{U}(x_{i-1})w_{i-1}, x_{i-1})\cdot\nn
&&\quad\quad\quad\cdot 
\mathcal{Y}_{i}(\mathcal{U}(x_{i})u_{m}w_{i}, x_{i})
\mathcal{Y}_{i+1}(\mathcal{U}(x_{i+1})w_{i+1}, x_{i+1})\cdot\nn
&&\quad\quad\quad\cdots \mathcal{Y}_{n}(\mathcal{U}(x_{n})w_{n}, x_{n})q^{L(0)}\nn
&&=\sum_{i=1}^{n}\sum_{m\ge 0}\sum_{l>0}\left(\frac{(2\pi i)^{m+1}}{m!}
\left(\left(x_{i}\frac{\partial}{\partial x_{i}}\right)^{m}
\left(\frac{-q^{l}x^{l}}{(1-q^{l})x_{i}^{l}}\right)\right)\right.\nn
&&\left.\quad\quad\quad\quad\quad\quad\quad\quad\quad\quad\quad
+\frac{(2\pi i)^{m+1}}{m!}
\left(\left(x_{i}\frac{\partial}{\partial x_{i}}\right)^{m}
\left(\frac{x^{-l}}{(1-q^{l})x_{i}^{-l}}\right)\right)\right)\cdot\nn
&&\quad\quad\quad\cdot\tr_{\tilde{W}_{n}}
\mathcal{Y}_{1}(\mathcal{U}(x_{1})w_{1}, x_{1})
\cdots \mathcal{Y}_{i-1}(\mathcal{U}(x_{i-1})w_{i-1}, x_{i-1})\cdot\nn
&&\quad\quad\quad\cdot \mathcal{Y}_{i}(\mathcal{U}(x_{i})u_{m}w_{i}, x_{i})
\mathcal{Y}_{i+1}(\mathcal{U}(x_{i+1})w_{i+1}, x_{i+1})\cdot\nn
&&\quad\quad\quad\cdots \mathcal{Y}_{n}(\mathcal{U}(x_{n})w_{n}, x_{n})q^{L(0)}\nn
&&=\sum_{i=1}^{n}\sum_{m\ge 0}\sum_{l>0}
\left(-\frac{(-1)^{m}(2\pi i)^{m+1}l^{m}}{m!}
\left(\frac{q^{l}\left(\frac{x_{i}}{x}\right)^{-l}}{1-q^{l}}\right)\right.\nn
&&\left.\quad\quad\quad\quad\quad\quad\quad\quad\quad\quad\quad\quad
\quad\quad\quad\quad+\frac{(2\pi i)^{m+1}l^{m}}{m!}
\left(\frac{\left(\frac{x_{i}}{x}\right)^{l}}{1-q^{l}}
\right)\right)\cdot\nn
&&\quad\quad\quad\cdot\tr_{\tilde{W}_{n}}
\mathcal{Y}_{1}(\mathcal{U}(x_{1})w_{1}, x_{1})
\cdots \mathcal{Y}_{i-1}(\mathcal{U}(x_{i-1})w_{i-1}, x_{i-1})\cdot\nn
&&\quad\quad\quad\cdot \mathcal{Y}_{i}(\mathcal{U}(x_{i})u_{m}w_{i}, x_{i})
\mathcal{Y}_{i+1}(\mathcal{U}(x_{i+1})w_{i+1}, x_{i+1})\cdot\nn
&&\quad\quad\quad\cdots \mathcal{Y}_{n}(\mathcal{U}(x_{n})w_{n}, x_{n})q^{L(0)}\nn
&&=\sum_{i=1}^{n}\sum_{m\ge 0}
P_{m+1}\left(\frac{x_{i}}{x}; q\right)\tr_{\tilde{W}_{n}}
\mathcal{Y}_{1}(\mathcal{U}(x_{1})w_{1}, x_{1})\cdot\nn 
&&\quad\quad\quad
\cdots\mathcal{Y}_{i-1}(\mathcal{U}(x_{i-1})w_{i-1}, x_{i-1})
\mathcal{Y}_{i}(\mathcal{U}(x_{i})u_{m}w_{i}, x_{i})\cdot\nn
&&\quad\quad\quad\cdot
\mathcal{Y}_{i+1}(\mathcal{U}(x_{i+1})w_{i+1}, x_{i+1})
\cdots \mathcal{Y}_{n}(\mathcal{U}(x_{n})w_{n}, x_{n})q^{L(0)}.
\end{eqnarray}
The identity (\ref{identity0-5}) is equivalent to (\ref{identity0}).
\epfv

The identities above are  proved using 
purely formal variable approach. But starting from the next 
identity, we need to use both formal variables and 
complex variables.

Let 
\begin{eqnarray*}
\wp_{1}(z; \tau)&=&\frac{1}{z}+\sum_{(k, l)\ne (0, 0)}
\left(\frac{1}{z-(k\tau+l)}+\frac{1}{k\tau+l}
+\frac{z}{(k\tau+l)^{2}}\right),\\
\wp_{2}(z; \tau)&=&\frac{1}{z^{2}}+\sum_{(k, l)\ne (0, 0)}
\left(\frac{1}{(z-(k\tau+l))^{2}}-\frac{1}{(k\tau+l)^{2}}\right)\nn
&=&-\frac{\partial}{\partial z}\wp_{1}(z; \tau)
\end{eqnarray*}
be the Weierstrass zeta function and the 
Weierstrass $\wp$-function, respectively, and let $\wp_{m}(z;\tau)$ for $m>2$ 
be the elliptic functions defined recursively by
\begin{equation}\label{wp-relation}
\wp_{m+1}(z, \tau)=-\frac{1}{m}\frac{\partial}{\partial z} \wp_{m}(z;
\tau).
\end{equation}
These functions have the following 
Laurent expansions in the region $0<|z|<\min(1, |\tau|)$: For $m\ge 1$
\begin{equation}\label{wp0}
\wp_{m}(z; \tau)=\frac{1}{z^{m}}+(-1)^{m}\sum_{k\ge 1}
\binom{2k+1}{m-1}G_{2k+2}(\tau)z^{2k+2-m}
\end{equation}
where 
$$G_{2k+2}(\tau)=\sum_{(m, l)\ne (0, 0)}
\frac{1}{(m\tau+l)^{2k+2}}$$ 
for $k\ge 1$ are the Eisenstein series.
We also have another Eisenstein
series
$$G_{2}(\tau)=\frac{\pi^{2}}{3}
+\sum_{m\in \mathbb{Z}\setminus \{0\}}\sum_{l\in \mathbb{Z}}
\frac{1}{(m\tau+l)^{2}}.$$ 

The Eisenstein series have the following $q$-expansions:
$$G_{2k+2}(\tau)=2\zeta(2k+2)+\frac{2(2\pi i)^{2k+2}}{(2k+1)!}
\sum_{l\ge 1}\frac{l^{2k+1}q_{\tau}^{l}}{1-q_{\tau}^{l}}$$
where 
$$\zeta(2k+2)=\sum_{l\ge 1}\frac{1}{l^{2k+2}}$$
and $(1-q_{\tau}^{l})^{-1}$ is understood as 
$\sum_{j\ge 0}q_{\tau}^{jl}$. 
For $k\ge 1$, let 
$$\tilde{G}_{2k+2}(q)=2\zeta(2k+2)+\frac{2(2\pi i)^{2k+2}}{(2k+1)!}
\sum_{l\ge 1}\frac{l^{2k+1}q^{l}}{1-q^{l}}$$
where $(1-q^{l})^{-1}$ is the formal power series
$\sum_{j\ge 0}q^{jl}$. 
Then $G_{2k+2}(\tau)=\tilde{G}_{2k+2}(q_{\tau})$.

The functions $\wp_{m}(z; \tau)$ for $m\ge 1$
also have the following $q$-expansions: For $z\in \mathbb{C}$ satisfying 
$0<|z|<\min(1, |\tau|)$ and $|q_{\tau}|<|q_{z}|<\frac{1}{|q_{\tau}|}$
and for $m\ge 1$,
\begin{equation}\label{wp.5}
\wp_{m}(z; \tau)=(-1)^{m}\left(P_{m}(q_{z}; q_{\tau})
-\frac{\partial^{m-1}}{\partial z^{m-1}}
(\tilde{G}_{2}(q_{\tau})z-\pi i)\right).
\end{equation}
For $m\ge 1$, let 
\begin{equation}\label{wp1}
\tilde{\wp}_{m}(y; q)=\frac{1}{y^{m}}+(-1)^{m}\sum_{k\ge 1}
\binom{2k+1}{m-1}\tilde{G}_{2k+2}(q)y^{2k+2-m}
\in y^{-m}\mathbb{C}[[y, q]].
\end{equation}
Then for $z\in \mathbb{C}$ satisfying 
$0<|z|<1$ and $m\ge 1$,
\begin{equation}\label{wp2}
\tilde{\wp}_{m}(z; q)=(-1)^{m}\left(P_{m}(q_{z}; q)
-\frac{\partial^{m-1}}{\partial z^{m-1}}
(\tilde{G}_{2}(q)z-\pi i)\right).
\end{equation}
We also have 
$$\tilde{\wp}_{m}(z; q_{\tau})=\wp_{m}(z; \tau)$$
when $m\ge 1$ and $0<|z|<\min(1, |\tau|)$ and 
$|q_{\tau}|<|q_{z}|<\frac{1}{|q_{\tau}|}$.
For more details on the Weierstrass zeta function, the 
Weierstrass $\wp$-function and the Eisenstein series,
see \cite{La} and \cite{Ko}.

Now we assume
that for any $W_{i}$ and $\tilde{W}_{i}$,  
$\mathcal{Y}_{i}$, $i=1, \dots, n$,
as in Proposition \ref{identity0-pr}, 
any $w_{i}\in W_{i}$, $i=1, \dots, n$, 
$\tilde{w}_{n}\in 
\tilde{W}_{n}$ and $\tilde{w}'_{n}\in  \tilde{W}'_{n}$,
$$\langle \tilde{w}'_{n}, \mathcal{Y}_{1}\left(w_{1}, 
z_{1}\right)\cdots 
\mathcal{Y}_{n}\left(w_{n}, 
z_{n}\right)\tilde{w}_{n}\rangle$$
is absolutely convergent when $|z_{1}|>\cdots>|z_{n}|>0$. We also
assume that commutativity and associativity 
for intertwining operators hold. (For more details on 
convergence, commutativity and associativity for intertwining operators,
see \cite{H1}, \cite{H3}, \cite{H4} and \cite{H9}.)

\begin{thm}
Let  $W_{i}$ and $\tilde{W}_{i}$,  $\mathcal{Y}_{i}$, $i=1, \dots, n$,
be as in Proposition \ref{identity0-pr}.
Then for any $u\in V$, $w_{i}\in W_{i}$, $i=1, \dots, n$ and any 
integer $j$ satisfying $1\le j\le n$, we have 
\begin{eqnarray}\label{identity1}
\lefteqn{\tr_{\tilde{W}_{n}}\mathcal{Y}_{1}(\mathcal{U}(x_{1})w_{1}, x_{1})
\cdots \mathcal{Y}_{j-1}(\mathcal{U}(x_{j-1})w_{j-1}, x_{j-1})\cdot
\mathcal{Y}_{j}(\mathcal{U}(x_{j})Y(u, y)w_{j}, x_{j})}\nn
&&\quad\quad\quad\quad\quad\cdot 
\mathcal{Y}_{j+1}(\mathcal{U}(x_{j+1})w_{j+1}, x_{j+1})
\cdots
\mathcal{Y}_{n}(\mathcal{U}(x_{n})w_{n}, x_{n})q^{L(0)}\nn
&&=\sum_{m\ge 0}(-1)^{m+1}
\left(\tilde{\wp}_{m+1}(-y; q)+\frac{\partial^{m}}{\partial y^{m}}
(\tilde{G}_{2}(q)y+\pi i)\right)
\cdot\nn
&&\quad\quad\quad\quad\quad\cdot\tr_{\tilde{W}_{n}}
\mathcal{Y}_{1}(\mathcal{U}(x_{1})w_{1}, x_{1})\cdots 
\mathcal{Y}_{j-1}(\mathcal{U}(x_{j-1})w_{j-1}, x_{j-1})\cdot\nn
&&\quad\quad\quad\quad\quad\cdot
\mathcal{Y}_{j}(\mathcal{U}(x_{j})u_{m}w_{j}, x_{j})
\mathcal{Y}_{j+1}(\mathcal{U}(x_{j+1})w_{j+1}, x_{j+1})\cdot\nn
&&\quad\quad\quad\quad\quad
\cdots \mathcal{Y}_{n}(\mathcal{U}(x_{n})w_{n}, x_{n})q^{L(0)}\nn
&&\quad+\sum_{i\ne j}^{}\sum_{m\ge 0}
P_{m+1}\left(\frac{x_{i}}{x_{j}e^{2\pi iy}}; q\right)\cdot\nn
&&\quad\quad\quad\quad\quad\cdot\tr_{\tilde{W}_{n}}
\mathcal{Y}_{1}(\mathcal{U}(x_{1})w_{1}, x_{1})\cdots 
\mathcal{Y}_{i-1}(\mathcal{U}(x_{i-1})w_{i-1}, x_{i-1})\cdot\nn
&&\quad\quad\quad\quad\quad\cdot
\mathcal{Y}_{i}(\mathcal{U}(x_{i})u_{m}w_{i}, x_{i})
\mathcal{Y}_{i+1}(w_{i+1}, x_{i+1})
\cdots \mathcal{Y}_{n}(\mathcal{U}(x_{n})w_{n}, x_{n})q^{L(0)}\nn
&&\quad+\tr_{\tilde{W}_{n}}o(\mathcal{U}(1)u)
\mathcal{Y}_{1}(\mathcal{U}(x_{1})w_{1}, x_{1})
\cdots \mathcal{Y}_{n}(\mathcal{U}(x_{n})w_{n}, x_{n})q^{L(0)}.\nn
&&
\end{eqnarray}
\end{thm}
\pf
We first prove the case $j=1$.
Note that by  (\ref{chg-var0}), 
\begin{eqnarray}\label{identity1-1}
\mathcal{U}(x_{1})Y(u, y)
&=&x_{1}^{L(0)}\mathcal{U}(1)Y(u, y)(\mathcal{U}(1))^{-1}\mathcal{U}(1)\nn
&=&x_{1}^{L(0)}Y(\mathcal{U}(e^{2\pi iy})u, e^{2\pi iy}-1)\mathcal{U}(1)\nn
&=&Y(x_{1}^{L(0)}\mathcal{U}(e^{2\pi iy})u, x_{1}e^{2\pi iy}-1)
x_{1}^{L(0)}\mathcal{U}(1)\nn
&=&Y(\mathcal{U}(x_{1}e^{2\pi iy})u, 
x_{1}(e^{2\pi iy}-1))
\mathcal{U}(x_{1})
\end{eqnarray}
(recall (\ref{ux-u1})).

Using (\ref{identity1-1}), we see that 
in the case $j=1$, the left-hand side of (\ref{identity1}) is equal to
\begin{eqnarray}\label{identity1-2}
\lefteqn{\tr_{\tilde{W}_{n}}\mathcal{Y}_{1}(\mathcal{U}(x_{1})
Y(u, y)w_{1}, x_{1})
\mathcal{Y}_{2}(\mathcal{U}(x_{2})w_{2}, x_{2})\cdots
\mathcal{Y}_{n}(\mathcal{U}(x_{n})w_{n}, x_{n})q^{L(0)}}\nn
&&=\tr_{\tilde{W}_{n}}\mathcal{Y}_{1}(
Y(\mathcal{U}(x_{1}e^{2\pi iy})u, x_{1}(e^{2\pi iy}-1))
\mathcal{U}(x_{1})w_{1}, x_{1})\cdot\nn
&&\quad\quad\quad\quad\quad\quad\quad\quad\quad\quad\cdot
\mathcal{Y}_{2}(\mathcal{U}(x_{2})w_{2}, x_{2})\cdots
\mathcal{Y}_{n}(\mathcal{U}(x_{n})w_{n}, x_{n})q^{L(0)}.\nn
&&
\end{eqnarray}

By the properties of intertwining operators, we know that the following 
associativity holds: For any $\tilde{w}_{1}\in \tilde{W}_{1}$,
$\tilde{w}'_{n}\in \tilde{W}'_{n}$,
\begin{eqnarray*}
\lefteqn{\langle \tilde{w}'_{n}, \mathcal{Y}_{1}(Y(\mathcal{U}(z_{1}q_{z})u, 
z_{1}(q_{z}-1))
\mathcal{U}(z_{1})w_{1}, z_{1})\tilde{w}_{1}\rangle}\nn
&&=\langle \tilde{w}'_{n}, 
Y(\mathcal{U}(z_{1}q_{z})u, z_{1}q_{z})
\mathcal{Y}_{1}(\mathcal{U}(z_{1})w_{1}, z_{1})\tilde{w}_{1}\rangle
\end{eqnarray*}
for any complex numbers $z$ and $z_{1}$ satisfying 
$|z_{1}q_{z}|>|z_{1}|>|z_{1}(q_{z}-1)|>0$
(or equivalently, $|q_{z}|>1>|q_{z}-1|>0$ and 
$z_{1}\ne 0$), where
we choose the values of $z_{1}^{n}$ for $n\in \mathbb{C}$
to be $e^{n\log z_{1}}$, $0\le \arg z_{1}<2\pi$ (and we shall 
choose this value throughout this paper). 
Thus for any complex numbers $z$ and $z_{1}$ satisfying 
$|q_{z}|>1>|q_{z}-1|>0$ and $z_{1}\ne 0$
and for the values $e^{n\log z_{1}}$ of $z_{1}^{n}$ 
for $n\in \mathbb{C}$, 
\begin{eqnarray}\label{identity1-3}
\lefteqn{\tr_{\tilde{W}_{n}}\mathcal{Y}_{1}(Y(\mathcal{U}(z_{1}q_{z})u, 
z_{1}(q_{z}-1))
\mathcal{U}(z_{1})w_{1}, z_{1})\cdot}\nn
&&\quad\quad\quad\quad\quad\quad\cdot
\mathcal{Y}_{2}(\mathcal{U}(x_{2})w_{2}, x_{2})\cdots
\mathcal{Y}_{n}(\mathcal{U}(x_{n})w_{n}, x_{n})q^{L(0)}\nn
&&=\tr_{\tilde{W}_{n}}Y(\mathcal{U}(z_{1}q_{z})u, z_{1}q_{z})
\mathcal{Y}_{1}(\mathcal{U}(z_{1})w_{1}, z_{1})\cdot\nn
&&\quad\quad\quad\quad\quad\quad\cdot
\mathcal{Y}_{2}(\mathcal{U}(x_{2})w_{2}, x_{2})\cdots
\mathcal{Y}_{n}(\mathcal{U}(x_{n})w_{n}, x_{n})q^{L(0)}
\end{eqnarray}
as formal series in $x_{2}, \dots, x_{n}$ and $q$. Since 
(\ref{identity1-3}) holds for any complex number $z$ satisfying 
$|q_{z}|>1>|q_{z}-1|>0$ and any nonzero complex number $z_{1}$,
the expansions of both sides of (\ref{identity1-3}) as Laurent 
series in $z_{1}$ are also equal for any complex number $z$ satisfying 
$|q_{z}|>1>|q_{z}-1|>0$. Thus
we can replace $z_{1}$ by the formal variable $x_{1}$ 
in (\ref{identity1-3}) so that
for any complex number $z$ satisfying 
$|q_{z}|>1>|q_{z}-1|>0$,
\begin{eqnarray}\label{identity1-4}
\lefteqn{\tr_{\tilde{W}_{n}}\mathcal{Y}_{1}(Y(\mathcal{U}(x_{1}q_{z})u, 
x_{1}(q_{z}-1))
\mathcal{U}(x_{1})w_{1}, x_{1})\cdot}\nn
&&\quad\quad\quad\quad\quad\quad\cdot
\mathcal{Y}_{2}(\mathcal{U}(x_{2})w_{2}, x_{2})\cdots
\mathcal{Y}_{n}(\mathcal{U}(x_{n})w_{n}, x_{n})q^{L(0)}\nn
&&=\tr_{\tilde{W}_{n}}Y(\mathcal{U}(x_{1}q_{z})u, x_{1}q_{z})
\mathcal{Y}_{1}(\mathcal{U}(x_{1})w_{1}, x_{1})\cdot\nn
&&\quad\quad\quad\quad\quad\quad\cdot
\mathcal{Y}_{2}(\mathcal{U}(x_{2})w_{2}, x_{2})\cdots
\mathcal{Y}_{n}(\mathcal{U}(x_{n})w_{n}, x_{n})q^{L(0)}
\end{eqnarray}
as formal series in $x_{1}, \dots, x_{n}$ and $q$.
Note that the left- and right-hand sides of 
(\ref{identity1-4}) are well-defined formal 
series in $x_{1}, \dots, x_{n}$ and $q$ for $z\in\mathbb{C}$ satisfying 
$1>|q_{z}-1|>0$ and $|q_{z}|>1$, respectively.

By (\ref{ux-u1}) and (\ref{identity0}), for any complex number $z$, we have 
\begin{eqnarray}\label{identity1-5}
\lefteqn{\tr_{\tilde{W}_{n}}Y(\mathcal{U}(x_{1}q_{z})u, x)
\mathcal{Y}_{1}(\mathcal{U}(x_{1})w_{1}, x_{1})\cdot}\nn
&&\quad\quad\quad\quad\quad\cdot
\mathcal{Y}_{2}(\mathcal{U}(x_{2})w_{2}, x_{2})\cdots
\mathcal{Y}_{n}(\mathcal{U}(x_{n})w_{n}, x_{n})q^{L(0)}\nn
&&=\tr_{\tilde{W}_{n}}Y(\mathcal{U}(x)(\mathcal{U}(x))^{-1}
\mathcal{U}(x_{1}q_{z})u, x)
\mathcal{Y}_{1}(\mathcal{U}(x_{1})w_{1}, x_{1})\cdot\nn
&&\quad\quad\quad\quad\quad\cdot
\mathcal{Y}_{2}(\mathcal{U}(x_{2})w_{2}, x_{2})\cdots
\mathcal{Y}_{n}(\mathcal{U}(x_{n})w_{n}, x_{n})q^{L(0)}\nn
&&=\sum_{i=1}^{n}\sum_{m\ge 0}
P_{m+1}\left(\frac{x_{i}}{x}; q\right)\cdot\nn
&&\quad\quad\quad\quad\quad\cdot\tr_{\tilde{W}_{n}}
\mathcal{Y}_{1}(\mathcal{U}(x_{1})w_{1}, x_{1})
\cdots\mathcal{Y}_{i-1}(\mathcal{U}(x_{i-1})w_{i-1}, x_{i-1})\cdot\nn 
&&\quad\quad\quad\quad\quad\cdot
\mathcal{Y}_{i}(\mathcal{U}(x_{i})(2\pi i)^{L(0)}
((\mathcal{U}(x))^{-1}\mathcal{U}(x_{1}q_{z})u)_{m}
w_{i}, x_{i})\cdot\nn
&&\quad\quad\quad\quad\quad\cdot
\mathcal{Y}_{i+1}(\mathcal{U}(x_{i+1})w_{i+1}, x_{i+1})
\cdots \mathcal{Y}_{n}(\mathcal{U}(x_{n})w_{n}, x_{n})q^{L(0)}\nn
&&+\tr_{\tilde{W}_{n}}o(\mathcal{U}(1)
(\mathcal{U}(x))^{-1}\mathcal{U}(x_{1}q_{z})u)\cdot\nn
&&\quad\quad\quad\quad\quad\cdot
\mathcal{Y}_{1}(\mathcal{U}(x_{1})w_{1}, x_{1})
\cdots \mathcal{Y}_{n}(\mathcal{U}(x_{n})w_{n}, x_{n})q^{L(0)}
\end{eqnarray}
as formal series in $x, x_{1}, \dots, x_{n}$ and $q$.

Now for any complex number $z$ satisfying
$|q_{z}|>1$ and $0<|z|<1$, 
substituting $x_{1}q_{z}$ for $x$
in (\ref{identity1-5}) and using (\ref{wp2}), we obtain
\begin{eqnarray}\label{identity1-6}
\lefteqn{\tr_{\tilde{W}_{n}}Y(\mathcal{U}(x_{1}q_{z})u, x_{1}q_{z})
\mathcal{Y}_{1}(\mathcal{U}(x_{1})w_{1}, x_{1})\cdot}\nn
&&\quad\quad\quad\quad\cdot
\mathcal{Y}_{2}(\mathcal{U}(x_{2})w_{2}, x_{2})\cdots
\mathcal{Y}_{n}(\mathcal{U}(x_{n})w_{n}, x_{n})q^{L(0)}\nn
&&=\sum_{m\ge 0}
P_{m+1}\left(\frac{1}{q_{z}}; q\right)\tr_{\tilde{W}_{n}}
\mathcal{Y}_{1}(\mathcal{U}(x_{1})u_{m}w_{1}, x_{1})\cdot\nn
&&\quad\quad\quad\quad\cdot
\mathcal{Y}_{2}(\mathcal{U}(x_{2})w_{2}, x_{2})\cdots
\mathcal{Y}_{n}(\mathcal{U}(x_{n})w_{n}, x_{n})q^{L(0)}\nn
&&\quad+\sum_{i=2}^{n}\sum_{m\ge 0}
P_{m+1}\left(\frac{x_{i}}{x_{1}q_{z}}; q\right)\tr_{\tilde{W}_{n}}
\mathcal{Y}_{1}(\mathcal{U}(x_{1})w_{1}, x_{1})\cdot\nn
&&\quad\quad\quad\quad
\cdots\mathcal{Y}_{i-1}(\mathcal{U}(x_{i-1})w_{i-1}, x_{i-1})
\mathcal{Y}_{i}(\mathcal{U}(x_{i})u_{m}w_{i}, x_{i})\cdot\nn
&&\quad\quad\quad\quad\cdot
\mathcal{Y}_{i+1}(\mathcal{U}(x_{i+1})w_{i+1}, x_{i+1})
\cdots \mathcal{Y}_{n}(\mathcal{U}(x_{n})w_{n}, x_{n})q^{L(0)}\nn
&&\quad+\tr_{\tilde{W}_{n}}o(\mathcal{U}(1)u)
\mathcal{Y}_{1}(\mathcal{U}(x_{1})w_{1}, x_{1})
\cdots \mathcal{Y}_{n}(\mathcal{U}(x_{n})w_{n}, x_{n})q^{L(0)}\nn
&&=\sum_{m\ge 0}(-1)^{m+1}
\left(\tilde{\wp}_{m+1}(-z; q)
+\frac{\partial^{m}}{\partial z^{m}}(\tilde{G}_{2}(q)z+\pi i)\right)
\cdot\nn
&&\quad\quad\quad\quad\cdot \tr_{\tilde{W}_{n}}
\mathcal{Y}_{1}(\mathcal{U}(x_{1})
u_{m}w_{1}, x_{1})\cdot\nn 
&&\quad\quad\quad\quad\cdot 
\mathcal{Y}_{2}(\mathcal{U}(x_{2})w_{2}, x_{2})\cdots
\mathcal{Y}_{n}(\mathcal{U}(x_{n})w_{n}, x_{n})q^{L(0)}\nn
&&\quad+\sum_{i=2}^{n}\sum_{m\ge 0}
P_{m+1}\left(\frac{x_{i}}{x_{1}q_{z}}; q\right)\cdot\nn
&&\quad\quad\quad\quad\cdot\tr_{\tilde{W}_{n}}
\mathcal{Y}_{1}(\mathcal{U}(x_{1})w_{1}, x_{1})
\cdots\mathcal{Y}_{i-1}(\mathcal{U}(x_{i-1})w_{i-1}, x_{i-1})\cdot\nn 
&&\quad\quad\quad\quad\cdot
\mathcal{Y}_{i}(\mathcal{U}(x_{i})u_{m}w_{i}, x_{i})\cdot\nn
&&\quad\quad\quad\quad\cdot
\mathcal{Y}_{i+1}(\mathcal{U}(x_{i+1})w_{i+1}, x_{i+1})
\cdots \mathcal{Y}_{n}(\mathcal{U}(x_{n})w_{n}, x_{n})q^{L(0)}\nn
&&\quad+\tr_{\tilde{W}_{n}}o(\mathcal{U}(1)u)
\mathcal{Y}_{1}(\mathcal{U}(x_{1})w_{1}, x_{1})
\cdots \mathcal{Y}_{n}(\mathcal{U}(x_{n})w_{n}, x_{n})q^{L(0)}
\end{eqnarray}
as a formal series in $x_{1}, \dots, x_{n}$ and $q$.

Combining (\ref{identity1-2}), (\ref{identity1-4}) and 
(\ref{identity1-6}), we see that for any complex number $z$ satisfying
$|q_{z}|>1$ and $0<|z|<1$, 
\begin{eqnarray}\label{identity1-7}
\lefteqn{\tr_{\tilde{W}_{n}}
\mathcal{Y}_{1}(\mathcal{U}(x_{1})Y(u, z)w_{1}, x_{1})
\mathcal{Y}_{2}(\mathcal{U}(x_{2})w_{2}, x_{2})
\cdots
\mathcal{Y}_{n}(\mathcal{U}(x_{n})w_{n}, x_{n})q^{L(0)}}\nn
&&=\sum_{m\ge 0}(-1)^{m+1}
\left(\tilde{\wp}_{m+1}(-z; q)
+\frac{\partial^{m}}{\partial z^{m}}(\tilde{G}_{2}(q)z+\pi i)\right)
\cdot\nn
&&\quad\quad\quad\quad\quad\quad\quad\cdot\tr_{\tilde{W}_{n}}
\mathcal{Y}_{1}(\mathcal{U}(x_{1})w_{1}, x_{1})\cdots 
\mathcal{Y}_{j-1}(\mathcal{U}(x_{j-1})w_{j-1}, x_{j-1})\cdot\nn
&&\quad\quad\quad\quad\quad\quad\quad\cdot
\mathcal{Y}_{j}(\mathcal{U}(x_{j})u_{m}w_{j}, x_{j})
\mathcal{Y}_{j+1}(\mathcal{U}(x_{j+1})w_{j+1}, x_{j+1})\cdot\nn
&&\quad\quad\quad\quad\quad\quad\quad
\cdots \mathcal{Y}_{n}(\mathcal{U}(x_{n})w_{n}, x_{n})q^{L(0)}\nn
&&\quad+\sum_{i=2}^{n}\sum_{m\ge 0}
P_{m+1}\left(\frac{x_{i}}{x_{1}q_{z}}; q\right)\cdot\nn
&&\quad\quad\quad\quad\quad\quad\quad\cdot\tr_{\tilde{W}_{n}}
\mathcal{Y}_{1}(\mathcal{U}(x_{1})w_{1}, x_{1})\cdots 
\mathcal{Y}_{i-1}(\mathcal{U}(x_{i-1})w_{i-1}, x_{i-1})\cdot\nn
&&\quad\quad\quad\quad\quad\quad\quad\cdot
\mathcal{Y}_{i}(\mathcal{U}(x_{i})u_{m}w_{i}, x_{i})
\mathcal{Y}_{i+1}(\mathcal{U}(x_{i+1})w_{i+1}, x_{i+1})\cdot\nn
&&\quad\quad\quad\quad\quad\quad\quad
\cdots \mathcal{Y}_{n}(\mathcal{U}(x_{n})w_{n}, x_{n})q^{L(0)}\nn
&&\quad+\tr_{\tilde{W}_{n}}o(\mathcal{U}(1)u)
\mathcal{Y}_{1}(\mathcal{U}(x_{1})w_{1}, x_{1})
\cdots \mathcal{Y}_{n}(\mathcal{U}(x_{n})w_{n}, x_{n})q^{L(0)}.\nn
&&
\end{eqnarray}
as a formal series in $x_{1}, \dots, x_{n}$ and $q$.
Since both sides of (\ref{identity1-7}) are
formal series in $x_{1}, \dots, x_{n}$ and $q$ whose coefficients
are convergent Laurent series in $z$ near $z=0$, the truth 
of (\ref{identity1-7})
in the region given by $|q_{z}|>1$  and $0<|z|<1$
implies that 
(\ref{identity1-7}) holds near $z=0$. Thus (\ref{identity1})
holds in the case of $j=1$.

We now prove the theorem in the general case. We use induction on $j$. 
Assume that for $j=k<n$, (\ref{identity1}) holds. We want to prove
that (\ref{identity1}) holds when $j=k+1$. Fix  
$z^{0}_{i}\in \mathbb{C}$, $i=1, \dots, n$, satisfying 
$|z^{0}_{1}|>\cdots>|z^{0}_{n}|>0$, let $\gamma_{k}$ be the path in 
$$M^{n}=\{(z_{1}, \dots, z_{n})\in \mathbb{C}^{n}\;|\;
z_{i}\ne 0, \; z_{i}\ne z_{l}, i\ne l\}$$
from $(z^{0}_{1}, \dots, z^{0}_{k}, z^{0}_{k+1}, \dots, z^{0}_{n})$ to 
$(z^{0}_{1}, \dots, z^{0}_{k+1}, z^{0}_{k}, \dots, z^{0}_{n})$ given by 
$$\gamma_{k}(t)=
(z^{0}_{1}, \dots, e^{(1-t)\log z^{0}_{k}+t\log z^{0}_{k+1}}, 
e^{t\log z^{0}_{k}+(1-t)\log z^{0}_{k+1}}, \dots, z^{0}_{n}).$$
(Actually just as the choices of $z^{0}_{i}\in \mathbb{C}$, $i=1, \dots, n$,
the choice of $\gamma_{k}$ can also be arbitrary. 
Here to be more specific, we 
explicitly give one particular choice.)
Given any branch of a multivalued analytic function of 
$z_{1}, \dots, z_{n}$ in the simply connected
region $|z_{1}|>\cdots>|z_{n}|>0$, $0\le \arg z_{i}<2\pi$, 
$i=1, \dots, n$, we have a unique analytic extension of this branch
to the region 
\begin{eqnarray}\label{region}
&|z_{1}|>\cdots >|z_{k-1}|>|z_{k+1}|>|z_{k}|>|z_{k+2}|>
\cdots >|z_{n}|>0, &\nn
&0\le \arg z_{i}<2\pi,\;\;\; i=1, \dots, n,
\end{eqnarray}
determined by the path $\gamma_{k}$. 
By commutativity for 
intertwining operators, there exist $V$-modules $\hat{W}_{k}$
and intertwining operators $\hat{\mathcal{Y}}_{k}$ and 
$\hat{\mathcal{Y}}_{k+1}$ of types $\binom{\hat{W}_{k}}{W_{k}\tilde{W}_{k+1}}$
and $\binom{\tilde{W}_{k-1}}{W_{k+1}\hat{W}_{k}}$, respectively, such that
for any $w_{i}\in W_{i}$,  $i=1, \dots, n$,
$\tilde{w}_{n}\in \tilde{W}_{n}$ and $\tilde{w}'_{n}\in \tilde{W}'_{n}$,
the branch of the analytic extension along $\gamma_{k}$ of 
$$\langle 
\tilde{w}'_{n}, \mathcal{Y}_{1}(w_{1}, z_{1})
\cdots \mathcal{Y}_{k}(w_{k}, z_{k})\mathcal{Y}_{k+1}(w_{k+1}, z_{k+1})\cdots
\mathcal{Y}_{n}(w_{n}, z_{n})\tilde{w}_{n}\rangle$$
to the region (\ref{region}),
is equal to 
\begin{eqnarray*}
\lefteqn{\langle \tilde{w}'_{n}, \mathcal{Y}_{1}(w_{1}, z_{1})
\cdots \mathcal{Y}_{k-1}(w_{k-1}, z_{k-1})
\hat{\mathcal{Y}}_{k+1}(w_{k+1}, z_{k+1})\cdot}\nn
&&\quad\quad\quad\quad\quad\quad\quad\quad\quad\quad\quad\cdot 
\hat{\mathcal{Y}}_{k}(w_{k}, z_{k})\mathcal{Y}_{k+2}(w_{k+2}, z_{k+2})\cdots
\mathcal{Y}_{n}(w_{n}, z_{n})\tilde{w}_{n}\rangle
\end{eqnarray*}
in the same region.

By induction assumption,  (\ref{identity1}) holds when $j=k$. So
(\ref{identity1})
holds when  $w_{j}, w_{j+1}, x_{j}, x_{j+1}$, 
$\mathcal{Y}_{j}$ and $\mathcal{Y}_{j+1}$ in (\ref{identity1})
are replaced
by $w_{k+1}, w_{k}, x_{k+1}, x_{k}$, $\hat{\mathcal{Y}}_{k+1}$ and 
$\hat{\mathcal{Y}}_{k}$, respectively. 
Then when we substitute $z_{1}, \dots, z_{n}\in \mathbb{C}$
in the region (\ref{region})
for $x_{1}, \dots, x_{n}$, respectively, (\ref{identity1}) becomes
an identity for formal Laurent series in $q$ whose coefficients
are single-valued analytic branches in the same region. 
By the commutativity stated above, we know that 
all these coefficients can be analytically extended back to the 
region $|z_{1}|>\cdots>|z_{n}|>0$,
$0\le \arg z_{i}<2\pi$, 
$i=1, \dots, n$,
along the path $\gamma_{k}^{-1}$ to 
the single-valued analytic branches which are nothing but the
coefficients of the formal Laurent series in $q$ obtained 
by substituting in (\ref{identity1})
(with $j=k+1$) $z_{1},  \cdots, z_{n}\in \mathbb{C}$
for $x_{1}, \dots, x_{n}$, respectively, when $z_{1},  \cdots, z_{n}$
are in the region $|z_{1}|>\cdots>|z_{n}|>0$,
$0\le \arg z_{i}<2\pi$, 
$i=1, \dots, n$. 
Since these analytic extensions are unique, 
we see that (\ref{identity1}) with $j=k+1$ holds. 
\epfv

\begin{rema}
{\rm In the proof of the theorem above, the use of the
complex variable approach is essential. If we use formal 
variable approach, in the case of $j=1$, the proof seems to be more 
complicated and in the case of $j>1$, 
certain identities can still be obtained but they
would involve $P(\frac{x_{i}q}{x_{j}e^{2\pi iy}}; q)$ (note the extra $q$
in the first argument)
for $i<j$. (After an early version of the present paper was finished and
circulated, Milas informed
the author that he also obtained independently some
identities similar to (\ref{identity0}), (\ref{identity0.5})
and the $j=1$ case of (\ref{identity1}).)
In fact, eventually we will
have multivalued analytic functions corresponding to flat sections of
certain vector bundles with flat connections over the moduli spaces 
of genus-one Riemann surfaces with punctures and standard local coordinates
vanishing at the punctures.
So if the convergence is proved,  those identities with extra $q$'s and 
the identities in the theorem above are related by analytic extensions.
But if we have only the identities involving extra $q$'s,
we would not even be able to prove the convergence of the series and 
therefore cannot construct the multivalued analytic 
functions. }
\end{rema}

Let $\mathbb{G}_{|z_{1}|>\cdots>|z_{n}|>0}$ be the space
of all multivalued analytic functions in $z_{1}, \dots, z_{n}$
defined on the region
$|z_{1}|>\cdots>|z_{n}|>0$ with preferred branches 
in the simply-connected region $|z_{1}|>\cdots>|z_{n}|>0$,
$0\le \arg z_{i}<2\pi$,
$i=1, \dots, n$. For any $f(z_{1}, \dots, z_{n})\in
\mathbb{G}_{|z_{1}|>\cdots>|z_{n}|>0}$, 
we have a multivalued analytic functions $f(q_{z_{1}}, \dots, q_{z_{n}})$
in $z_{1}, \dots, z_{n}$
defined on the region $|q_{z_{1}}|>\cdots>|q_{z_{n}}|>0$.
All such functions form a space
$\mathbb{G}_{|q_{z_{1}}|>\cdots>|q_{z_{n}}|>0}$. 
We see that 
$$\tr_{\tilde{W}_{n}}\mathcal{Y}_{1}(\mathcal{U}(q_{z_{1}})w_{1}, q_{z_{1}})
\cdots
\mathcal{Y}_{n}(\mathcal{U}(q_{z_{n}})w_{n}, q_{z_{n}})q^{L(0)}$$
and 
\begin{eqnarray*}
&\tilde{G}_{2k+2}(q)
\tr_{\tilde{W}_{n}}\mathcal{Y}_{1}(\mathcal{U}(q_{z_{1}})w_{1}, q_{z_{1}})
\cdots
\mathcal{Y}_{n}(\mathcal{U}(q_{z_{n}})w_{n}, q_{z_{n}})q^{L(0)}
\end{eqnarray*}
for $k\ge 1$ are elements of 
$\mathbb{G}_{|q_{z_{1}}|>\cdots>|q_{z_{n}}|>0}((q))$.

\begin{thm}
Let  $W_{i}$ and $\tilde{W}_{i}$,  $\mathcal{Y}_{i}$, $i=1, \dots, n$,
be as in Proposition \ref{identity0-pr}.
Then for any $u\in V$, $w_{i}\in W_{i}$, $i=1, \dots, n$, any 
integer $j$ satisfying $1\le j\le n$ and any $l\in \mathbb{Z}_{+}$, 
in $\mathbb{G}_{|q_{z_{1}}|>\cdots>|q_{z_{n}}|>0}((q))$, 
we have
\begin{eqnarray}\label{identity2}
\lefteqn{\tr_{\tilde{W}_{n}}
\mathcal{Y}_{1}(\mathcal{U}(q_{z_{1}})w_{1}, q_{z_{1}})
\cdots \mathcal{Y}_{j-1}(\mathcal{U}(q_{z_{j-1}})w_{j-1}, q_{z_{j-1}})
\mathcal{Y}_{j}(\mathcal{U}(q_{z_{j}})u_{-l}w_{j}, q_{z_{j}})\cdot}\nn
&&\quad\quad\cdot 
\mathcal{Y}_{j+1}(\mathcal{U}(q_{z_{j+1}})w_{j+1}, q_{z_{j+1}})
\cdots
\mathcal{Y}_{n}(\mathcal{U}(q_{z_{n}})w_{n}, q_{z_{n}})q^{L(0)}\nn
&&=\sum_{k\ge 1}(-1)^{l+1}\binom{2k+1}{l-1}\tilde{G}_{2k+2}(q)
\tr_{\tilde{W}_{n}}\mathcal{Y}_{1}(\mathcal{U}(q_{z_{1}})w_{1}, q_{z_{1}})
\cdot \nn
&& \quad\quad\cdots \mathcal{Y}_{j-1}(\mathcal{U}(q_{z_{j-1}})w_{j-1}, 
q_{z_{j-1}})
\mathcal{Y}_{j}(\mathcal{U}(q_{z_{j}})u_{2k+2-l}w_{j}, q_{z_{j}})\cdot\nn
&&\quad\quad\cdot 
 \mathcal{Y}_{j+1}(\mathcal{U}(q_{z_{j+1}})w_{j+1}, q_{z_{j+1}})
\cdots
\mathcal{Y}_{n}(\mathcal{U}(q_{z_{n}})w_{n}, q_{z_{n}})q^{L(0)}\nn
&&\quad+\sum_{i\ne j}\sum_{m\ge 0}(-1)^{m+l}
\binom{-m-1}{l-1}\tilde{\wp}_{m+l}(z_{i}-z_{j}; q)
\tr_{\tilde{W}_{n}}\mathcal{Y}_{1}(\mathcal{U}(q_{z_{1}})w_{1}, 
q_{z_{1}})\cdot\nn
&&\quad\quad
\cdots \mathcal{Y}_{i-1}(\mathcal{U}(q_{z_{i-1}})w_{i-1}, q_{z_{i-1}})
\mathcal{Y}_{i}(\mathcal{U}(q_{z_{i}})u_{m}w_{i}, q_{z_{i}})\cdot\nn
&&\quad\quad\cdot 
\mathcal{Y}_{i+1}(\mathcal{U}(q_{z_{i+1}})w_{i+1}, q_{z_{i+1}})
\cdots
\mathcal{Y}_{n}(\mathcal{U}(q_{z_{n}})w_{n}, q_{z_{n}})q^{L(0)}\nn
&&\quad +\delta_{l, 1}\tilde{G}_{2}(q)\sum_{i=1}^{n}
\tr_{\tilde{W}_{n}}\mathcal{Y}_{1}(\mathcal{U}(q_{z_{1}})w_{1}, q_{z_{1}})
\cdots \mathcal{Y}_{i-1}(\mathcal{U}(q_{z_{i-1}})w_{i-1}, q_{z_{i-1}})\cdot\nn
&&\quad\quad\cdot 
\mathcal{Y}_{i}(\mathcal{U}(q_{z_{i}})(u_{1}+u_{0}z_{i})w_{i}, 
q_{z_{i}})
\mathcal{Y}_{i+1}(\mathcal{U}(q_{z_{i+1}})w_{i+1}, q_{z_{i+1}})\cdot\nn
&&\quad\quad\cdots
\mathcal{Y}_{n}(\mathcal{U}(q_{z_{n}})w_{n}, q_{z_{n}})q^{L(0)}\nn
&&\quad+\delta_{l, 1}\tr_{\tilde{W}_{n}}o(\mathcal{U}(1)u)
\mathcal{Y}_{1}(\mathcal{U}(q_{z_{1}})w_{1}, q_{z_{1}})
\cdots \mathcal{Y}_{n}(\mathcal{U}(q_{z_{n}})w_{n}, q_{z_{n}})q^{L(0)}.
\end{eqnarray}
\end{thm}
\pf
Let $z_{1}, \dots, z_{n}\in \mathbb{C}$ satisfying 
$|q_{z_{1}}|>\cdots>|q_{z_{n}}|>0$. Then  from (\ref{identity1-1})
and the domain of convergence of the $q$-coefficients of 
$P_{m}(\frac{q_{z_{i}}}{q_{z_{j}}q_{z}}, q)$ for $m\ge 1$, we see that
for $z$ in a sufficiently small neighborhood of $0$, we can substitute 
$z, q_{z_{1}}, \dots, q_{z_{n}}$ for $y, x_{1}, \dots, x_{n}$
in (\ref{identity1}) so that both sides of (\ref{identity1}) become
formal series in $q$ whose coefficients are 
analytic functions of $z$. Using (\ref{wp2}) after 
these substitutions, 
we obtain 
\begin{eqnarray}\label{identity2-1}
\lefteqn{\tr_{\tilde{W}_{n}}\mathcal{Y}_{1}(\mathcal{U}(q_{z_{1}})w_{1}, 
q_{z_{1}})
\cdots \mathcal{Y}_{j-1}(\mathcal{U}(q_{z_{j-1}})w_{j-1}, q_{z_{j-1}})
\mathcal{Y}_{j}(\mathcal{U}(q_{z_{j}})Y(u, z)w_{j}, q_{z_{j}})\cdot}\nn
&&\quad\quad\quad\quad\cdot 
\mathcal{Y}_{j+1}(\mathcal{U}(q_{z_{j+1}})w_{j+1}, q_{z_{j+1}})
\cdots
\mathcal{Y}_{n}(\mathcal{U}(q_{z_{n}})w_{n}, q_{z_{n}})q^{L(0)}\nn
&&=\sum_{m\ge 0}(-1)^{m+1}
\left(\tilde{\wp}_{m+1}(-z; q)+\frac{\partial^{m}}{\partial z^{m}}
(\tilde{G}_{2}(q)z+\pi i)\right)
\cdot\nn
&&\quad\quad\quad\quad\cdot\tr_{\tilde{W}_{n}}
\mathcal{Y}_{1}(\mathcal{U}(q_{z_{1}})w_{1}, q_{z_{1}})\cdots 
\mathcal{Y}_{j-1}(\mathcal{U}(q_{z_{j-1}})w_{j-1}, q_{z_{j-1}})\cdot\nn
&&\quad\quad\quad\quad\cdot
\mathcal{Y}_{j}(\mathcal{U}(q_{z_{j}})u_{m}w_{j}, q_{z_{j}})
\mathcal{Y}_{j+1}(\mathcal{U}(q_{z_{j+1}})w_{j+1}, q_{z_{j+1}})\cdot\nn
&&\quad\quad\quad\quad
\cdots \mathcal{Y}_{n}(\mathcal{U}(q_{z_{n}})w_{n}, q_{z_{n}})q^{L(0)}\nn
&&\quad+\sum_{i\ne j}^{}\sum_{m\ge 0}(-1)^{m+1}
\Biggl(\tilde{\wp}_{m+1}(z_{i}-z_{j}-z; q)\nn
&&\quad\quad\quad\quad\quad\quad\quad\quad\quad\quad\quad
+(-1)^{m+1}\frac{\partial^{m}}{\partial z_{i}^{m}}
(\tilde{G}_{2}(q)(z_{i}-z_{j}-z)-\pi i)\Biggr)\cdot\nn
&&\quad\quad\quad\quad\cdot\tr_{\tilde{W}_{n}}
\mathcal{Y}_{1}(\mathcal{U}(q_{z_{1}})w_{1}, q_{z_{1}})\cdots 
\mathcal{Y}_{i-1}(\mathcal{U}(q_{z_{i-1}})w_{i-1}, q_{z_{i-1}})\cdot\nn
&&\quad\quad\quad\quad\cdot
\mathcal{Y}_{i}(\mathcal{U}(q_{z_{i}})u_{m}w_{i}, q_{z_{i}})
\mathcal{Y}_{i+1}(\mathcal{U}(q_{z_{i+1}})w_{i+1}, q_{z_{i+1}})\cdot\nn
&&\quad\quad\quad\quad
\cdots \mathcal{Y}_{n}(\mathcal{U}(q_{z_{n}})w_{n}, q_{z_{n}})q^{L(0)}\nn
&&\quad+\tr_{\tilde{W}_{n}}o(\mathcal{U}(1)u)
\mathcal{Y}_{1}(\mathcal{U}(q_{z_{1}})w_{1}, q_{z_{1}})
\cdots \mathcal{Y}_{n}(\mathcal{U}(q_{z_{n}})w_{n}, q_{z_{n}})q^{L(0)}.
\end{eqnarray}
For $l>2$, taking the $z^{l-1}$-coefficients of 
both sides of (\ref{identity2-1})
and using (\ref{wp-relation}), (\ref{wp1}) and 
the Taylor expansion formulas for
$\tilde{\wp}_{m+1}(z_{i}-z_{j}-z; q)$ as functions of $z$, 
we obtain (\ref{identity2}) in this case. 
For $l=2$, the right-hand side contains an
extra term 
\begin{eqnarray*}
\lefteqn{-\tilde{G}_{2}(q)\sum_{i=1}^{n}\tr_{\tilde{W}_{n}}
\mathcal{Y}_{1}(\mathcal{U}(q_{z_{1}})w_{1}, q_{z_{1}})\cdots 
\mathcal{Y}_{i-1}(\mathcal{U}(q_{z_{i-1}})w_{i-1}, q_{z_{i-1}})\cdot}\nn
&&\cdot
\mathcal{Y}_{i}(\mathcal{U}(q_{z_{i}})u_{0}w_{i}, q_{z_{i}})
\mathcal{Y}_{i+1}(\mathcal{U}(q_{z_{i+1}})w_{i+1}, q_{z_{i+1}})
\cdots \mathcal{Y}_{n}(\mathcal{U}(q_{z_{n}})w_{n}, q_{z_{n}})q^{L(0)}
\end{eqnarray*}
which is equal to $0$ by (\ref{identity0.5}). 
For $l=1$, the right-hand side contains the 
two extra terms in the right-hand side of (\ref{identity2}) 
and also contains terms 
\begin{eqnarray*}
\lefteqn{-\pi \sum_{i=1}^{n}\tr_{\tilde{W}_{n}}
\mathcal{Y}_{1}(\mathcal{U}(q_{z_{1}})w_{1}, q_{z_{1}})\cdots 
\mathcal{Y}_{i-1}(\mathcal{U}(q_{z_{i-1}})w_{i-1}, q_{z_{i-1}})\cdot}\nn
&&\cdot
\mathcal{Y}_{i}(\mathcal{U}(q_{z_{i}})u_{0}w_{i}, q_{z_{i}})
\mathcal{Y}_{i+1}(\mathcal{U}(q_{z_{i+1}})w_{i+1}, q_{z_{i+1}})
\cdots \mathcal{Y}_{n}(\mathcal{U}(q_{z_{n}})w_{n}, q_{z_{n}})q^{L(0)}
\end{eqnarray*}
and 
\begin{eqnarray*}
\lefteqn{-z_{j} \sum_{i=1}^{n}\tr_{\tilde{W}_{n}}
\mathcal{Y}_{1}(\mathcal{U}(q_{z_{1}})w_{1}, q_{z_{1}})\cdots 
\mathcal{Y}_{i-1}(\mathcal{U}(q_{z_{i-1}})w_{i-1}, q_{z_{i-1}})\cdot}\nn
&&\cdot
\mathcal{Y}_{i}(\mathcal{U}(q_{z_{i}})u_{0}w_{i}, q_{z_{i}})
\mathcal{Y}_{i+1}(\mathcal{U}(q_{z_{i+1}})w_{i+1}, q_{z_{i+1}})
\cdots \mathcal{Y}_{n}(\mathcal{U}(q_{z_{n}})w_{n}, q_{z_{n}})q^{L(0)}
\end{eqnarray*}
which are again equal to $0$ by (\ref{identity0.5}). 

Since both sides of (\ref{identity2})
are in $\mathbb{G}_{|q_{z_{1}}|>\cdots>|q_{z_{n}}|>0}((q))$, 
(\ref{identity2}) hold in this space.
\epfv

\renewcommand{\theequation}{\thesection.\arabic{equation}}
\renewcommand{\thethm}{\thesection.\arabic{thm}}
\setcounter{equation}{0}
\setcounter{thm}{0}

\section{Differential equations}

We shall fix a positive integer $n$ in this section. Let 
$$R=\mathbb{C}[\tilde{G}_{4}(q), \tilde{G}_{6}(q),
\tilde{\wp}_{2}(z_{i}-z_{j}; q), 
\tilde{\wp}_{3}(z_{i}-z_{j}; q)]_{i, j=1, \dots, n,\;i<j},$$
that is, the commutative associative 
algebra over $\mathbb{C}$ generated by the series
$\tilde{G}_{4}(q)$, $\tilde{G}_{6}(q)$, 
$\tilde{\wp}_{2}(z_{i}-z_{j}; q)$ and $\tilde{\wp}_{3}(z_{i}-z_{j}; q)$
for $i, j=1, \dots, n$ satisfying $i<j$. (Note that by definition, 
when $n=1$,
$R$ is the commutative associative 
algebra over $\mathbb{C}$ generated by the series
$\tilde{G}_{4}(q)$ and $\tilde{G}_{6}(q)$.) 
By definition, $R$ is finitely generated over the field $\mathbb{C}$
and thus  is Noetherian. 
We have:

\begin{lemma}
For $k\ge 2$, $m\ge 2$, $i, j=1, \dots, n$, $i\ne j$,
$$\tilde{G}_{2k}(q), \tilde{\wp}_{m}(z_{i}-z_{j}; q)\in R.$$
\end{lemma}
\pf
It is known (see, for example, \cite{Ko}) that  any 
modular form can be written
as a polynomial of $G_{4}(\tau)$ and $G_{6}(\tau)$. In particular, 
$G_{2k}(\tau)$ for $k\ge 2$ can be written as polynomials of
$G_{4}(\tau)$ and $G_{6}(\tau)$. Since $\tilde{G}_{2k}(q)$ for 
$k\ge 2$ are nothing but the $q$-expansions of $G_{2k}(\tau)$,
they can be written as polynomials of 
$\tilde{G}_{4}(q)$ and $\tilde{G}_{6}(q)$ and thus are elements of 
$R$. 
Also using the relation
$$\frac{\partial}{\partial z}
\wp_{2}(z; \tau)^{2}=4\wp_{2}(z; \tau)^{3}-60G_{4}(\tau)\wp_{2}(z; \tau)
-140G_{6}(\tau)$$
(see, for example, \cite{Ko} or \cite{La}) and induction, 
it is easy to see that $\wp_{m}(z, \tau)$ for $m\ge 2$ 
can be written as polynomials of $G_{4}(\tau)$, $G_{6}(\tau)$,
$\wp_{2}(z; \tau)$ and 
$$\frac{\partial}{\partial z}
\wp_{2}(z; \tau)=-2\wp_{3}(z; \tau).$$
Since $\tilde{G}_{4}(q)$,  $\tilde{G}_{6}(q)$ and 
$\tilde{\wp}_{m}(z; q)$ for 
$m\ge 2$ are nothing but the $q$-expansion of 
$G_{4}(\tau)$, $G_{6}(\tau)$ and $\wp_{m}(z; \tau)$, respectively,
$\tilde{\wp}_{m}(z; q)$ for 
$m\ge 2$ can be written as polynomials of 
$\tilde{G}_{4}(q)$,  $\tilde{G}_{6}(q)$,
$\tilde{\wp}_{2}(z; q)$ and $\tilde{\wp}_{3}(z; q)$.
Consequently for $m\ge 2$ and $i, j=1, \dots, n$, $i< j$, 
$\tilde{\wp}_{m}(z_{i}-z_{j}; q)$ can be written as polynomials of 
$\tilde{G}_{4}(q)$, $\tilde{G}_{6}(q)$,
$\tilde{\wp}_{2}(z_{i}-z_{j}; q)$ and $\tilde{\wp}_{3}(z_{i}-z_{j}; q)$
and thus are elements of $R$. Since $\tilde{\wp}_{m}(z; q)$ for 
$m\ge 2$ are either even or odd in $z$, 
$\tilde{\wp}_{m}(z_{i}-z_{j}; q)$ for 
$m\ge 2$, $i, j=1, \dots, n$, $i> j$ are also in $R$.
\epfv

As in the preceding two sections, we fix a vertex operator algebra
$V$ of central charge $c$ in this section. 
We say that a $V$-module $W$ is {\it $C_{2}$-cofinite} or
satisfies the {\it $C_{2}$-cofiniteness condition} if 
$\dim W/C_{2}(W)<\infty$, where
$C_{2}(W)$ is the subspace of $W$ spanned by elements of the form
$u_{-2}w$ for $u\in V$ and $w\in W$. 
In this section, we assume that all $V$-modules are
$\mathbb{R}$-graded and $C_{2}$-cofinite. We also assume that
for any $W_{i}$ and $\tilde{W}_{i}$,  
$\mathcal{Y}_{i}$, $i=1, \dots, n$,
as in Proposition \ref{identity0-pr}, 
any $w_{i}\in W_{i}$ ($i=1, \dots, n$), $\tilde{w}_{n}\in 
\tilde{W}_{n}$ and $\tilde{w}'_{n}\in  \tilde{W}'_{n}$,
$$\langle \tilde{w}'_{n}, \mathcal{Y}_{1}\left(w_{1}, 
z_{1}\right)\cdots 
\mathcal{Y}_{n}\left(w_{n}, 
z_{1}\right)\tilde{w}_{n}\rangle_{\tilde{W}_{n}}$$
is absolutely convergent when $|z_{1}|>\cdots>|z_{n}|>0$. We also
assume that the commutativity and associativity 
for intertwining operators hold. Note that all
the conditions on the vertex operator algebra $V$ assumed in the 
preceding two sections hold. So we can use the results of these sections.

Let  $T=R \otimes W_{1}\otimes \cdots \otimes W_{n}$.
Then $T$ and also $\mathbb{G}_{|q_{z_{1}}|>\cdots>|q_{z_{n}}|>0}((q))$ 
are $R$-modules.

Let $J$ be the $R$-submodule of 
$T$ generated by elements of the 
form
\begin{eqnarray*}
\lefteqn{\mathcal{A}_{j}(u; w_{1}, \dots, w_{n})}\nn
&&=1\otimes w_{1}\otimes \cdots \otimes w_{j-1}\otimes u_{-2}w_{j}
\otimes w_{j+1}\otimes \cdots \otimes w_{n}\nn
&&+\sum_{k\ge 1}(2k+1)\tilde{G}_{2k+2}(q)\otimes 
w_{1}\otimes \cdots \otimes w_{j-1}\otimes u_{2k}w_{j}
\otimes w_{j+1}\otimes \cdots \otimes w_{n}\nn
&&+\sum_{i\ne j}\sum_{m\ge 0}(-1)^{m}
(m+1)\nn
&&\quad\quad\quad
\tilde{\wp}_{m+2}(z_{i}-z_{j}; q)\otimes w_{1}\otimes \cdots\otimes 
w_{i-1}\otimes u_{m}w_{i}\otimes w_{i+1}\otimes
\cdots \otimes w_{n}
\end{eqnarray*}
for $j=1, \dots, n$ and $w_{i}\in W_{i}$, $i=1, \dots, n$.

The gradings (by conformal weights)
on $W_{i}$ for $i=1, \dots, n$ induce a grading on 
$T$ and this grading on $T$ induces a grading on $J$. We shall use 
$T_{(r)}$ and $J_{(r)}$ to denote the homogeneous subspaces of 
$T$ and $J$, respectively, of conformal weight $r\in \mathbb{R}$. Let 
$F_{r}(T)=\coprod_{s\le r}T_{(s)}$ and 
$F_{r}(J)=\coprod_{s\le r}J_{(s)}$. We also introduce another 
grading on $R$ and $T$. We say that the elements 
$\tilde{G}_{2k}(q)$ for any $k\ge 1$  have {\it modular weight}
$2k$ and the element $\tilde{\wp}_{m}(z_{i}-z_{j}; q)$ 
for any $m\ge 2$ and $i<j$ have {\it modular weight} $m$. 
These modular weights of the generators of $R$ give a 
grading which is called the grading by {\it modular weights}.
For $m\in \mathbb{Z}$, let $R_{m}$ be the 
subspace of $R$ consisting of all elements of 
modular weight $m$. Then we have $R_{m}=0$ if $m<0$ 
and $R=\oplus_{m\in \mathbb{N}}R_{m}$.
An element of $T$ is 
said to have {\it modular weight $m$} if it is a linear combination of 
elements of the form $f\otimes w_{1}\otimes \cdots\otimes w_{n}$
where $w_{1}\in W_{1}, \dots, w_{n}\in W_{n}$ are 
homogeneous and $f\in R$ has modular weight $m-\sum_{i=1}^{n}\wt w_{i}$.
We have:

\begin{prop}\label{md-wt-aj}
Let $w_{1}\in W_{1}, \dots, w_{n}\in W_{n}$ be homogeneous.
For $j=1, \dots, n$, $\mathcal{A}_{j}$ have modular weight
$\sum_{i=1}^{n}\wt w_{i} +2$. In particular, the grading by modular 
weights on $T$ induce a grading, also called the grading by 
modular weights, 
on $T/J$.
\end{prop}
\pf
The conclusion follows directly from the definition of $\mathcal{A}_{j}$.
\epfv

\begin{prop}\label{decomposition}
There exists $N\in \mathbb{R}$ such that 
for any $r\in \mathbb{R}$, $F_{r}(T)\subset 
F_{r}(J)+F_{N}(T)$.
In particular, $T=J+F_{N}(T)$.
\end{prop}
\pf
Since $\dim W_{i}/C_{2}(W_{i})<\infty$ for $i=1, \dots, n$,
there exists $N\in \mathbb{R}$ such that 
\begin{equation}\label{subset}
\coprod_{r>N}T_{(r)}\subset \sum_{j=1}^{n}
R \otimes W_{1}
\otimes \cdots \otimes W_{j-1}\otimes C_{2}(W_{j})
\otimes W_{j+1}\otimes 
W_{n}.
\end{equation}

We use induction on $r\in \mathbb{R}$. 
If $r$ is equal to $N$, 
$F_{N}(T)\subset F_{N}(J)+F_{N}(T)$. Now we assume that 
$F_{r}(T)\subset F_{r}(J)+F_{N}(T)$ for $r< s$ where $s>N$.
We want to show that any homogeneous element of 
$T_{(s)}$ can  be written as
a sum of an element of $F_{s}(J)$ and an element of $F_{N}(T)$.
Since $s>N$, by (\ref{subset}), any element of $T_{(s)}$ is an element of 
the right-hand side of (\ref{subset}). Thus this element of $T_{(s)}$
is a sum of elements of 
$$R \otimes W_{1}
\otimes \cdots \otimes W_{j-1}\otimes C_{2}(W_{j})
\otimes W_{j+1}\otimes 
W_{n}$$
for $j=1, \dots, n$. 
So we need only discuss elements of the form 
\begin{equation}\label{decomp-1}
1\otimes w_{1}\otimes \cdots \otimes w_{j-1}\otimes u_{-2}w_{j}\otimes w_{j+1}
\otimes \cdots\otimes w_{n}
\end{equation}
where 
$w_{i}\in W_{i}$ for $i=1, \dots, n$ and $u\in V$.  
By assumption, the conformal weight of (\ref{decomp-1})
is $s$
and the conformal weights of 
$$1\otimes w_{1}\otimes \cdots \otimes w_{j-1}\otimes u_{2k}w_{j}
\otimes w_{j+1}\otimes \cdots \otimes w_{n}$$
for $k\in \mathbb{Z}_{+}$
and
$$1\otimes w_{1}\otimes \cdots\otimes 
w_{i-1}\otimes u_{m}w_{i}\otimes w_{i+1}\otimes
\cdots \otimes w_{n}$$
for $i\ne j$ and $m\ge 0$
are all less than the conformal weight $s$ of (\ref{decomp-1}).
So $\mathcal{A}_{j}(u;  w_{1}, \dots, w_{n})
\in F_{s}(J)$.
Thus (\ref{decomp-1})
can be written as a sum of an element of $F_{s}(J)$ and 
elements of $T$ of conformal weights less than $s$. 
Then by the induction assumption, 
we know that (\ref{decomp-1})
can be written as a sum of an element of $F_{s}(J)$ and 
an element of $F_{N}(T)$. 

Now we have 
\begin{eqnarray*}
T&=&\cup_{r\in \mathbb{R}}F_{r}(T)\nn
&\subset &\cup_{r\in \mathbb{R}}
F_{r}(J)+F_{N}(T)\nn
&=&J+F_{N}(T).
\end{eqnarray*}
But we know that 
$J+F_{N}(T)\subset T$. So we have 
$T=J+F_{N}(T)$.
\epfv

\begin{cor}
The quotient $R$-module $T/J$ is finitely generated.
\end{cor}
\pf
Since $T=J+F_{N}(T)$ and $F_{N}(T)$
is finitely-generated, $T/J$ is finitely-generated.
\epfv

For $V$-modules $W_{i}$
and $\tilde{W}_{i}$, $i=1, \dots, n$,  and intertwining operators
$\mathcal{Y}_{i}$, $i=1, \dots, n$,  
as in Proposition \ref{identity0-pr}, and 
for any $w_{i}\in W_{i}$, $i=1, \dots, n$,
we shall consider the element
\begin{eqnarray}\label{correl-fn}
\lefteqn{F_{\mathcal{Y}_{1}, \dots, \mathcal{Y}_{n}}(w_{1}, \dots, w_{n};
z_{1}, \dots, z_{n}; q)}\nn
&&=\tr_{\tilde{W}_{n}}\mathcal{Y}_{1}(\mathcal{U}(q_{z_{1}})w_{1}, q_{z_{1}})
\cdots
\mathcal{Y}_{1}(\mathcal{U}(q_{z_{n}})w_{n}, q_{z_{n}})
q^{L(0)-\frac{c}{24}}
\end{eqnarray}
of $\mathbb{G}_{|q_{z_{1}}|>\cdots>|q_{z_{n}}|>0}((q))$.
The map from 
$W_{1}\otimes \cdots\otimes W_{n}$ to $N$ given by 
$$w_{1}\otimes \cdots\otimes w_{n}\mapsto 
F_{\mathcal{Y}_{1}, \dots, \mathcal{Y}_{n}}(w_{1}, \dots, w_{n};
z_{1}, \dots, z_{n}; q)$$
can be naturally extended to an $R$-module map
$$\psi_{\mathcal{Y}_{1}, \dots, \mathcal{Y}_{n}}:
T\to \mathbb{G}_{|q_{z_{1}}|>\cdots>|q_{z_{n}}|>0}((q)).$$

\begin{rema}\label{geometry}
{\rm We consider 
(\ref{correl-fn})
instead of 
$$\tr_{\tilde{W}_{n}}\mathcal{Y}_{1}\left(
q_{z_{1}}^{L(0)}w_{1}, 
q_{z_{1}}\right)\cdots 
\mathcal{Y}_{n}\left(q_{z_{n}}^{L(0)}w_{n}, 
q_{z_{1}}\right)q^{L(0)-\frac{c}{24}}$$
because (\ref{correl-fn}) satisfies simple identities and 
have simple duality properties. Geometrically,
(\ref{correl-fn})
corresponds to a torus with punctures and local coordinates
whose local coordinates at the punctures are the standard ones
in terms of the global coordinates on the torus obtained from
the parallelogram defining the torus. In fact, to construct 
and study chiral genus-one correlation functions, we need to 
use the description of tori in terms of parallelograms. But geometrically
$q$-traces of products of intertwining operators correspond to
annuli, not such parallelograms. Thus we have to use the conformal 
transformation $w\mapsto \tilde{w}=q_{w}$ to map parallelograms to 
annuli. Under this map, the standard local coordinates
$\tilde{w}-q_{z_{i}}$ vanishing at $q_{z_{i}}$ 
for $i=1, \dots, n$ on an annulus 
are not pulled
back to the standard local coordinates $w-z_{i}$ vanishing at $z_{i}$,
respectively,
on the corresponding parallelogram. The operators $\mathcal{U}(q_{z_{i}})$
for $i=1, \dots, n$ correspond exactly to the local coordinates vanishing 
at $q_{z_{i}}$, respectively, on the annulus such that their pull-backs
to the parallelogram are the standard local coordinates 
vanishing at $z_{i}$, respectively.}
\end{rema}

We have:

\begin{prop}
The $R$-submodule $J$ of $T$ is in the kernel of 
$\psi_{\mathcal{Y}_{1}, \dots, \mathcal{Y}_{n}}$.
\end{prop}
\pf
This result follows from the definitions of $J$ and
$\psi_{\mathcal{Y}_{1}, \dots, \mathcal{Y}_{n}}$
and the case $l=2$ of (\ref{identity2})
immediately.
\epfv

\begin{prop}[$L(-1)$-derivative property]\label{l-1-der}
Let $W_{i}$ and $\tilde{W}_{i}$
for $i=1, \dots, n$ be $V$-modules and 
$\mathcal{Y}_{i}$ intertwining operators  of 
types $\binom{\tilde{W}_{i-1}}{W_{i}\tilde{W}_{i}}$ ($i=1, \dots, n$,
$\tilde{W}_{0}=\tilde{W}_{n}$), 
respectively. Then 
\begin{eqnarray*}
\lefteqn{\frac{\partial}{\partial z_{j}}
F_{\mathcal{Y}_{1}, \dots, \mathcal{Y}_{n}}(w_{1}, 
\dots, w_{n};
z_{1}, \dots, z_{n}; q)}\nn
&&=F_{\mathcal{Y}_{1}, \dots, \mathcal{Y}_{n}}(w_{1}, \dots, w_{j-1},
L(-1)w_{j}, w_{j+1},
\dots, w_{n};
z_{1}, \dots, z_{n}; q)
\end{eqnarray*}
for $1\le j\le n$ and $w_{i}\in W_{i}$, $i=1, \dots, n$.
\end{prop}
\pf
The result follows immediately from the definition (\ref{correl-fn})
and Proposition \ref{l-1-e(y)}.
\epfv

\begin{lemma}
Let $W_{i}$ and $\tilde{W}_{i}$
for $i=1, \dots, n$ be $V$-modules and 
$\mathcal{Y}_{i}$ intertwining operators  of 
types $\binom{\tilde{W}_{i-1}}{W_{i}\tilde{W}_{i}}$ ($i=1, \dots, n$,
$\tilde{W}_{0}=\tilde{W}_{n}$), 
respectively. Then for $j=1, \dots, n$ and any homogeneous 
elements $w_{1}\in W_{1}, \dots, w_{n}\in W_{n}$,
\begin{eqnarray}\label{mod-inv-der}
\lefteqn{\left((2\pi i)^{2}q\frac{\partial}{\partial q}
+\tilde{G}_{2}(q)\sum_{i=1}^{n}\wt w_{i}
+\tilde{G}_{2}(q)\sum_{i=1}^{n}z_{i}\frac{\partial}{\partial z_{i}}
-\sum_{i\ne j}\tilde{\wp}_{1}(z_{i}-z_{j}; q)
\frac{\partial}{\partial z_{i}}\right)\cdot}\nn
&& \quad\quad\quad\cdot
F_{\mathcal{Y}_{1}, \dots, \mathcal{Y}_{n}}(w_{1}, 
\dots, w_{n};
z_{1}, \dots, z_{n}; q)\nn
&&=F_{\mathcal{Y}_{1}, \dots, \mathcal{Y}_{n}}(w_{1}, 
\dots, w_{j-1}, L(-2)w_{j}, w_{j+1}, \dots, w_{n};
z_{1}, \dots, z_{n}; q)\nn
&&\quad -\sum_{k\ge 1}\tilde{G}_{2k+2}(q)\cdot\nn
&& \quad\quad\quad\cdot
F_{\mathcal{Y}_{1}, \dots, \mathcal{Y}_{n}}(w_{1}, 
\dots, w_{j-1}, L(2k)w_{j}, w_{j+1}, \dots, w_{n};
z_{1}, \dots, z_{n}; q)\quad\quad\quad\nn
&&\quad+\sum_{i\ne j}\sum_{m\ge 1}(-1)^{m}
\tilde{\wp}_{m+1}(z_{i}-z_{j}; q)\cdot\nn
&& \quad\quad\quad\cdot
F_{\mathcal{Y}_{1}, \dots, \mathcal{Y}_{n}}(w_{1}, 
\dots, w_{i-1}, L(m-1)w_{i}, w_{i+1}, \dots, w_{n};
z_{1}, \dots, z_{n}; q).\nn
&&
\end{eqnarray}
\end{lemma}
\pf
In (\ref{identity2}), we take $u=\omega$ and $l=1$. Using
$\omega_{m}=L(m-1)$ for $m\in \mathbb{Z}$ and Lemma
\ref{u1-omega}, we obtain (\ref{mod-inv-der}).
\epfv

For simplicity, we introduce, for any $\alpha\in \mathbb{C}$,
the notation
$$\mathcal{O}_{j}(\alpha)=(2\pi i)^{2}q\frac{\partial}{\partial q}
+\tilde{G}_{2}(q)\alpha
+\tilde{G}_{2}(q)\sum_{i=1}^{n}z_{i}\frac{\partial}{\partial z_{i}}
-\sum_{i\ne j}\tilde{\wp}_{1}(z_{i}-z_{j}; q)
\frac{\partial}{\partial z_{i}}$$
for $j=1, \dots, n$. We shall also use the notation
$$\prod_{j=1}^{m}\mathcal{O}(\alpha_{j})$$
to denote
$$\mathcal{O}(\alpha_{1})\cdots \mathcal{O}(\alpha_{m}).$$
Note the order of the product.
We have:

\begin{thm}\label{sys1}
Let $V$ be a vertex operator algebra satisfying the conditions 
stated in the beginning of this section and let 
$W_{i}$ for $i=1, \dots, n$ be $V$-modules. Then 
for any homogeneous $w_{i}\in W_{i}$ ($i=1, \dots, n$), there exist
$$a_{p, \;i}(z_{1}, \dots, z_{n}; q)\in 
R_{p}, b_{p, \;i}(z_{1}, \dots, z_{n}; q)\in R_{2p}$$
for $p=1, \dots, m$ and $i=1, \dots, n$ such that 
for any $V$-modules $\tilde{W}_{i}$ ($i=1, \dots, n$) and
intertwining operators
$\mathcal{Y}_{i}$   of 
types $\binom{\tilde{W}_{i-1}}{W_{i}\tilde{W}_{i}}$ ($i=1, \dots, n$,
$\tilde{W}_{0}=\tilde{W}_{n}$), 
respectively, the series 
(\ref{correl-fn})
satisfies the expansion of the system of differential equations
\begin{eqnarray}
&{\displaystyle \frac{\partial^{m}\varphi}{\partial z_{i}^{m}}+
\sum_{p=1}^{m}a_{p, \; i}(z_{1}, 
\dots, z_{n};q)
\frac{\partial^{m-p}\varphi}{\partial z_{i}^{m-p}}
=0,}\label{eqn1}\\
&{\displaystyle \prod_{k=1}^{m}\mathcal{O}_{i}\left(\sum_{i=1}^{n}
\wt w_{i}+2(m-k)\right)\varphi\quad\quad\quad\quad\quad\quad\quad
\quad\quad\quad\quad\quad\quad\quad\quad\quad}\nn
&{\displaystyle  
+\sum_{p=1}^{m}
b_{p, \;i}(z_{1}, \dots, z_{n}; q)
\prod_{k=1}^{m-p}\mathcal{O}_{i}\left(\sum_{i=1}^{n}
\wt w_{i}+2(m-p-k)\right)
\varphi=0,}\label{eqn2}
\end{eqnarray}
$i=1, \dots, n$, 
in the regions $1>|q_{ z_{1}}|>\cdots >|q_{z_{n}}|>|q|>0$.
\end{thm}
\pf
For fixed $w_{i}\in W_{i}$, $i=1, \dots, n$,
let $\Pi_{i}$ for $i=1, \dots, n$ be the $R$-submodules  of
$T/J$ generated by 
\begin{equation}\label{eqn1-0}
[1\otimes w_{1}\otimes \cdots\otimes w_{i-1}\otimes L^{k}(-1)w_{i}
\otimes w_{i+1}\otimes \cdots \otimes w_{n}]
\end{equation}
for $k\in \mathbb{N}$, respectively.
Since $R$ is Noetherian and 
$T/J$ is a finitely-generated $R$-module, 
the $R$-submodules $\Pi_{i}$ for $i=1, \dots, n$ are all 
finitely-generated. Thus there exist 
$a_{p, \;i}(z_{1}, \dots, z_{n}; q)\in R$
for $p=1, \dots, m$ and $i=1, \dots, n$  such that in $\Pi_{i}$
\begin{eqnarray}\label{eqn1-1}
\lefteqn{[1\otimes 
w_{1}\otimes \cdots\otimes w_{i-1}\otimes L^{m}(-1)w_{i}
\otimes w_{i+1}\otimes\cdots \otimes w_{n}]}\nn
&&+\sum_{p=1}^{m}a_{p, \; i}(z_{1}, 
\dots, z_{n};q)\cdot\nn
&&\quad\quad\quad\cdot
[1\otimes w_{1}\otimes \cdots\otimes w_{i-1}\otimes 
L^{m-p}(-1)w_{i}
\otimes w_{i+1}\otimes\cdots \otimes w_{n}]\nn
&&=0.
\end{eqnarray}
(Note that in the argument above, we actually first obtain 
$a_{p, \;i}(z_{1}, \dots, z_{n}; q)$ for $p=1, \dots, m$ 
for any fixed $i$. Thus $m$ in fact depends on $i$. But since there
are only finitely many $i$, we can always 
choose a sufficiently large $m$ such that it is independent of $i$.)
Since  for any $k\ge 0$, (\ref{eqn1-0})
has modular weight $\sum_{i=1}^{n}w_{i}+k$, we see from 
(\ref{eqn1-1}) that $a_{p, \; i}(z_{1}, 
\dots, z_{n};q)$ for $p=1, \dots, m$ and $i=1, \dots, n$
can be chosen to have modular weights $p$.
Applying $\psi_{\mathcal{Y}_{1}, \dots, \mathcal{Y}_{n}}$
to both sides of (\ref{eqn1-1}) and then
using the $L(-1)$-derivative property (Proposition \ref{l-1-der}) for 
(\ref{correl-fn}),
we see that  (\ref{correl-fn}) satisfies the equations (\ref{eqn1}). 

On the other hand, 
let $\mathcal{Q}_{i}: T\to T$ for $i=1, \dots, n$ be the linear 
map defined
by 
\begin{eqnarray*}
\lefteqn{\mathcal{Q}_{i}(1\otimes 
w_{1}\otimes \cdots \otimes w_{n})}\nn
&&=1\otimes w_{1}\otimes \cdots \otimes w_{i-1}\otimes L(-2)w_{i}
\otimes w_{i+1}\otimes \cdots \otimes w_{n}\nn
&&\quad-\sum_{k\ge 1}\tilde{G}_{2k+2}(q)
\otimes w_{1}\otimes \cdots \otimes w_{i-1}\otimes L(2k)w_{i}
\otimes w_{i+1}\otimes \cdots \otimes w_{n}\nn
&&\quad +\sum_{j\ne i}\sum_{m\ge 1}(-1)^{m}
\tilde{\wp}_{m+1}(z_{j}-z_{i}; q)\nn
&&\quad\quad\quad\quad\quad\otimes 
w_{1}\otimes \cdots \otimes w_{i-1}\otimes L(m-1)w_{i}
\otimes w_{i+1}\otimes \cdots \otimes w_{n}.
\end{eqnarray*}
For the same
fixed $w_{i}\in W_{i}$, $i=1, \dots, n$, as above, let $\Lambda_{i}$ 
for $i=1, \dots, n$ be the 
$R$-submodules  of $T/J$ generated by 
$[\mathcal{Q}_{i}^{k}(1\otimes w_{1}\otimes \cdots \otimes w_{n})]$ 
for $k\ge 0$. Since $R$ is Noetherian and 
$T/J$ is a finitely-generated $R$-module, 
the $R$-submodule $\Lambda_{i}$ is also finitely generated. 
Thus  there exist 
$b_{p, \;i}(z_{1}, \dots, z_{n}; q)
\in R$
for $p=1, \dots, m$ and for $i=1, \dots, n$ such that in $\Lambda_{i}$
\begin{equation}\label{eqn2-1}
[\mathcal{Q}_{i}^{m}(1\otimes w_{1}\otimes \cdots \otimes w_{n})]
+\sum_{p=1}^{t}b_{p, \;i}(z_{1}, \dots, z_{n}; q)
[\mathcal{Q}_{i}^{t-p}(1\otimes w_{1}\otimes \cdots \otimes w_{n})]=0.
\end{equation}
(Note that using the same argument above we can always find 
$m$ sufficiently large such that it is independent of $i$. But 
in general it might be different from the $m$ in (\ref{eqn1-1}). 
But we can always take the $m$ in (\ref{eqn1-1}) and the $m$ obtained 
here to be sufficiently large so that these two $m$ are equal.)
Since $[\mathcal{Q}_{i}^{k}(1\otimes w_{1}\otimes \cdots \otimes w_{n})]$
for any $k\ge 0$
has modular weight $\sum_{i=1}^{n}w_{i}+2k$, 
we see from 
(\ref{eqn2-1}) that $b_{p, \;i}(z_{1}, \dots, z_{n}; q)$
for $p=1, \dots, m$ and $i=1, \dots, n$
can be chosen to  have modular weights $2p$.
Applying $\psi_{\mathcal{Y}_{1}, \dots, \mathcal{Y}_{n}}$
to both sides of (\ref{eqn2-1}) and then
using (\ref{mod-inv-der}), 
we see that 
(\ref{correl-fn})
satisfies (\ref{eqn2}).
\epfv

\renewcommand{\theequation}{\thesection.\arabic{equation}}
\renewcommand{\thethm}{\thesection.\arabic{thm}}
\setcounter{equation}{0}
\setcounter{thm}{0}

\section{Chiral genus-one correlation functions and genus-one duality}

In this section, using the systems of differential equations 
obtained in the preceding section, 
we construct chiral genus-one correlation functions
and establish their duality properties.

We still fix a vertex operator algebra $V$ satisfying the same conditions
assumed in the preceding section.

\begin{thm}\label{conv-ext}
In the region $1>|q_{z_{1}}|>\cdots 
>|q_{z_{n}}|>|q_{\tau}|>0$, the series
\begin{equation}\label{correl-fn-tau}
F_{\mathcal{Y}_{1}, \dots, \mathcal{Y}_{n}}(w_{1}, 
\dots, w_{n};
z_{1}, \dots, z_{n}; q_{\tau})
\end{equation}
is absolutely convergent and can be analytically extended
to a (multivalued) analytic function in the region given by
$\Im(\tau)>0$ (here $\Im(\tau)$ is the imaginary part of $\tau$), 
$z_{i}\ne z_{j}+k\tau+l$ for $i, j=1, \dots, n$,
$i\ne j$, $k, l\in \mathbb{Z}$.
\end{thm}
\pf
We know that the coefficients of (\ref{correl-fn-tau}) 
as a series in powers of $q_{\tau}$ are absolutely 
convergent when $1>|q_{z_{1}}|>\cdots 
>|q_{z_{n}}|>0$. So for fixed $z_{1}, \dots, z_{n}$
satisfying $1>|q_{z_{1}}|>\cdots 
>|q_{z_{n}}|>0$, (\ref{correl-fn-tau}) is a well-defined 
series in powers of $q_{\tau}$. For fixed $z_{1}, \dots, z_{n}$ satisfying 
$1>|q_{z_{1}}|>\cdots 
>|q_{z_{n}}|>|q_{\tau}|>0$, the ordinary differential equation 
(\ref{eqn2}) with the variable $q_{\tau}$ has a regular singular
point at $q_{\tau}=0$. Since the series (\ref{correl-fn-tau})
satisfies (\ref{eqn2}) with $q=q_{\tau}$, 
it is absolutely convergent as a series 
in powers of $q_{\tau}$. Since the coefficients of the 
equation (\ref{eqn2}) are analytic in $z_{1}, \dots, z_{n}$,
the sum of (\ref{correl-fn-tau}) as a series 
in powers of $q_{\tau}$ is also analytic in $z_{1}, \dots, z_{n}$.
In particular, (\ref{correl-fn-tau}) as the expansion of 
an analytic function in the region $1>|q_{z_{1}}|>\cdots 
>|q_{z_{n}}|>|q_{\tau}|>0$ must be absolutely convergent as
a series of multiple sums.
Thus the sum of this series is analytic in $z_{1}, \dots, z_{n}$ and 
$q_{\tau}$ and give a (multivalued) 
analytic function in the region $1>|q_{ z_{1}}|>\cdots 
>|q_{z_{n}}|>|q_{\tau}|>0$. So the first part of the theorem is proved. 

Now we know that the coefficients of the system 
(\ref{eqn1})--(\ref{eqn2}) with $q=q_{\tau}$ 
are analytic in $z_{1}, \dots, z_{n}$ and $\tau$ 
with the only possible singularities 
$z_{i}\ne z_{j}+k\tau+l$ for $i, j=1, \dots, n$,
$i\ne j$, $k, l\in \mathbb{Z}$. So 
(\ref{correl-fn-tau})
as a solution of the system in the region $1>|q_{ z_{1}}|>\cdots 
>|q_{z_{n}}|>|q|>0$ can be analytically extended to the region 
given by
$\Im(\tau)>0$, 
$z_{i}\ne z_{j}+k\tau+l$ for $i, j=1, \dots, n$,
$i\ne j$, $k, l\in \mathbb{Z}$.
\epfv

We shall call  functions in the region 
$\Im(\tau)>0$,
$z_{i}\ne z_{j}+k\tau+l$ for $i, j=1, \dots, n$,
$i\ne j$, $k, l\in \mathbb{Z}$
obtained by analytically extending 
(\ref{correl-fn-tau})
{\it (chiral) genus-one correlation functions} or simply {\it genus-zero 
correlation functions} and we shall use 
\begin{equation}\label{correl-ext}
\overline{F}_{\mathcal{Y}_{1}, \dots, \mathcal{Y}_{n}}(w_{1}, 
\dots, w_{n};
z_{1}, \dots, z_{n}; \tau)
\end{equation}
to denote these functions.

\begin{thm}[Genus-one commutativity]
Let $W_{i}$ and $\tilde{W}_{i}$ 
be $V$-modules and
$\mathcal{Y}_{i}$ intertwining operators  of 
types $\binom{\tilde{W}_{i-1}}{W_{i}\tilde{W}_{i}}$ ($i=1, \dots, n$,
$\tilde{W}_{0}=\tilde{W}_{n}$), 
respectively. Then for any $1\le k\le n-1$, there exist
$V$-modules $\hat{W}_{k}$ and intertwining operators
$\hat{\mathcal{Y}}_{k}$ and $\hat{\mathcal{Y}}_{k+1}$
of types $\binom{\hat{W}_{k}}{W_{k}\tilde{W}_{k+1}}$
and $\binom{\tilde{W}_{k-1}}{W_{k+1}\hat{W}_{k}}$, respectively,
such that 
$$F_{\mathcal{Y}_{1}, \dots, \mathcal{Y}_{n}}(w_{1}, 
\dots, w_{n};
z_{1}, \dots, z_{n}; q_{\tau})$$
and 
\begin{eqnarray*}
\lefteqn{F_{\mathcal{Y}_{1}, \dots, \mathcal{Y}_{k-1},
\hat{\mathcal{Y}}_{k+1}, \hat{\mathcal{Y}}_{k}, \mathcal{Y}_{k+2}
\dots, \mathcal{Y}_{n}}(w_{1}, 
\dots, w_{k-1}, w_{k+1}, w_{k}, w_{k+2}, \dots, w_{n};}\nn
&&\quad\quad\quad\quad\quad\quad\quad\quad\quad\quad\quad\quad
\quad\quad\quad
z_{1}, \dots, z_{k-1}, z_{k+1}, z_{k}, z_{k+2}, 
\dots, z_{n}; q_{\tau})
\end{eqnarray*}
are analytic extensions of each other, or equivalently, 
\begin{eqnarray*}
\lefteqn{\overline{F}_{\mathcal{Y}_{1}, \dots, \mathcal{Y}_{n}}(w_{1}, 
\dots, w_{n};
z_{1}, \dots, z_{n}; \tau)}\nn
&&=\overline{F}_{\mathcal{Y}_{1}, \dots, \mathcal{Y}_{k-1},
\hat{\mathcal{Y}}_{k+1}, \hat{\mathcal{Y}}_{k}, \mathcal{Y}_{k+2}
\dots, \mathcal{Y}_{n}}(w_{1}, 
\dots, w_{k-1}, w_{k+1}, w_{k}, w_{k+2}, \dots, w_{n};\nn
&&\quad\quad\quad\quad\quad\quad\quad\quad\quad\quad\quad\quad
\quad\quad\quad
z_{1}, \dots, z_{k-1}, z_{k+1}, z_{k}, z_{k+2}, 
\dots, z_{n}; \tau).
\end{eqnarray*}
More generally, for any $\sigma\in S_{n}$, 
there exist $V$-modules $\hat{W}_{i}$ ($i=1, \dots, n$) and 
intertwining operators $\hat{\mathcal{Y}}_{i}$ of 
types $\binom{\hat{W}_{i-1}}{W_{\sigma(i)}\hat{W}_{i}}$ ($i=1, \dots, n$,
$\hat{W}_{0}=\hat{W}_{n}=\tilde{W}_{n}$), respectively, 
such that 
$$F_{\mathcal{Y}_{1}, \dots, \mathcal{Y}_{n}}(w_{1}, 
\dots, w_{n};
z_{1}, \dots, z_{n}; q_{\tau})$$
and 
$$F_{\hat{\mathcal{Y}}_{1}, \dots, \hat{\mathcal{Y}}_{n}}(w_{\sigma(1)}, 
\dots, w_{\sigma(n)};
z_{\sigma(1)}, \dots, z_{\sigma(n)}; q_{\tau})$$ 
are analytic extensions of each other, or equivalently,
\begin{eqnarray*}
\lefteqn{\overline{F}_{\mathcal{Y}_{1}, \dots, \mathcal{Y}_{n}}(w_{1}, 
\dots, w_{n};
z_{1}, \dots, z_{n}; \tau)}\nn
&&=\overline{F}_{\hat{\mathcal{Y}}_{1}, \dots, 
\hat{\mathcal{Y}}_{n}}(w_{\sigma(1)}, 
\dots, w_{\sigma(n)};
z_{\sigma(1)}, \dots, z_{\sigma(n)}; \tau).
\end{eqnarray*}
\end{thm}
\pf
This follows immediately from  commutativity for 
intertwining operators. 
\epfv

For any $V$-module $W$ and $r\in \mathbb{R}$, we use $P_{r}$
to denote the projection from $W$ or $\overline{W}$ to $W_{(r)}$.

\begin{thm}[Genus-one associativity]\label{g1-asso}
Let $W_{i}$ and $\tilde{W}_{i}$ for $i=1, \dots, n$ be 
$V$-modules and
$\mathcal{Y}_{i}$ intertwining operators  of 
types $\binom{\tilde{W}_{i-1}}{W_{i}\tilde{W}_{i}}$ ($i=1, \dots, n$,
$\tilde{W}_{0}=\tilde{W}_{n}$), 
respectively. Then for any $1\le k\le n-1$, 
there exist a $V$-module $\hat{W}_{k}$ and 
intertwining operators $\hat{\mathcal{Y}}_{k}$ and 
$\hat{\mathcal{Y}}_{k+1}$ of 
types $\binom{\hat{W}_{k}}{W_{k} W_{k+1}}$ and 
$\binom{\tilde{W}_{k-1}}{\hat{W}_{k}\tilde{W}_{k+1}}$, respectively, 
such that 
\begin{eqnarray}\label{g-1-iter}
\lefteqn{\overline{F}_{\mathcal{Y}_{1}, \dots, \mathcal{Y}_{k-1},
\hat{\mathcal{Y}}_{k+1}, \mathcal{Y}_{k+2}, \dots,
\mathcal{Y}_{n}}(w_{1}, 
\dots, w_{k-1}, \hat{\mathcal{Y}}(w_{k}, z_{k}-z_{k+1})
w_{k+1}, }\nn
&&\quad\quad\quad\quad\quad\quad\quad\quad\quad\quad\quad\quad\quad
w_{k+2}, \dots, w_{n};
z_{1}, \dots, z_{k-1}, z_{k+1}, \dots, z_{n}; \tau)\nn
&&=\sum_{r\in \mathbb{R}}
\overline{F}_{\mathcal{Y}_{1}, \dots, \mathcal{Y}_{k-1},
\hat{\mathcal{Y}}_{k+1}, \mathcal{Y}_{k+2}, \dots,
\mathcal{Y}_{n}}(w_{1}, 
\dots, w_{k-1}, P_{r}(\hat{\mathcal{Y}}(w_{k}, z_{k}-z_{k+1})
w_{k+1}), \nn
&&\quad\quad\quad\quad\quad\quad\quad\quad\quad\quad\quad\quad\quad
w_{k+2}, \dots, w_{n};
z_{1}, \dots, z_{k-1}, z_{k+1}, \dots, z_{n}; \tau)\nn
&&
\end{eqnarray}
is absolutely convergent when $1>|q_{z_{1}}|>\cdots 
>|q_{z_{k-1}}|>|q_{z_{k+1}}|>
\dots >|q_{z_{n}}|>|q_{\tau}|>0$ and 
$1>|q_{(z_{k}-z_{k+1})}-
1|>0$
and is convergent to 
$$\overline{F}_{\mathcal{Y}_{1}, \dots, \mathcal{Y}_{n}}(w_{1}, 
\dots, w_{n};
z_{1}, \dots, z_{n}; \tau)$$ 
when $1>|q_{ z_{1}}|>\cdots 
>|q_{z_{n}}|>|q_{\tau}|>0$ and 
$|q_{(z_{k}-z_{k+1})}|>1>|q_{(z_{k}-z_{k+1})}-
1|>0$.
\end{thm}
\pf
Using Proposition \ref{e-assoc}, we know that for any $1\le k\le n-1$, 
there exist a $V$-module $\tilde{W}_{k}$ and 
intertwining operators $\hat{\mathcal{Y}}_{k}$ and 
$\hat{\mathcal{Y}}_{k+1}$ of 
$\binom{\hat{W}_{k}}{W_{k} W_{k+1}}$ and 
$\binom{\tilde{W}_{k-1}}{\hat{W}_{k}\tilde{W}_{k+1}}$, respectively, 
such that for any $z_{1}, \dots, z_{n}\in \mathbb{C}$ satisfying 
$1>|q_{z_{1}}|>\dots >|q_{z_{n}}|>0$ and 
$|q_{z_{k+1}}|>|q_{z_{k}}-q_{z_{k+1}}|>0$, we have
\begin{eqnarray}\label{assoc-1}
\lefteqn{\langle \tilde{w}_{n}', 
\mathcal{Y}_{1}(\mathcal{U}(q_{z_{1}})w_{1}, q_{z_{1}})
\cdots \mathcal{Y}_{k-1}(\mathcal{U}(q_{z_{k-1}})w_{k-1}, 
q_{z_{k-1}})\cdot}\nn
&&\quad\quad\quad\quad\cdot\hat{\mathcal{Y}}_{k+1}(\mathcal{U}(q_{z_{k+1}})
\hat{\mathcal{Y}}_{k}(w_{k}, z_{k}-z_{k+1})
w_{k+1}, q_{z_{k+1}})\cdot\nn
&&\quad\quad\quad\quad\cdot
\mathcal{Y}_{k+2}(\mathcal{U}(q_{z_{k+2}})w_{k+2}, q_{z_{k+2}})\cdots
\mathcal{Y}_{n}(\mathcal{U}(q_{z_{n}})w_{n}, q_{z_{n}})\tilde{w}_{n}
\rangle\nn
&&=\langle \tilde{w}_{n}',
\mathcal{Y}_{1}(\mathcal{U}(q_{z_{1}})w_{1}, q_{z_{1}})
\cdots
\mathcal{Y}_{n}(\mathcal{U}(q_{z_{n}})w_{n}, q_{z_{n}})\tilde{w}_{n}
\rangle
\end{eqnarray}
for any $\tilde{w}_{n}\in \tilde{W}_{n}$ and $\tilde{w}'_{n}
\in \tilde{W}'_{n}$. Thus as  series in $q$,
\begin{eqnarray}\label{assoc-2}
\lefteqn{\tr_{\tilde{W}_{n}}
\mathcal{Y}_{1}(\mathcal{U}(q_{z_{1}})w_{1}, q_{z_{1}})
\cdots \mathcal{Y}_{k-1}(\mathcal{U}(q_{z_{k-1}})w_{k-1}, 
q_{z_{k-1}})\cdot}\nn
&&\quad\quad\quad\quad\cdot\hat{\mathcal{Y}}_{k+1}(\mathcal{U}(q_{z_{k+1}})
\hat{\mathcal{Y}}_{k}(w_{k}, z_{k}-z_{k+1})
w_{k+1}, q_{z_{k+1}})\cdot\nn
&&\quad\quad\quad\quad\cdot
\mathcal{Y}_{k+2}(\mathcal{U}(q_{z_{k+2}})w_{k+2}, q_{z_{k+2}})\cdots
\mathcal{Y}_{n}(\mathcal{U}(q_{z_{n}})w_{n}, q_{z_{n}})
q^{L(0)-\frac{c}{24}}\nn
&&=\tr_{\tilde{W}_{n}}
\mathcal{Y}_{1}(\mathcal{U}(q_{z_{1}})w_{1}, q_{z_{1}})
\cdots
\mathcal{Y}_{n}(\mathcal{U}(q_{z_{n}})w_{n}, q_{z_{n}})
q^{L(0)-\frac{c}{24}}.
\end{eqnarray}
Since the right-hand side of (\ref{assoc-2}) is convergent absolutely
when $q=q_{\tau}$ and $1>|q_{z_{1}}|>\cdots >|q_{z_{n}}|>|q_{\tau}|>0$,
the left-hand side is also convergent absolutely when 
$q=q_{\tau}$, $1>|q_{z_{1}}|>\cdots >|q_{z_{n}}|>|q_{\tau}|>0$
and $|q_{z_{k+1}}|>|q_{z_{k}}-q_{z_{k+1}}|>0$. 
Also, since the right-hand side of (\ref{assoc-2}) satisfies the system
(\ref{eqn1})--(\ref{eqn2}), so does the left-hand side. 
Thus the left-hand side of (\ref{assoc-2}) with $q=q_{\tau}$
is convergent absolutely to an
analytic function, in the region $1>|q_{z_{1}}|>\cdots >|q_{z_{n}}|
>|q_{\tau}|>0$
and $|q_{z_{k+1}}|>|q_{z_{k}}-q_{z_{k+1}}|>0$, which can be 
analytically extended to the multivalued function 
$$\overline{F}_{\mathcal{Y}_{1}, \dots, \mathcal{Y}_{n}}(w_{1}, 
\dots, w_{n};
z_{1}, \dots, z_{n}; \tau)$$ 
on 
$\mathbb{C}^{n}\times \mathbb{H}$. Moreover, in the region 
$1>|q_{z_{1}}|>\cdots 
>|q_{z_{k-1}}|>|q_{z_{k+1}}|>
\dots >|q_{z_{n}}|>|q_{\tau}|>0$ and 
$1>|q_{(z_{k}-z_{k+1})}-
1|>0$, this function has the expansion (\ref{g-1-iter}),
proving the theorem.
\epfv

\renewcommand{\theequation}{\thesection.\arabic{equation}}
\renewcommand{\thethm}{\thesection.\arabic{thm}}
\setcounter{equation}{0}
\setcounter{thm}{0}

\section{The regularity of the singular points for fixed $q=q_{\tau}$}

In this section, we prove that for fixed $\tau\in 
\mathbb{H}$, 
we can actually obtain the coefficients in the system
(\ref{eqn1}) with $q=q_{\tau}$ 
such that (\ref{eqn1}) with $q=q_{\tau}$  
is regular at  its  singular points of the form 
$z_{i}=z_{j}+\alpha \tau+\beta$ for $1\le i<j\le n$. In particular, 
at these singular points, all the 
chiral genus-one correlation functions as functions of 
$z_{1}, \dots,z_{n}$ are regular.
The method used here is an adaption of the method
used in \cite{H11}. Though it will be important for the further 
study of chiral genus-one correlation functions, 
the result in this section will not be needed
in the next two sections.

As in the genus-zero case discussed in \cite{H11}, 
we need  certain  filtrations on $R$
and on the $R$-module $T$.
For $m\in \mathbb{Z}_{+}+1$, let $F_{m}^{\mbox{\scriptsize \text{sing}}}(R)$
be the vector subspace of $R$ spanned by elements 
of the form 
$$f(q)\prod_{1\le i<j\le n} 
(\tilde{\wp}_{2}(z_{i}-z_{j}; q))^{k_{i, j}}
(\tilde{\wp}_{3}(z_{i}-z_{j}; q))^{l_{i, j}}$$ 
where $k_{i, j}, l_{i, j}\in \mathbb{Z}_{+}$ satisfying
$\sum_{1\le i<j\le n}2k_{i, j}+3l_{i, j}\le m$ 
and $f(q)$ is a polynomial
of $\tilde{G}_{4}(q)$ and $\tilde{G}_{6}(q)$.
Clearly these subspaces of $R$ give a filtration of 
$R$. With respect to this
filtration, $R$ is a filtered algebra, that is,
$F_{m_{1}}^{\mbox{\scriptsize \text{sing}}}(R)\subset 
F_{m_{2}}^{\mbox{\scriptsize \text{sing}}}(R)$ for $m_{1}\le m_{2}$,
$R=\cup_{m\in \mathbb{Z}}F_{m}^{\mbox{\scriptsize \text{sing}}}(R)$ and
$F_{m_{1}}^{\mbox{\scriptsize \text{sing}}}(R)F_{m_{2}}^{\mbox{\scriptsize 
\text{sing}}}(R)\subset
F_{m_{1}+m_{2}}^{\mbox{\scriptsize \text{sing}}}(R)$ for any $m_{1}, 
m_{2}\in \mathbb{Z}_{+}+1$. 

\begin{lemma}\label{wp-filtr}
For $m\in \mathbb{Z}_{+}+1$ and $i, j=1, \dots, n$ satisfying $i<j$, we have
$\tilde{\wp}_{m}(z_{i}-z_{j}; q)\in 
F_{m}^{\mbox{\scriptsize \text{sing}}}(R)$.
\end{lemma}
\pf
Using induction, this lemma follows easily from (\ref{wp-relation}),
the fact that 
$\tilde{\wp}_{m}(z_{i}-z_{j}; q)$ is a polynomial of 
$\tilde{\wp}_{2}(z_{i}-z_{j}; q)$ and
$\tilde{\wp}_{3}(z_{i}-z_{j}; q)$ and the fact that 
$\tilde{\wp}_{m}(z_{i}-z_{j}; q)$ is either even or odd in the variable
$z_{i}-z_{j}$.
\epfv

For convenience, we shall use $\sigma$ to denote 
$\sum_{i=1}^{n}\wt w_{i}$ for homogeneous
$w_{i}\in W_{i}$, $i=1, \dots, n$, when the dependence 
on $w_{i}$ is clear. 
Let $F_{r}^{\mbox{\scriptsize \text{sing}}}(T)$ for $r\in \mathbb{R}$ 
be the subspace of $T$
spanned by elements of the form 
$$\left(f(q)\prod_{1\le i<j\le n} 
(\tilde{\wp}_{2}(z_{i}-z_{j}; q))^{k_{i, j}}
(\tilde{\wp}_{3}(z_{i}-z_{j}; q))^{l_{i, j}}\right)
\otimes  w_{1}
\otimes \cdots \otimes w_{n}$$
where $k_{i, j}, l_{i, j}\in \mathbb{Z}_{+}$ satisfying
$\sum_{1\le i<j\le n}2k_{i, j}+3l_{i, j}+\sigma \le r$ 
and $f(q)$ is a polynomial
of $\tilde{G}_{4}(q)$ and $\tilde{G}_{6}(q)$.
These subspaces give a filtration of $T$ in the following sense:
$F_{r}^{\mbox{\scriptsize \text{sing}}}(T)\subset F_{s}^{\mbox{\scriptsize \text{sing}}}(T)$ for 
$r\le s$; 
$T=\cup_{r\in \mathbb{R}}F_{r}^{\mbox{\scriptsize \text{sing}}}(T)$;
$F_{m}^{\mbox{\scriptsize \text{sing}}}(R)F_{r}^{\mbox{\scriptsize \text{sing}}}(T)\subset
F_{r+m}^{\mbox{\scriptsize \text{sing}}}(T)$.

Let $F_{r}^{\mbox{\scriptsize \text{sing}}}(J)
=F_{r}^{\mbox{\scriptsize \text{sing}}}(T)\cap J$ 
for $r\in \mathbb{R}$. Then
we have the following:

\begin{prop}\label{filtr-zk=zl}
For any $r\in \mathbb{R}$, there exists $N\in \mathbb{R}$ such that
$F_{r}^{\mbox{\scriptsize \text{sing}}}(T)
\subset F_{r}^{\mbox{\scriptsize \text{sing}}}(J)+F_{N}(T)$.
\end{prop}
\pf
The proof is a refinement of the proof of Proposition 
\ref{decomposition}. The only additional property we need 
is that the elements  $\mathcal{A}_{j}(u, w_{1}, \dots, w_{n})$,
are all in 
$F^{\mbox{\scriptsize \text{sing}}}_{\swt u+\sigma}(J)$.
By Lemma \ref{wp-filtr}, this is true. 
\epfv

Let $R^{\mbox{\scriptsize \text{reg}}}$ be the commutative associative algebra 
over $\mathbb{C}$ generated by the series
$\tilde{G}_{4}(q)$, $\tilde{G}_{6}(q)$, 
$(z_{i}-z_{j})^{2}\tilde{\wp}_{2}(z_{i}-z_{j}; q)$,
$(z_{i}-z_{j})^{3}\tilde{\wp}_{3}(z_{i}-z_{j}; q)$ and 
$z_{i}-z_{j}$ for $i, j=1, \dots, n$ satisfying $i<j$. 
Note that $R^{\mbox{\scriptsize \text{reg}}}$ 
is a subalgebra of the algebra $R[z_{i}-z_{j}]_{1\le i<j\le n}$. 
Since $R^{\mbox{\scriptsize \text{reg}}}$ is finitely 
generated over the field $\mathbb{C}$, it is a Noetherian
ring.
We also consider the $R^{\mbox{\scriptsize \text{reg}}}$-module 
$$T^{\mbox{\scriptsize \text{reg}}}=R^{\mbox{\scriptsize \text{reg}}}\otimes 
W_{1}\otimes \cdots\otimes W_{n}.$$
Note that $T^{\mbox{\scriptsize \text{reg}}}$ is a subspace of 
$$R[z_{i}-z_{j}]_{1\le i<j\le n}\otimes 
W_{1}\otimes \cdots\otimes W_{n}.$$
The grading by conformal weights on $W_{1}\otimes \cdots\otimes W_{n}$
induces a grading (by conformal weights) on $T^{\mbox{\scriptsize \text{reg}}}$.
Let $T_{(r)}^{\mbox{\scriptsize \text{reg}}}$ for $r\in \mathbb{R}$ be the space 
of elements of $T^{\mbox{\scriptsize \text{reg}}}$ of conformal weight $r$. Then 
$T^{\mbox{\scriptsize \text{reg}}}=\coprod_{r\in \mathbb{R}}T_{(r)}^{\mbox{\scriptsize \text{reg}}}$.
 
Let $w_{i}\in W_{i}$ for $i=1, \dots, n$ be homogeneous. 
Then by 
Proposition \ref{filtr-zk=zl}, 
$$w_{1}\otimes \cdots\otimes w_{n}
=\mathcal{W}_{1}+\mathcal{W}_{2}$$
where $\mathcal{W}_{1}\in 
F_{\sigma}^{\mbox{\scriptsize \text{sing}}}(J)$
and $\mathcal{W}_{2}\in 
F_{N}(T)$ (and $\sigma=\sum_{i=1}^{n}\wt w_{i}$). 

\begin{lemma}\label{lemma1}
For any $s\in [0, 1)$, there exist
$S\in \mathbb{R}$ such that $s+S\in \mathbb{Z}_{+}$ and for any homogeneous
$w_{i}\in W_{i}$, $i=1, \dots, n$, satisfying $\sigma\in s+\mathbb{Z}$,
$\prod_{1\le i<j\le n}(z_{i}-z_{j})^{\sigma+S}
\mathcal{W}_{2}\in T^{\mbox{\scriptsize \text{reg}}}$.
\end{lemma}
\pf
Let $S$ be a real number such that $s+S\in \mathbb{Z}_{+}$
and such that for any $r\in \mathbb{R}$ satisfying
$r\le -S$, $T_{(r)}=0$. By definition, elements of 
$F_{r}^{\mbox{\scriptsize \text{sing}}}(T)$ for any $r\in \mathbb{R}$
are sums of elements of the form 
$$\left(f(q)\prod_{1\le i<j\le n} 
(\tilde{\wp}_{2}(z_{i}-z_{j}; q))^{k_{i, j}}
(\tilde{\wp}_{3}(z_{i}-z_{j}; q))^{l_{i, j}}\right)
\otimes  \tilde{w}_{1}
\otimes \cdots \otimes \tilde{w}_{n}$$
where $k_{i, j}, l_{i, j}\in \mathbb{Z}_{+}$ satisfying
\begin{equation}\label{ineq1}
\sum_{1\le i<j\le n}2k_{i, j}+3l_{i, j}
+ \sum_{i=1}^{n}\wt \tilde{w}_{i}\le r
\end{equation}
and $f(q)$ is a polynomial
of $\tilde{G}_{4}(q)$ and $\tilde{G}_{6}(q)$.
Since for nonzero $\tilde{w}_{1}
\otimes \cdots \otimes \tilde{w}_{n}$,
$\sum_{i=1}^{n}\wt \tilde{w}_{i}>-S$ which together with (\ref{ineq1})
implies 
$$r>\sum_{1\le i<j\le n}2k_{i, j}+3l_{i, j}-S$$ 
or 
$$r+S-\sum_{1\le i<j\le n}2k_{i, j}+3l_{i, j}>0.$$ 
Consequently, 
$$r+S-(2k_{i, j}+3l_{i, j})>0, \;\;\;1\le i<j\le n.$$
Thus 
$\prod_{1\le i<j\le n}(z_{i}-z_{j})^{r+S}
F_{r}^{\mbox{\scriptsize \text{sing}}}(T)\in 
T^{\mbox{\scriptsize \text{reg}}}$.

By definition,
$$\mathcal{W}_{2}=w_{1}\otimes\cdots \otimes w_{n}
-\mathcal{W}_{1},$$
where 
$$\mathcal{W}_{1}\in 
F_{\sigma}^{\mbox{\scriptsize \text{sing}}}(J)
\subset F_{\sigma}^{\mbox{\scriptsize \text{sing}}}(T).$$
By the discussion above, 
$\prod_{1\le i<j\le n}(z_{i}-z_{j})^{\sigma+S}
\mathcal{W}_{1}\in T^{\mbox{\scriptsize \text{reg}}}$ and by definition, 
$$w_{0}\otimes w_{1}\otimes w_{2}\otimes w_{3}\in T^{\mbox{\scriptsize \text{reg}}}.$$
Thus $\prod_{1\le i<j\le n}(z_{i}-z_{j})^{\sigma+S}
\mathcal{W}_{2}\in T^{\mbox{\scriptsize \text{reg}}}$.
\epfv

\begin{thm}\label{regular}
Let $W_{i}$ and $w_{i}\in W_{i}$ for $i=1, \dots, n$
be the same as in Theorem \ref{sys1} and let $\tau\in 
\mathbb{H}$. Then there exist
$$a_{p, \;j}(z_{1}, \dots, z_{n}; q)\in 
R_{p}$$
for $p=1, \dots, m$ and $j=1, \dots, n$
such that the (possible) singular points 
of the form $z_{i}=z_{j}+\alpha \tau+\beta$ 
for  $1\le i<j\le n$ and $\alpha, \beta\in \mathbb{Z}$
of the system (\ref{eqn1}) with $q=q_{\tau}$
satisfied by 
$$F_{\mathcal{Y}_{1}, \dots, \mathcal{Y}_{n}}(w_{1}, 
\dots, w_{n};
z_{1}, \dots, z_{n}; q_{\tau})$$
are regular. 
\end{thm}
\pf
We need only prove that for any fixed integers $i, j$ satisfying
$1\le i<j\le n\}$, the singular point $z_{i}=z_{j}$
is regular
because the coefficients of the system (\ref{eqn1})
are periodic with periods $1$ and $\tau$.

By Proposition \ref{filtr-zk=zl}, for $k=1, \dots, n$, 
$$1\otimes w_{1}\otimes \cdots \otimes w_{k-1}\otimes
L^{p}(-1)w_{k}\otimes w_{k+1}\otimes \cdots\otimes w_{n}
=\mathcal{W}_{1}^{(p)}+\mathcal{W}_{2}^{(p)}$$
for $p\ge 0$,
where $\mathcal{W}^{(p)}_{1}\in 
F_{\sigma+p}^{\mbox{\scriptsize \text{sing}}}(J)$
and $\mathcal{W}^{(p)}_{2}\in 
F_{N}(T)$. 

By  Lemma \ref{lemma1}, there exists $S\in \mathbb{R}$ such that 
$\sigma+S\in \mathbb{Z}_{+}$ and
$$\prod_{1\le i<j\le n}(z_{i}-z_{j})^{\sigma+p+S}
\mathcal{W}^{(p)}_{2}\in T^{\mbox{\scriptsize \text{reg}}}$$
and 
thus 
$$\prod_{1\le i<j\le n}(z_{i}-z_{j})^{\sigma+p+S}
\mathcal{W}^{(p)}_{2}\in \coprod_{r\le N}T_{(r)}^{\mbox{\scriptsize \text{reg}}}$$
for $p\ge 0$ and $1\le i<j\le n$. 
Since $R^{\mbox{\scriptsize \text{reg}}}$ is a Noetherian ring and 
$\coprod_{r\le N}T_{(r)}^{\mbox{\scriptsize \text{reg}}}$ is a finitely-generated 
$R^{\mbox{\scriptsize \text{reg}}}$-module, 
the submodule of $\coprod_{r\le N}T_{(r)}^{\mbox{\scriptsize \text{reg}}}$ generated by 
$\prod_{1\le i<j\le n}(z_{i}-z_{j})^{\sigma+p+S}
\mathcal{W}^{(p)}_{2}$ for $p\ge 0$ is also finitely generated. 
Let 
$$\prod_{1\le i<j\le n}(z_{i}-z_{j})^{\sigma+m-p+S}
\mathcal{W}^{(m-p)}_{2}$$
for $p=1, \dots, m$
be a set of generators of this submodule. (Note that as in the 
proof of Theorem \ref{sys1}, we can always choose $m$ to be 
independent of $i$.)
Then there exist $c_{p, \;l}(z_{1}, \dots, z_{n}; q)\in 
R^{\mbox{\scriptsize \text{reg}}}$
for $p=1, \dots, m$ such that 
\begin{eqnarray}\label{linear-comb1}
\lefteqn{\prod_{1\le i<j\le n}(z_{i}-z_{j})^{\sigma+m+S}
\mathcal{W}^{(m)}_{2}}\nn
&&=-\sum_{p=1}^{m}
c_{p, \;i}(z_{1}, \dots, z_{n}; q)
\prod_{1\le i<j\le n}(z_{i}-z_{j})^{\sigma+m-p+S}
\mathcal{W}^{(m-p)}_{2}.
\end{eqnarray}

Since $\mathcal{W}^{(m-p)}_{2}\in T$ for $p=1, \dots, m$, 
$$\prod_{1\le i<j\le n}(z_{i}-z_{j})^{\sigma+m-p+S}\mathcal{W}^{(m-p)}_{2}\in 
\prod_{1\le i<j\le n}(z_{i}-z_{j})^{\sigma+m-p+S}T.$$
Projecting both sides of 
(\ref{linear-comb1}) to $\prod_{1\le i<j\le n}(z_{i}-z_{j})^{\sigma+m+S}T$, 
we obtain
\begin{eqnarray}\label{linear-comb2}
\lefteqn{\prod_{1\le i<j\le n}(z_{i}-z_{j})^{\sigma+m+S}
\mathcal{W}^{(m)}_{2}}\nn
&&=-\sum_{p=1}^{m}
d_{p, \;i}(z_{1}, \dots, z_{n}; q)
\prod_{1\le i<j\le n}(z_{i}-z_{j})^{\sigma+m-p+S}
\mathcal{W}^{(m-p)}_{2},
\end{eqnarray}
where
$d_{p, \;i}(z_{1}, \dots, z_{n}; q)$ for $p=1, \dots, 
m$ are the projection images in 
$$\prod_{1\le i<j\le n}
(z_{i}-z_{j})^{p}R\cap R^{\mbox{\scriptsize \text{reg}}}$$
of $c_{p, \;i}(z_{1}, \dots, z_{n}; q)$. The equality (\ref{linear-comb1})
is equivalent to
\begin{equation}\label{linear-comb3}
\mathcal{W}^{(m)}_{2}+\sum_{p=1}^{m}
d_{p, \;i}(z_{1}, \dots, z_{n}; q)
\prod_{1\le i<j\le n}(z_{i}-z_{j})^{-p}
\mathcal{W}^{(m-p)}_{2}=0.
\end{equation}

Let 
\begin{equation}\label{coeff}
a_{p, \;i}(z_{1}, \dots, z_{n}; q)=
d_{p, \;i}(z_{1}, \dots, z_{n}; q)
\prod_{1\le i<j\le n}(z_{i}-z_{j})^{-p}\in R
\end{equation}
for $p=1, \dots, 
m$.
Then (\ref{linear-comb3}) gives
\begin{eqnarray}\label{w1}
\lefteqn{1\otimes w_{1}\otimes \cdots \otimes w_{i-1}
L^{m}(-1)w_{i}\otimes w_{i+1}\otimes \cdots\otimes w_{n}}\nn
&&+\sum_{p=1}^{m}a_{p, \;i}(z_{1}, \dots, z_{n}; q)\cdot\nn
&&\quad\quad\quad\quad\quad
\cdot (1\otimes w_{1}\otimes \cdots \otimes w_{i-1}\otimes 
L^{m-p}(-1)w_{i}\otimes w_{i+1}\otimes \cdots\otimes w_{n})\nn
&&\quad =\mathcal{W}^{(m)}_{1}+\sum_{p=1}^{m}
a_{p, \;i}(z_{1}, \dots, z_{n}; q)\mathcal{W}^{(m-p)}_{1}.
\end{eqnarray}
Since $\mathcal{W}^{(m-p)}_{1}\in 
F_{\sigma+p}^{\mbox{\scriptsize \text{sing}}}(J)
\subset J$ for $p=0, \dots, m$, the right-hand side of (\ref{w1})
is in $J$. Thus we obtain
\begin{eqnarray*}
\lefteqn{[1\otimes w_{1}\otimes \cdots \otimes w_{i-1}
L^{m}(-1)w_{i}\otimes w_{i+1}\otimes \cdots\otimes w_{n}]}\nn
&&+\sum_{p=1}^{m}a_{p, \;i}(z_{1}, \dots, z_{n}; q)\cdot\nn
&&\quad\quad\quad\quad
\cdot [1\otimes w_{1}\otimes \cdots \otimes w_{i-1}\otimes 
L^{m-p}(-1)w_{i}\otimes w_{i+1}\otimes \cdots\otimes w_{n}]=0.
\end{eqnarray*}

Now using the fact that $d_{p, \;i}(z_{1}, \dots, z_{n}; q)
\in R^{\mbox{\scriptsize \text{reg}}}$
for $p=1, 
\dots, m$, (\ref{coeff}),  we see
that the singular point $z_{i}=z_{j}$ for $1\le i< j\le n$ of 
the system (\ref{eqn1}) is regular.
As in the proof of Theorem \ref{sys1}, 
$a_{p, \;i}(z_{1}, \dots, z_{n}; q)$ can be further chosen 
to be also in $R_{p}$ for $p=1, \dots, m$.
\epfv

\section{Associative algebras and
vertex operator algebras}

In this section, we introduce a new product in a vertex operator
algebra and using this new product, we construct an associative algebra.
We study this associative algebra and its representations. These 
results are needed in the next section. 
The new product we introduce is very different from the one introduced
by Zhu in \cite{Z} and thus our algebra looks very different from Zhu's
algebra. It turns out that our algebra is in fact 
isomorphic to Zhu's algebra and the results for our algebra
are parallel to those results 
in \cite{Z}, \cite{FZ} and \cite{L2} for Zhu's
algebra. Thus many of the results in this section can also be proved
using the isomorphism between these two algebras and 
the results for Zhu's algebra. 

On the other hand, the product introduced in this section is conceptual
and geometric (see Remark \ref{bullet-geom} below) and many of the
formulas are simpler than those in \cite{Z}.  Also in the present paper,
we assume that the reader is familiar only with some basic notions and
results on vertex operator algebras, their representations and
intertwining operator algebras. Because of these reasons, we shall give
direct proofs of all results.

We define a product $\bullet$ in $V$ as follows: For $u, v\in V$,
\begin{eqnarray*}
u\bullet v&=&\res_{y}y^{-1}Y\left(u, \frac{1}{2\pi i}\log(1+y)\right)v\nn
&=&\res_{x}\frac{2\pi ie^{2\pi i x}}{e^{2\pi ix}-1}
Y(u, x)v.
\end{eqnarray*}
Let $\tilde{O}(V)$ be the subspace of
$V$ spanned by elements of the form 
\begin{eqnarray*}
\res_{y}y^{-n}Y\left(u, \frac{1}{2\pi i}\log(1+y)\right)v
&=&\res_{x}\frac{2\pi ie^{2\pi i x}}{(e^{2\pi ix}-1)^{n}}
Y(u, x)v
\end{eqnarray*}
for $n>1$ and $u, v\in V$. 
Let $\tilde{A}(V)=V/\tilde{O}(V)$.

\begin{prop}
The product $\bullet$ in $V$ induces a product (denoted still by $\bullet$)
in $\tilde{A}(V)$ such that $\tilde{A}(V)$ together with this product, 
the equivalence class
of the vacuum $\mathbf{1}\in V$ is an associative 
algebra with identity. Moreover,
for any $u, v\in V$, $(L(-1)u)\bullet v\in \tilde{O}(V)$ and
$u\bullet v\equiv v\bullet u-2\pi i u_{0}v \mod \tilde{O}(V)$. In particular, 
$\omega+\tilde{O}(V)$ is in the center of $\tilde{A}(V)$.
\end{prop}
\pf
Let $u_{1}, u_{2}, u_{3}\in V$ and $n\in \mathbb{Z}_{+}+1$. 
Using the commutator formula for vertex operators, we have
\begin{eqnarray}\label{algebra1}
\lefteqn{\res_{x_{1}}\res_{x_{2}}
\frac{2\pi ie^{2\pi i x_{1}}}{e^{2\pi ix_{1}}-1}
\frac{2\pi ie^{2\pi i x_{2}}}{(e^{2\pi ix_{2}}-1)^{n}}
Y(u_{1}, x_{1})Y(u_{2}, x_{2})u_{3}}\nn
&&=\res_{x_{1}}\res_{x_{2}}
\frac{2\pi ie^{2\pi i x_{1}}}{e^{2\pi ix_{1}}-1}
\frac{2\pi ie^{2\pi i x_{2}}}{(e^{2\pi ix_{2}}-1)^{n}}
Y(u_{2}, x_{2})Y(u_{1}, x_{1})u_{3}\nn
&&\quad +\res_{x_{1}}\res_{x_{2}}\res_{x_{0}}
\frac{2\pi ie^{2\pi i x_{1}}}{e^{2\pi ix_{1}}-1}
\frac{2\pi ie^{2\pi i x_{2}}}{(e^{2\pi ix_{2}}-1)^{n}}\cdot\nn
&&\quad\quad\quad\quad\quad\quad\cdot 
x_{2}^{-1}\delta\left(\frac{x_{1}-x_{0}}{x_{2}}\right)
Y(Y(u_{1}, x_{0})u_{2}, x_{2})u_{3}.
\end{eqnarray}
Since the left-hand side of (\ref{algebra1}) spans $V\bullet \tilde{O}(V)$
and the right-hand side of (\ref{algebra1})
is in $\tilde{O}(V)$, 
we conclude that $V\bullet \tilde{O}(V)\subset \tilde{O}(V)$. 

Now let $u_{1}, u_{2}, u_{3}\in V$ and $n\in \mathbb{Z}_{+}$. 
Using the Jacobi identity, we have 
\begin{eqnarray}\label{algebra2}
\lefteqn{\res_{x_{0}}\res_{x_{2}}
\frac{2\pi ie^{2\pi i x_{0}}}{(e^{2\pi ix_{0}}-1)^{n}}
\frac{2\pi ie^{2\pi i x_{2}}}{e^{2\pi ix_{2}}-1}
Y(Y(u_{1}, x_{0})u_{2}, x_{2})u_{3}}\nn
&&=\res_{x_{1}}\res_{x_{2}}\res_{x_{0}}
\frac{2\pi ie^{2\pi i x_{0}}}{(e^{2\pi ix_{0}}-1)^{n}}
\frac{2\pi ie^{2\pi i x_{2}}}{e^{2\pi ix_{2}}-1}\cdot\nn
&&\quad\quad\quad\quad\quad\quad\cdot 
x_{0}^{-1}\delta\left(\frac{x_{1}-x_{2}}{x_{0}}\right)
Y(u_{1}, x_{1})Y(u_{2}, x_{2})u_{3}\nn
&&\quad -\res_{x_{1}}\res_{x_{2}}\res_{x_{0}}
\frac{2\pi ie^{2\pi i x_{0}}}{(e^{2\pi ix_{0}}-1)^{n}}
\frac{2\pi ie^{2\pi i x_{2}}}{e^{2\pi ix_{2}}-1}\cdot\nn
&&\quad\quad\quad\quad\quad\quad\cdot 
x_{0}^{-1}\delta\left(\frac{x_{2}-x_{1}}{-x_{0}}\right)
Y(u_{2}, x_{2})Y(u_{1}, x_{1})u_{3}\nn
&&=\res_{x_{1}}\res_{x_{2}}\res_{x_{0}}
\frac{2\pi ie^{2\pi i (x_{1}-x_{2})}}{(e^{2\pi i(x_{1}-x_{2})}-1)^{n}}
\frac{2\pi ie^{2\pi i x_{2}}}{e^{2\pi ix_{2}}-1}\cdot\nn
&&\quad\quad\quad\quad\quad\quad\cdot 
x_{0}^{-1}\delta\left(\frac{x_{1}-x_{2}}{x_{0}}\right)
Y(u_{1}, x_{1})Y(u_{2}, x_{2})u_{3}\nn
&&\quad -\res_{x_{1}}\res_{x_{2}}\res_{x_{0}}
\frac{2\pi ie^{-2\pi i (x_{2}-x_{1})}}{(e^{-2\pi i(x_{2}-x_{1})}-1)^{n}}
\frac{2\pi ie^{2\pi i x_{2}}}{e^{2\pi ix_{2}}-1}\cdot\nn
&&\quad\quad\quad\quad\quad\quad\cdot 
x_{0}^{-1}\delta\left(\frac{x_{2}-x_{1}}{-x_{0}}\right)
Y(u_{2}, x_{2})Y(u_{1}, x_{1})u_{3}\nn
&&=\res_{x_{1}}\res_{x_{2}}
\frac{2\pi ie^{2\pi i (x_{1}-x_{2})}}{((e^{2\pi ix_{1}}-1)
e^{-2\pi ix_{2}}+(e^{-2\pi ix_{2}}-1))^{n}}
\frac{2\pi ie^{2\pi i x_{2}}}{e^{2\pi ix_{2}}-1}\cdot\nn
&&\quad\quad\quad\quad\quad\quad\cdot 
Y(u_{1}, x_{1})Y(u_{2}, x_{2})u_{3}\nn
&&\quad -\res_{x_{1}}\res_{x_{2}}
\frac{2\pi ie^{-2\pi i (x_{2}-x_{1})}}{(e^{-2\pi ix_{2}}-1)
e^{2\pi ix_{1}}-(e^{2\pi ix_{1}}-1))^{n}}
\frac{2\pi ie^{2\pi i x_{2}}}{e^{2\pi ix_{2}}-1}\cdot\nn
&&\quad\quad\quad\quad\quad\quad\cdot 
Y(u_{2}, x_{2})Y(u_{1}, x_{1})u_{3}\nn
&&=\sum_{k\in \mathbb{N}}\binom{-n}{k}\res_{x_{1}}\res_{x_{2}}
2\pi ie^{2\pi i (x_{1}-x_{2})}(e^{2\pi ix_{1}}-1)^{-n-k}
e^{-2\pi i(-n-k)x_{2}}\cdot\nn
&&\quad\quad\quad\quad\quad\quad\cdot (e^{-2\pi ix_{2}}-1)^{k}
\frac{2\pi ie^{2\pi i x_{2}}}{e^{2\pi ix_{2}}-1}
Y(u_{1}, x_{1})Y(u_{2}, x_{2})u_{3}\nn
&&\quad -\sum_{k\in \mathbb{N}}\binom{-n}{k}\res_{x_{1}}\res_{x_{2}}
2\pi ie^{-2\pi i (x_{2}-x_{1})}(e^{-2\pi ix_{2}}-1)^{-n-k}
e^{2\pi i(-n-k)x_{1}}\cdot\nn
&&\quad\quad\quad\quad\quad\quad\cdot (e^{2\pi ix_{1}}-1)^{k}
\frac{2\pi ie^{2\pi i x_{2}}}{e^{2\pi ix_{2}}-1}
Y(u_{2}, x_{2})Y(u_{1}, x_{1})u_{3}.
\end{eqnarray}
When $n>1$, the left-hand side of (\ref{algebra2})
spans $\tilde{O}(V)\bullet V$ and
the right-hand side of (\ref{algebra2}) 
is in $\tilde{O}(V)$. So $\tilde{O}(V)\bullet V
\subset \tilde{O}(V)$. Thus $\tilde{O}(V)$ is a two-side
ideal for the product $\bullet$ so that $\bullet$ induces 
a product, denoted still by $\bullet$, in $\tilde{A}(V)$.

When $n=1$, the left-hand side of (\ref{algebra2})
is $(u_{1}\bullet u_{2})\bullet u_{3}$ and the right-hand side
of (\ref{algebra2}) is a sum of $u_{1}\bullet (u_{2}\bullet u_{3})$
and elements of $\tilde{O}(V)$. Thus the product $\bullet$ in 
$\tilde{A}(V)$ is associative. 

It is obvious that $\mathbf{1}+\tilde{O}(V)$ is the identity 
of $\tilde{A}(V)$. 

For $u, v\in V$, using the $L(-1)$-derivative property and 
the property of $\res_{x}$, we have
\begin{eqnarray*}
(L(-1)u)\bullet v&=&\res_{x}\frac{2\pi ie^{2\pi ix}}{e^{2\pi ix}-1}
Y(L(-1)u, x)v\nn
&=&\res_{x}\frac{2\pi ie^{2\pi ix}}{e^{2\pi ix}-1}
\frac{d}{dx}Y(u, x)v\nn
&=&-\res_{x}\frac{d}{dx}\frac{2\pi ie^{2\pi ix}}{e^{2\pi ix}-1}
Y(u, x)v\nn
&=&\res_{x}\frac{(2\pi i)^{2}e^{2\pi ix}}{(e^{2\pi ix}-1)^{2}}Y(u, x)v\nn
&\in & \tilde{O}(V).
\end{eqnarray*}
In particular, we have $L(-1)u=(L(-1)u)\bullet \mathbf{1}\in \tilde{O}(V)$.

For $u, v\in V$, using this fact and skew-symmetry, we have
\begin{eqnarray*}
u\bullet v&=&\res_{x}\frac{2\pi ie^{2\pi ix}}{e^{2\pi ix}-1}
Y(u, x)v\nn
&=&\res_{x}\frac{2\pi ie^{2\pi ix}}{e^{2\pi ix}-1}
e^{xL(-1)}Y(v, -x)u\nn
&\equiv&\res_{x}\frac{2\pi ie^{2\pi ix}}{e^{2\pi ix}-1}
Y(v, -x)u\mod \tilde{O}(V)\nn
&=&\res_{y}\frac{-2\pi ie^{-2\pi iy}}{e^{-2\pi iy}-1}Y(v, y)u\nn
&=&\res_{y}\left(\frac{2\pi ie^{2\pi iy}}{e^{2\pi iy}-1}-2\pi i\right)
Y(v, y)u\nn
&=&v\bullet u-2\pi i u_{0}v.
\end{eqnarray*}
Since $\omega_{0}=L(-1)$, we have
\begin{eqnarray*}
\omega\bullet u&\equiv&u\bullet \omega-2\pi i \omega_{0}u\mod \tilde{O}(V)\nn
&\equiv&u\bullet \omega -2\pi i L(-1)u\mod \tilde{O}(V)\nn
&\equiv&u\bullet \omega \mod  \tilde{O}(V).
\end{eqnarray*}
\epfv

\begin{rema}\label{bullet-geom}
{\rm The product $\bullet$ has a very clear geometric meaning.
We know that vertex operators correspond to the sphere 
$\mathbb{C}\cup \{\infty\}$ with the negatively oriented puncture 
$\infty$ and the  positively oriented and ordered puncture $z$ and $0$ and 
with the standard locally coordinates vanishing at these punctures
(see \cite{H3.5}). 
The variable in the vertex operators corresponds to the position of the 
first oriented puncture $z$. We can restrict $z$ to be in some subsets of 
$\mathbb{C}\setminus \{0\}$ to get vertex operators 
defined only locally. Since Riemann surfaces are constructed
locally using subsets of $\mathbb{C}$, 
we have vertex operators locally on any Riemann surfaces (see, for example,
\cite{FB}). 
The map $y\mapsto \frac{1}{2\pi i}\log(1+y)$ in fact maps 
an annulus in the sphere $\mathbb{C}\cup \{\infty\}$ to 
a parallelogram in the universal covering of 
the torus corresponding to the parallelogram. 
So the product $\bullet$ can be understood as 
the constant term of the pull-back of the local vertex operator map
on the torus by this map to the annulus. Such pull-
backs of local vertex operators are needed
because we want to construct (global) correlation 
functions on one Riemann surface (a torus) 
from the (global) correlation functions 
on another Riemann surface (the sphere), not just want to study 
local vertex operators on a single 
Riemann surface or sheaves of vertex operators obtained from 
local vertex operators.}
\end{rema}

\begin{prop}
The associative algebra $\tilde{A}(V)$ is isomorphic 
to the associative algebra $A(V)$ constructed by Zhu in 
\cite{Z}.
\end{prop}
\pf
For $u, v\in V$ $n\in \mathbb{Z}_{+}$, by (\ref{chg-var0}), we have
\begin{eqnarray*}
\lefteqn{\mathcal{U}(1)\res_{y}y^{-n}
Y\left(u, \frac{1}{2\pi i}\log(1+y)\right)v}\nn
&&=\res_{y}y^{-n}
Y\left((1+y)^{L(0)}\mathcal{U}(1)u, y\right)\mathcal{U}(1)v.
\end{eqnarray*}
By this formula in the cases of $n\ge 2$ 
and the fact that $\mathcal{U}(1)$ is
invertible, we have
$$\mathcal{U}(1)(\tilde{O}(V))=O(V).$$
Thus $\mathcal{U}(1)$ induces a linear isomorphism 
from $\tilde{A}(V)$ to $A(V)$. 

For $u, v\in V$, using the formula above in the case of $n=1$, we have
\begin{eqnarray*}
\mathcal{U}(1)(u\bullet v)&=&\res_{y}y^{-1}
\mathcal{U}(1)Y\left(u, \frac{1}{2\pi i}\log(1+y)\right)v\nn
&=&\res_{y}y^{-1}
Y\left((1+y)^{L(0)}\mathcal{U}(1)u, y\right)\mathcal{U}(1)v\nn
&=&(\mathcal{U}(1)u)*(\mathcal{U}(1)v),
\end{eqnarray*}
which
shows that the linear isomorphism induced by 
$\mathcal{U}(1)$ is in fact an isomorphism of algebras.
\epfv

We need the notion of $\mathbb{N}$-gradable weak $V$-module (which was
called $V$-module in \cite{Z}). A {\it weak $V$-module} is a 
vector space $W$ equipped with a vertex operator map 
\begin{eqnarray*}
Y: V\otimes W&\to& W[[x, x^{-1}]],\nn
u\otimes w&\mapsto & Y(u, x)w=\sum_{n\in \mathbb{Z}}u_{n}wx^{-n-1}
\end{eqnarray*}
satisfying all the axioms for 
$V$-modules except for those involving the grading. 
An {\it $\mathbb{N}$-gradable weak $V$-module} is a weak $V$-module 
such that there exists
an $\mathbb{N}$-grading $W=\coprod_{n\in \mathbb{N}}W_{(n)}$
satisfying the condition that $u_{n}w\in W_{(\swt u -n-1 +\swt w)}$
for homogeneous $u\in V$, $w\in W$ and $n\in \mathbb{Z}$.

Let $W$ be an $\mathbb{N}$-gradable weak $V$-module and 
$$T(W)=\{w\in W\;|\; u_{n}w=0, u\in V, \wt u-n-1<0\}.$$
Let $P_{T(W)}: W\to T(W)$ be the projection from 
$W$ to $T(W)$. 
For any $u\in V$, we define
$\rho_{W}(u): T(W)\to T(W)$ by 
\begin{eqnarray*}
\rho_{W}(u)w&=&P_{T(W)}(o(\mathcal{U}(1)u)w)\nn
&=&\res_{x}x^{-1}P_{T(W)}(Y(\mathcal{U}(x)u, x)w)
\end{eqnarray*}
for $w\in T(W)$.

In remaining part of this section, we assume 
for simplicity that  $V_{(n)}=0$ for $n< 0$.

\begin{prop}\label{functor-f}
For any $\mathbb{N}$-gradable weak $V$-module $W$, $\rho_{W}(u)=0$ for
$u\in \tilde{O}(V)$ and the map given by $u+\tilde{O}(V)\mapsto
\rho_{W}(u)$ for $u\in V$ gives $T(W)$ an $\tilde{A}(V)$-module
structure. The functor $T$ from the category of $\mathbb{N}$-gradable weak
$V$-modules to the category of $\tilde{A}(V)$-modules given by $W\mapsto
T(W)$ has a right inverse, that is, there exists a functor $S$ from the
category of $\tilde{A}(V)$-modules to the category of 
$\mathbb{N}$-gradable weak $V$-modules such that $TS=I$, where 
$I$ is the identity functor on the category of $\tilde{A}(V)$-modules.
In particular, for any 
$\tilde{A}(V)$-module $M$, $T(S(M))=M$. Moreover we can find such an
$S$ such that for any 
$\mathbb{N}$-gradable weak
$V$-module $W$, there exists a natural surjective homomorphism of $V$-modules
from $S(T(W))$ to the $\mathbb{N}$-gradable weak $V$-submodule of $W$
generated by $T(W)$.  
\end{prop}
\pf
For $n\in \mathbb{Z}_{+}$, $u, v\in V$ and $w\in T(W)$,  
using the definitions, (\ref{chg-var0}), the 
$L(0)$-conjugation formula, the Jacobi identity and the assumption
that the weight of a nonzero element of $V$ is nonnegative, we have 
\begin{eqnarray}\label{module1}
\lefteqn{\rho\left(\res_{x_{0}}
\frac{2\pi ie^{2\pi ix_{0}}}{(e^{2\pi ix_{0}}-1)^{n}}Y(u, x_{0})v\right)w}\nn
&&=\res_{x_{2}}\frac{1}{x_{2}}\res_{x_{0}}
\frac{2\pi ie^{2\pi ix_{0}}}{(e^{2\pi ix_{0}}-1)^{n}}
P_{T(W)}(Y(\mathcal{U}(x_{2})Y(u, x_{0})v, x_{2})w)\nn
&&=\res_{x_{2}}\res_{x_{0}}
\frac{2\pi ie^{2\pi ix_{0}}}{x_{2}(e^{2\pi ix_{0}}-1)^{n}}\cdot \nn
&&\quad\quad\quad\cdot 
P_{T(W)}(Y(x_{2}^{L(0)}Y(\mathcal{U}(e^{2\pi ix_{0}})u, e^{2\pi ix_{0}}-1)
\mathcal{U}(1)v, x_{2})w)\nn
&&=\res_{x_{2}}\res_{x_{0}}
\frac{2\pi ie^{2\pi ix_{0}}}{x_{2}(e^{2\pi ix_{0}}-1)^{n}}\cdot \nn
&&\quad\quad\quad\cdot P_{T(W)}(Y(Y(x_{2}^{L(0)}\mathcal{U}(e^{2\pi ix_{0}})u, 
x_{2}(e^{2\pi ix_{0}}-1))
x_{2}^{L(0)}\mathcal{U}(1)v, x_{2})w)\nn
&&=\res_{x_{2}}\res_{y_{0}}
\frac{1}{x_{2}y_{0}^{n}}
P_{T(W)}(Y(Y(\mathcal{U}(x_{2}+y_{0})u, y_{0})
\mathcal{U}(x_{2})v, x_{2})w)\nn
&&=\res_{x_{2}}\res_{y_{0}}\res_{x_{1}}
\frac{1}{x_{2}y_{0}^{n}}y_{0}^{-1}\delta\left(\frac{x_{1}-x_{2}}{y_{0}}\right)
\cdot \nn
&&\quad\quad\quad\cdot P_{T(W)}(Y(\mathcal{U}(x_{2}+y_{0})u, x_{1})
Y(\mathcal{U}(x_{2})v, x_{2})w)\nn
&&\quad -\res_{x_{2}}\res_{y_{0}}\res_{x_{1}}
\frac{1}{x_{2}y_{0}^{n}}y_{0}^{-1}\delta\left(\frac{x_{2}-x_{1}}{-y_{0}}\right)
\cdot \nn
&&\quad\quad\quad\cdot P_{T(W)}(Y(\mathcal{U}(x_{2})v, x_{2})
Y(\mathcal{U}(x_{2}+y_{0})u, x_{1})w)\nn
&&=\res_{x_{2}}\res_{x_{1}}\frac{1}{x_{2}(x_{1}-x_{2})^{n}}
P_{T(W)}(Y(\mathcal{U}(x_{1})u, x_{1})
Y(\mathcal{U}(x_{2})v, x_{2})w)\nn
&&\quad -\res_{x_{2}}\res_{x_{1}}
\frac{1}{x_{2}(-x_{2}+x_{1})^{n}}
P_{T(W)}(Y(\mathcal{U}(x_{2})v, x_{2})
Y(\mathcal{U}(x_{1})u, x_{1})w)\nn
&&=\res_{x_{2}}\res_{x_{1}}\frac{1}{x_{2}(x_{1}-x_{2})^{n}}
P_{T(W)}(Y(\mathcal{U}(x_{1})u, x_{1})
Y(\mathcal{U}(x_{2})v, x_{2})w).\nn
&&
\end{eqnarray}

When $n\ge 2$, the left-hand  side of (\ref{module1})
spans $\rho(\tilde{O}(V))$ and the right-side of (\ref{module1})
is equal to $0$. 
When $n=1$, the left-hand  and right-hand sides
of (\ref{module1}) are equal to 
$\rho(u\bullet v)w$ and $\rho(u)\rho(v)w$, respectively. Thus the first 
conclusion is proved. 

To prove the second conclusion, we need to construct
an $\mathbb{N}$-gradable weak $V$-module from an $\tilde{A}(V)$-module.
Consider the affinization
$$V[t,t^{-1}] = V \otimes \mathbb{C}[t,t^{-1}]
\subset V \otimes \mathbb{C} ((t)) \subset V\otimes \mathbb{C} [[t,t^{-1}]]
\subset V[[t,t^{-1}]]$$
of $V$ and the tensor algebra $\mathcal{T}(V[t,t^{-1}])$ 
generated by $V[t,t^{-1}]$. For simplicity,
we shall denote $u\otimes t^{m}$ for $u\in V$ and $m\in \mathbb{Z}$ by
$u(m)$ and we shall omit the tensor product sign $\otimes$ when we 
write an element of $\mathcal{T}(V[t,t^{-1}])$. Thus $\mathcal{T}(V[t,t^{-1}])$
is spanned by elements of the form $u_{1}(m_{1})\cdots u_{k}(m_{k})$
for $u_{i}\in V$ and $m_{i}\in \mathbb{Z}$, $i=1, \dots, k$.

Let $M$ be an $\tilde{A}(V)$-module and let $\rho: \tilde{A}(V)
\to \edo M$ be the map giving the representation of $\tilde{A}(V)$
on $M$. Consider 
$\mathcal{T}(V[t,t^{-1}])\otimes M$. Again for simplicity we shall omit the 
tensor product sign. So $\mathcal{T}(V[t,t^{-1}])\otimes M$
is spanned by elements
of the form $u_{1}(m_{1})\cdots u_{k}(m_{k})w$ for
$u_{i}\in V$, $m_{i}\in \mathbb{Z}$, $i=1, \dots, k$, and $w\in M$
and for any $u\in V$, $m\in \mathbb{Z}$, $u(m)$ acts from the left on 
$\mathcal{T}(V[t,t^{-1}])\otimes M$.
For homogeneous $u_{i}\in V$, $m_{i}\in \mathbb{Z}$, $i=1, \dots, k$, 
and $w\in M$, 
we define the {\it degree} of $u_{1}(m_{1})\cdots u_{k}(m_{k})w$
to be $(\wt u_{1}-m_{1}-1)+\cdots +(\wt u_{k}-m_{k}-1)$. 
For any $u\in V$, let
$$Y_{t}(u, x): \mathcal{T}(V[t,t^{-1}])\otimes M
\to (\mathcal{T}(V[t,t^{-1}])\otimes M)[[x, x^{-1}]]$$
be defined by 
$$Y_{t}(u, x)=\sum_{m\in \mathbb{Z}}u(m)x^{-m-1}.$$
For a homogeneous element $u\in V$, let $o_{t}(u)=u(\wt u-1)$. 
Using linearity, we extend $o_{t}(u)$ to nonhomogeneous $u$.

Let $\mathcal{I}$ be the $\mathbb{Z}$-graded 
$\mathcal{T}(V[t,t^{-1}])$-submodule 
of $\mathcal{T}(V[t,t^{-1}])\otimes M$
generated by elements of the forms
$u(m)w$ ($u\in V$, $\wt u-m-1< 0$, $w\in M$), 
$o_{t}(\mathcal{U}(1)u)w-\rho(u+\tilde{O}(V))w$ ($u\in V$, $w\in M$) and
the coefficients in $x_{1}$ and $x_{2}$ of 
\begin{eqnarray*}
\lefteqn{Y_{t}(u, x_{1})Y_{t}(v, x_{2})w
-Y_{t}(v, x_{2})Y_{t}(u, x_{1})w}&\nn
&&-\res_{x_{0}}x_{2}^{-1}\delta\left(\frac{x_{1}-x_{0}}{x_{2}}\right)
Y_{t}(Y(u, x_{0})v, x_{2})w
\end{eqnarray*}
($u, v\in V$ and $w\in \mathcal{T}(V[t,t^{-1}])\otimes M$). 
(Note that the coefficients of the formal expression above 
are indeed in $\mathcal{T}(V[t,t^{-1}])\otimes M$.)
Let $S_{1}(M)=(\mathcal{T}(V[t,t^{-1}])\otimes M)/\mathcal{I}$. Then 
$S_{1}(M)$ is also  a $\mathbb{Z}$-graded 
$\mathcal{T}(V[t,t^{-1}])$-module. 
In fact, by definition of $\mathcal{I}$, we see that 
$S_{1}(M)$ is  spanned by elements
of the form $u_{1}(m_{1})\cdots u_{k}(m_{k})w+\mathcal{I}$ for homogeneous
$u_{i}\in V$, $m_{i}<\wt u_{i}-1$, $i=1, \dots, k$ and $w\in M$. In particular, we 
see that $S_{1}(M)$ has an $\mathbb{N}$-grading.
Thus for $u\in V$ and $w\in S_{1}(M)$,
$u(m)w=0$ when $m$ is sufficiently large.

Let $\mathcal{J}$ be the $\mathbb{N}$-graded 
$\mathcal{T}(V[t,t^{-1}])$-submodule of $S_{1}(M)$
generated by the coefficients in $x$ of  
$$Y_{t}(L(-1)u, x)w-\frac{d}{dx}Y_{t}(u, x)w$$
and the coefficients in $x_{0}$ and $x_{2}$ of 
\begin{eqnarray*}
&{\displaystyle Y_{t}(Y(u, x_{0})v, x_{2})w-\res_{x_{1}}
x_{0}^{-1}\delta\left(\frac{x_{1}-x_{2}}{x_{0}}\right)
Y_{t}(u, x_{1})Y_{t}(v, x_{2})w}&\nn
&{\displaystyle +x_{0}^{-1}\delta\left(\frac{x_{2}-x_{1}}{-x_{0}}\right)
Y_{t}(v, x_{2})Y_{t}(u, x_{1})w}&
\end{eqnarray*}
(which are indeed in $S_{1}(M)$)
for $u, v\in V$, $w\in S_{1}(M)$.

Let $S(M)=S_{1}(M)/\mathcal{J}$. Then $S(M)$ is also 
an $\mathbb{N}$-graded $\mathcal{T}(V[t,t^{-1}])$-module. We can
still use elements of $\mathcal{T}(V[t,t^{-1}])\otimes M$ 
to represent elements of $S(M)$. But note that these 
elements now satisfy relations. We equip $S(M)$ with 
the vertex operator map $Y: V\otimes S(M)\to S(M)[[x, x^{-1}]]$
given by $u\otimes w\mapsto Y(u, x)w=Y_{t}(u, x)w$. 
As in $S_{1}(M)$, for $u\in V$ and $w\in S(M)$,
we also have $u(m)w=0$ when $m>\wt u-1$. Clearly 
$Y(\mathbf{1}, x)=I_{S(M)}$ (where $I_{S(M)}$ is the identity 
operator on $S(M)$). 
By definition, we know that the commutator formula, the associator 
formula and the $L(-1)$-derivative property all hold. 
Thus $S(M)$ is an $\mathbb{N}$-gradable weak $V$-module 
such that $T(S(M))=M$. 

Let $W$ be an $\mathbb{N}$-gradable weak $V$-module. 
We define a linear map from $S(T(W))$ to $W$ by 
mapping $u_{1}(m_{1})\cdots u_{k}(m_{k})w$ of $S(M)$
to $(u_{1})_{m_{1}}\cdots (u_{k})_{m_{k}}w$ of $W$ for 
$u_{i}\in V$, $m_{i}\in \mathbb{Z}$ ($i=1, \dots, k$) and $w\in T(W)$.
Note that the relations among $u_{1}(m_{1})\cdots u_{k}(m_{k})w$
for $u_{i}\in V$, $m_{i}\in \mathbb{Z}$ ($i=1, \dots, k$) and $w\in T(W)$
are given by just the action of $V$ on $T(W)$, the commutator formula,
the associator formula and the $L(-1)$-derivative property 
for vertex operators. These relations also hold in $W$. Thus 
this map is well-defined. Clearly, this is a surjective homomorphism of 
$V$-modules from $S(T(W))$ to 
the $\mathbb{N}$-gradable weak $V$-submodule of $W$
generated by $T(W)$.
\epfv

\begin{prop}\label{reductive}
Assume that every $\mathbb{N}$-gradable weak $V$-module is 
completely reducible. Then $T$ and $S$ are equivalences of categories.
In particular, $\tilde{A}(V)$ is semi-simple and
$M$ is an irreducible $\tilde{A}(V)$-module if and only if $S(M)$ is an 
irreducible $V$-module.
\end{prop}
\pf
We need only prove that for any $\mathbb{N}$-gradable weak $V$-module
$W$, $S(T(W))$ is naturally isomorphic to $W$. Since 
every $\mathbb{N}$-gradable weak $V$-module is 
completely reducible, we need only to consider the case that 
$W$ is irreducible.

Let $W$ be an irreducible $\mathbb{N}$-gradable weak $V$-module.
Then by assumption, $S(T(W))$ is completely reducible. 
If $S(T(W))$ is not irreducible, then there exist nonzero
$V$-modules $W_{1}$ and
$W_{2}$ such that $S(T(W))=W_{1}\oplus W_{2}$. In particular, 
we have $T(W)=T(S(T(W)))=T(W_{1})\oplus T(W_{2})$. Here 
$T(W_{1})$ and $T(W_{2})$ are both nonzero because both 
$W_{1}$ and $W_{2}$ are nonzero.
we have a $V$-submodule $\tilde{W}_{1}$ of $W$ generated by $T(W_{1})$. 
This $V$-submodule $\tilde{W}_{1}$ is obviously nonzero because
$T(W_{1})$ is nonzero. It is also not $W$ since $T(W_{2})$ is not
in $\tilde{W}_{1}$. Thus $W$ is not irreducible. Contradiction. 
So $S(T(W))$ is irreducible. Since $W$ is irreducible, 
the natural homomorphism in Proposition \ref{functor-f}
from $S(T(W))$ to $W$ is surjective. This homomorphism 
must also be injective because $S(T(W))$ is irreducible. Thus
$S(T(W))$ is naturally isomorphic to $W$.
\epfv

\begin{cor}\label{finitely-many}
If every $\mathbb{N}$-gradable weak $V$-module is 
completely reducible, then there are only finitely many 
inequivalent irreducible $V$-modules.
\end{cor}
\pf
By Proposition \ref{reductive}, $\tilde{A}(V)$ is semisimple.
Thus there are only finitely many inequivalent irreducible 
$\tilde{A}(V)$-modules. By Proposition \ref{reductive} again, 
these finitely many irreducible $\tilde{A}(V)$-modules
are mapped under $S$ to a complete set of inequivalent 
$V$-modules. Thus there are only finitely many 
inequivalent irreducible $V$-modules.
\epfv

Given an $\mathbb{N}$-gradable weak $V$-module $W$, we can also construct
an  $\tilde{A}(V)$-bimodule: For $u\in V$ and $w\in W$, we define
\begin{eqnarray*}
u\bullet w&=&\res_{y}y^{-1}Y\left(u, \frac{1}{2\pi i}\log(1+y)\right)w\nn
&=&\res_{x}\frac{2\pi ie^{2\pi i x}}{e^{2\pi ix}-1}Y(u, x)w,\nn
w\bullet u&=&\res_{y}y^{-1}e^{\frac{1}{2\pi i}\log(1+y)L(-1)}
Y\left(u, -\frac{1}{2\pi i}\log(1+y)\right)w\nn
&=&\res_{x}\frac{2\pi ie^{2\pi i x}}{e^{2\pi ix}-1}
e^{xL(-1)}Y(u, -x)w.
\end{eqnarray*}
Let 
$\tilde{O}(W)$ be the subspace of
$V$ spanned by elements of the form 
$$\res_{y}y^{-n}Y\left(u, \frac{1}{2\pi i}\log(1+y)\right)w
=\res_{x}\frac{2\pi ie^{2\pi i x}}{(e^{2\pi ix}-1)^{n}}
Y(u, x)w$$
for $n\in \mathbb{Z}_{+}+1$ and $u\in V$ and $w\in W$. 
Let $\tilde{A}(W)=W/\tilde{O}(W)$.

\begin{prop}
The left and right actions of $V$ on $W$ 
induce an $\tilde{A}(V)$-bimodule structure on 
$\tilde{A}(W)$.
\end{prop}
\pf
The proof is the same as the proof that $\tilde{A}(V)$ is 
an associative algebra above.
\epfv

We now assume that every $\mathbb{N}$-gradable weak $V$-modules
is completely reducible. It is easy to show that an irreducible 
$\mathbb{N}$-gradable weak $V$-module must be a $V$-module (see
\cite{Z}).
Thus we can reduce the study of $\mathbb{N}$-gradable weak $V$-modules
to the study of $V$-modules. 

Let $W_{1}$, $W_{2}$ and $W_{3}$ be  $V$-modules
and $\mathcal{Y}$ an intertwining operator of 
type $\binom{W_{3}}{ W_{1}W_{2}}$. Then 
$\tilde{A}(W_{1})\otimes_{\tilde{A}(V)} T(W_{2})$ and $T(W_{3})$
are both left $\tilde{A}(V)$-modules. 
For homogeneous $w_{1}\in W_{1}$, let 
$$o_{\mathcal{Y}}(w_{1})=\mathcal{Y}_{\swt w_{1}-1}(w_{1}),$$
where for $n\in \mathbb{C}$ and $w_{1}\in W_{1}$, 
$\mathcal{Y}_{n}(w_{1})\in \hom(W_{2}, W_{3})$ is given by 
$$\mathcal{Y}(w_{1}, x)
=\sum_{n\in \mathbb{C}}\mathcal{Y}_{n}(w_{1})x^{-n-1}.$$
We define $o_{\mathcal{Y}}(w_{1})$ for general $w_{1}\in W_{1}$ 
by linearity. 

\begin{lemma}\label{o-y}
For $w_{1}\in \tilde{O}(W_{1})$, $o_{\mathcal{Y}}(\mathcal{U}(1)w_{1})=0$.
\end{lemma}
\pf
The proof is the same as the one for the first conclusion in
Proposition \ref{functor-f}.
\epfv

Let
$$\rho(\mathcal{Y}): \tilde{A}(W_{1})\otimes_{\tilde{A}(V)} T(W_{2})
\to W_{3}$$
be defined by 
\begin{eqnarray*}
\rho(\mathcal{Y})((w_{1}+\tilde{O}(W_{1}))\otimes w_{2})
&=&o_{\mathcal{Y}}(\mathcal{U}(1)w_{1})w_{2}\nn
&=&\res_{x}x^{-1}\mathcal{Y}(\mathcal{U}(x)w_{1}, x)w_{2}
\end{eqnarray*}
for $w_{1}\in W_{1}, w_{2}\in T(W_{2})$. By Lemma \ref{o-y}, 
$\rho(\mathcal{Y})$ is indeed well defined. We have:

\begin{prop}\label{rho-isom}
The image of $\rho(\mathcal{Y})$ is in fact in $T(W_{3})$ and 
$\rho(\mathcal{Y})$ is in fact in 
$$\hom_{\tilde{A}(V)}(\tilde{A}(W_{1})\otimes_{\tilde{A}(V)} T(W_{2}),
T(W_{3})).$$
The map 
\begin{eqnarray*}
\rho: \mathcal{V}_{W_{1}W_{2}}^{W_{3}}&\to& 
\hom_{\tilde{A}(V)}(\tilde{A}(W_{1})\otimes_{\tilde{A}(V)} T(W_{2}),
T(W_{3}))\nn
\mathcal{Y}&\mapsto& \rho(\mathcal{Y})
\end{eqnarray*}
is a linear isomorphism.
\end{prop}
\pf
For any $w_{1}\in W_{1}$, since the 
weight of $o_{\mathcal{Y}}(\mathcal{U}(1)w_{1})$ is $0$, 
it is clear that the image of $\rho(\mathcal{Y})$ is in fact in 
$T(W_{3})$. Thus $\rho$ is a linear map from 
$\mathcal{V}_{W_{1}W_{2}}^{W_{3}}$ to 
$\hom_{\tilde{A}(V)}(\tilde{A}(W_{1})\otimes_{\tilde{A}(V)} T(W_{2}),
T(W_{3}))$. We still need to 
show that it is in fact an isomorphism. 

To show that $\rho$ is injective, we need the following obvious fact:
For any $\tilde{w}_{1}\in W_{1}$, 
$\tilde{w}_{2}\in W_{2}$
and $\tilde{w}'_{3}$, using the Jacobi identity and (\ref{chg-var0}), 
we can always write 
$\langle \tilde{w}'_{3}, \mathcal{Y}(\tilde{w}_{1}, x)\tilde{w}_{2}\rangle$
as a linear combination of series of the form 
$\langle w'_{3}, \mathcal{Y}(w_{1}, x)w_{2}\rangle$
for $w_{1}\in W_{1}$, $w_{2}\in T(W_{2})$ and $w'_{3}\in T(W'_{3})$
with Laurent polynomials of $x$ as coefficients. 
If $\rho(\mathcal{Y})=0$, 
$\langle w'_{3}, \mathcal{Y}(w_{1}, x)w_{2}\rangle=0$
for $w_{1}\in W_{1}$, $w_{2}\in T(W_{2})$ and $w'_{3}\in T(W'_{3})$.
Thus the  fact above shows that
$\mathcal{Y}=0$. So $\rho$ is injective.

We now prove that $\rho$ is surjective. Given any element $f$ of 
$$\hom_{\tilde{A}(V)}(\tilde{A}(W_{1})\otimes_{\tilde{A}(V)} T(W_{2}),
T(W_{3})),$$
we want to construct an element $\mathcal{Y}^{f}$
of $\mathcal{V}_{W_{1}W_{2}}^{W_{3}}$ such that $\rho(\mathcal{Y}^{f})=f$.
We assume that $W_{1}$, $W_{2}$, $W_{3}$ and thus $W'_{3}$ are all
irreducible $V$-modules. The general case follows from this case
by using the assumption that every $\mathbb{N}$-gradable 
weak $V$-module is completely reducible. Since $W_{1}$, $W_{2}$ and $W_{3}$
are irreducible $V$-modules, there exists $h_{1}, h_{2},
h_{3}\in \mathbb{C}$ such that the weights of nonzero homogeneous 
elements of $W_{1}$,
$W_{2}$ and $W_{3}$ are in $h_{1}+\mathbb{N}$, $h_{2}+\mathbb{N}$
and $h_{3}+\mathbb{N}$, respectively, and 
$(W_{1})_{(h_{1})}$, $(W_{3})_{(h_{3})}$ and $(W_{1})_{(h_{3})}$
are not $0$. Let $h=h_{3}-h_{1}-h_{2}$.
Then we know that for any intertwining operator $\mathcal{Y}$ of type 
$\binom{W_{3}}{W_{1}W_{2}}$, $\mathcal{Y}(w_{1}, x)w_{2}\in 
x^{h}W_{3}[[x, x^{-1}]]$. 

We consider the affinization 
$V[t, t^{-1}]$ and also 
$$t^{-h}W_{1}[t, t^{-1}]=W_{1}\otimes t^{-h}\mathbb{C}[t, t^{-1}].$$
For simplicity, we shall use $u(m)$ and $w_{1}(n)$ to denote 
$u\otimes t^{m}$ and $w_{1}\otimes t^{n}$, respectively,
for $u\in V$, $w_{1}\in W_{1}$, $m\in \mathbb{Z}$ and $n\in h+\mathbb{Z}$.
We consider the tensor algebra $\mathcal{T}(V[t, t^{-1}]\oplus 
t^{-h}W_{1}[t, t^{-1}])$
generated by $V[t, t^{-1}]$ and $t^{-h}W_{1}[t, t^{-1}]$. The 
tensor algebra $\mathcal{T}(V[t, t^{-1}])$
is a subalgebra of $\mathcal{T}(V[t, t^{-1}]\oplus 
t^{-h}W_{1}[t, t^{-1}])$ and $t^{-h}W_{1}[t, t^{-1}]$ is a subspace. Let
$\mathcal{T}_{V; W_{1}}$ 
be the $\mathcal{T}(V[t, t^{-1}])$-sub-bimodule of
$\mathcal{T}(V[t, t^{-1}]\oplus 
t^{-h}W_{1}[t, t^{-1}])$ generated by $t^{-h}W_{1}[t, t^{-1}]$.
For simplicity we shall omit the 
tensor product sign for elements.

Consider 
$\mathcal{T}_{V; W_{1}}\otimes T(W_{2})$. For simplicity 
again we shall omit the 
tensor product sign for elements. So $\mathcal{T}_{V; W_{1}}\otimes T(W_{2})$ 
is spanned by elements
of the form 
$$u_{1}(m_{1})\cdots u_{k}(m_{k})w_{1}(n)
u_{k+1}(m_{k+1})\cdots u_{k+l}(m_{k+l})w_{2}$$
for
$u_{i}\in V$, $m_{i}\in \mathbb{Z}$, $i=1, \dots, k+l$, 
$w_{1}\in W_{1}$ and $w_{2}\in T(W_{2})$.
For any $u\in V$, $m\in \mathbb{Z}$, $u(m)$ acts from the left on 
$\mathcal{T}_{V; W_{1}}\otimes T(W_{2})$.
For homogeneous $u_{i}\in V$, $m_{i}\in \mathbb{Z}$, $i=1, \dots, k+l$, 
and homogeneous $w_{1}\in W_{1}$ and $w_{2}\in T(W_{2})$,
we define the weight of 
$$u_{1}(m_{1})\cdots u_{k}(m_{k})w_{1}(n)
u_{k+1}(m_{k+1})\cdots u_{k+l}(m_{k+l})w_{2}$$
to be 
$$(\wt u_{1}-m_{1}-1)+\cdots +(\wt u_{k}-m_{k+l}-1)
+(\wt w_{1}-n-1)+\wt w_{2}.$$
For any $u\in V$ and $w_{1}\in W_{1}$, let
\begin{eqnarray*}
Y_{t}(u, x): \mathcal{T}_{V; W_{1}}\otimes T(W_{2})
&\to& (\mathcal{T}_{V; W_{1}}\otimes T(W_{2}))[[x, x^{-1}]],\\
\mathcal{Y}_{t}(w_{1}, x): \mathcal{T}_{V; W_{1}}\otimes T(W_{2})
&\to& x^{h}(\mathcal{T}_{V; W_{1}}\otimes T(W_{2}))[[x, x^{-1}]]
\end{eqnarray*}
be defined by 
\begin{eqnarray*}
Y_{t}(u, x)&=&\sum_{m\in \mathbb{Z}}u(m)x^{-m-1},\\
Y_{t}(w_{1}, x)&=&\sum_{n\in h+\mathbb{Z}}w_{1}(n)x^{-n-1},
\end{eqnarray*}
respectively.
For homogeneous elements $u\in V$ and $w_{1}\in W_{1}$, 
let $o_{t}(u)=u(\wt u-1)$ and $o_{t}(w_{1})=w_{1}(\wt w_{1}-1)$. 
Using linearity, we extend $o_{t}(u)$ and $o_{t}(w_{1})$ 
to nonhomogeneous $u$ and $w_{1}$.

Let $\mathcal{I}_{V; W_{1}, W_{2}}$ be the $h_{1}+h_{2}+\mathbb{Z}$-graded 
$\mathcal{T}(V[t,t^{-1}])$-submodule 
of $\mathcal{T}_{V; W_{1}}\otimes T(W_{2})$
generated by elements of the forms
$u(m)w_{2}$ ($u\in V$, $m\in \mathbb{Z}$, 
$\wt u-m-1< 0$, $w_{2}\in T(W_{2})$), 
$w_{1}(n)w_{2}$ ($w_{1}\in W_{1}$, $n\in \mathbb{Z}-h$, 
$\wt w_{1}-n-1+\wt w_{2}< h_{3}$, $w_{2}\in T(W_{2})$), 
$o_{t}(\mathcal{U}(1)u)w_{2}-\rho(u+\tilde{O}(V))w_{2}$ 
($u\in V$, $w_{2}\in T(W_{2})$),  and
the coefficients in $x_{1}$ and $x_{2}$ of 
\begin{eqnarray*}
\lefteqn{Y_{t}(u, x_{1})Y_{t}(v, x_{2})w
-Y_{t}(v, x_{2})Y_{t}(u, x_{1})w}&\nn
&&-\res_{x_{0}}x_{2}^{-1}\delta\left(\frac{x_{1}-x_{0}}{x_{2}}\right)
Y_{t}(Y(u, x_{0})v, x_{2})w,\\
\lefteqn{Y_{t}(u, x_{1})\mathcal{Y}_{t}(w_{1}, x_{2})w_{2}
-\mathcal{Y}_{t}(w_{1}, x_{2})Y_{t}(u, x_{1})w}&\nn
&&-\res_{x_{0}}x_{2}^{-1}\delta\left(\frac{x_{1}-x_{0}}{x_{2}}\right)
\mathcal{Y}_{t}(Y(u, x_{0})w_{1}, x_{2})w_{2}
\end{eqnarray*}
($u, v\in V$, $w_{1}\in W_{1}$ and 
$w_{2}\in \mathcal{T}_{V; W_{1}}\otimes T(W_{2})$). 
Let 
$$S_{1}(V; W_{1}, W_{2})=(\mathcal{T}_{V; W_{1}}\otimes T(W_{2}))/
\mathcal{I}_{V; W_{1}, W_{2}}.$$
Then 
$S_{1}(V; W_{1}, W_{2})$ is also  an $h_{1}+h_{2}+\mathbb{Z}$-graded 
$\mathcal{T}(V[t,t^{-1}])$-module. 
In fact, by definition of $\mathcal{I}_{V; W_{1}, W_{2}}$, 
$S_{1}(V; W_{1}, W_{2})$ is spanned by elements
of the form 
$$u_{1}(m_{1})\cdots u_{k}(m_{k})w_{1}(n)
w_{2}+\mathcal{I}_{V; W_{1}, W_{2}}$$
for homogeneous
$u_{i}\in V$, $m_{i}<\wt u_{i}-1$, $i=1, \dots, k$, homogeneous 
$w_{1}\in W_{1}$, $n\le \wt w_{1}-1+\wt w_{2}-h_{3}$
and $w_{2}\in T(W_{2})$. In particular, we 
see that $S_{1}(V; W_{1}, W_{2})$ has an $\mathbb{N}$-grading.
Thus for $u\in V$, $w_{1}\in W_{1}$ and $w_{2}\in S_{1}(V; W_{1}, W_{2})$,
$u(m)w_{2}=0$ and $w_{1}(n)w_{2}=0$ when $m$ and $n$ are
sufficiently large.

Let $\mathcal{J}_{V; W_{1}, W_{2}}$ be the $h_{1}+h_{2}+\mathbb{Z}$-graded 
$\mathcal{T}(V[t,t^{-1}])$-submodule of $S_{1}(V; W_{1}, W_{2})$
generated by the coefficients in $x$ of  
\begin{eqnarray*}
&{\displaystyle Y_{t}(L(-1)u, x)w_{2}-\frac{d}{dx}Y_{t}(u, x)w_{2},}&\\
&{\displaystyle \mathcal{Y}_{t}(L(-1)w_{1}, x)w_{2}-\frac{d}{dx}
\mathcal{Y}_{t}(w_{1}, x)w_{2}}&
\end{eqnarray*}
and the coefficients in $x_{0}$ and $x_{2}$ of 
\begin{eqnarray*}
&{\displaystyle Y_{t}(Y(u, x_{0})v, x_{2})w_{2}-\res_{x_{1}}
x_{0}^{-1}\delta\left(\frac{x_{1}-x_{2}}{x_{0}}\right)
Y_{t}(u, x_{1})Y_{t}(v, x_{2})w_{2}}&\nn
&{\displaystyle +x_{0}^{-1}\delta\left(\frac{x_{2}-x_{1}}{-x_{0}}\right)
Y_{t}(v, x_{2})Y_{t}(u, x_{1})w_{2}},&\\
&{\displaystyle \mathcal{Y}_{t}(Y(u, x_{0})w_{1}, x_{2})w_{2}-\res_{x_{1}}
x_{0}^{-1}\delta\left(\frac{x_{1}-x_{2}}{x_{0}}\right)
Y_{t}(u, x_{1})\mathcal{Y}_{t}(w_{1}, x_{2})w_{2}}&\nn
&{\displaystyle +x_{0}^{-1}\delta\left(\frac{x_{2}-x_{1}}{-x_{0}}\right)
\mathcal{Y}_{t}(w_{1}, x_{2})Y_{t}(u, x_{1})w_{2}},&\\
\end{eqnarray*}
for $u, v\in V$, $w_{1}\in W_{1}$ and $w_{2}\in S_{1}(V; W_{1}, W_{2})$.

Let $S(V; W_{1}, W_{2})
=S_{1}(V; W_{1}, W_{2})/\mathcal{J}_{V; W_{1}, W_{2}}$. Then 
$S(V; W_{1}, W_{2})$ is also 
a $\mathcal{T}(V[t,t^{-1}])$-module with an $\mathbb{N}$-grading. We can
still use elements of $\mathcal{T}_{V; W_{1}}\otimes T(W_{2})$
to represent elements of $S(M)$. But note that these 
elements now satisfy relations. Now we have operators 
$Y_{t}(u, x)$ and $\mathcal{Y}_{t}(w_{1}, x)$ for $u\in V$ and 
$w_{1}\in W_{1}$ acting on 
$S(V; W_{1}, W_{2})$. By construction, these operators satisfy the lower 
truncation property, the identity property (for $Y_{t}$),
the commutator formula (for $Y_{t}$ and 
for $Y_{t}$ and $\mathcal{Y}_{t}$),
the associator formula (for $Y_{t}$ and $\mathcal{Y}_{t}$)
and the $L(-1)$-derivative property (for $Y_{t}$ and $\mathcal{Y}_{t}$).

We also have a linear map $\mu: S(V; W_{1}, W_{2})\to W_{3}$ 
defined as follows:
For $u_{i}\in V$, $m_{i}<\wt u_{i}-1$, $i=1, \dots, k$, 
$w_{1}\in W_{1}$, $n\le \wt w_{1}-1$
and $w_{2}\in T(W_{2})$, we define 
$$\mu(u_{1}(m_{1})\cdots u_{k}(m_{k})w_{1}(n)
w_{2})=(u_{1})_{m_{1}}\cdots (u_{k})_{m_{k}}
f((w_{1}+\tilde{O}(W_{1}))\otimes w_{2}).$$
Since all the relations in $S(V; W_{1}, W_{2})$ for 
representatives of the form 
$$u_{1}(m_{1})\cdots u_{k}(m_{k})w_{1}(n)
w_{2}$$
above are also satisfied by their 
images in $W_{3}$, $\mu$ is well defined.

Now we construct an intertwining operator $\mathcal{Y}^{f}$ 
of type $\binom{W_{3}}{W_{1}W_{2}}$ as follows: 
By Proposition \ref{reductive},  $W_{2}$ is isomorphic to $S(T(W_{2}))$. 
So we can work with 
$S(T(W_{2}))$ instead of $W_{2}$.  We know that $S(T(W_{2}))$ is spanned by 
elements of the form $u_{1}(m_{1})\cdots u_{k}(m_{k}) w_{2}$
for $u_{i}\in V$, $m_{i}<\wt u_{i}-1$ for $i=1, \dots, k$ and $w_{2}\in 
T(W_{2})$. 
Let $w_{1}\in W_{1}$. 
We define 
$$\mathcal{Y}^{f}(w_{1}, x)u_{1}(m_{1})\cdots u_{k}(m_{k}) w_{2}
=\mu(\mathcal{Y}_{t}(w_{1}, x)u_{1}(m_{1})\cdots u_{k}(m_{k})w_{2}).$$
Since $\mathcal{Y}_{t}$ satisfies the commutator formula, the 
associator formula and the $L(-1)$-derivative property,
so does $\mathcal{Y}^{f}$. Thus $\mathcal{Y}^{f}$ satisfies the Jacobi
identity and the $L(-1)$-derivative property. So it is an intertwining
operator of the desired type. It is clear from the construction that
$\rho(\mathcal{Y}^{f})=f$.
\epfv

\begin{lemma}\label{tr}
Let $A$ be a semisimple associative algebra over $\mathbb{C}$, $\omega$
an element in the center of $A$,
$M$ an $A$-bimodule and
$F: M\to \mathbb{C}$ a linear functional on $M$ satisfying the 
following property:
\begin{enumerate}

\item For $u\in A$ and $w\in M$, $F(uw)=F(wu)$.

\item There exist $h\in \mathbb{C}$ and $s\in \mathbb{Z}_{+}$
such that for $w\in M$, $F((\omega-h)^{s}w)=0$.

\end{enumerate}

Then there exist irreducible left $A$-modules $M_{i}$ and left module maps
\begin{eqnarray*}
f_{i}: M\otimes_{A}M_{i} &\to& M_{i}\nn
w\otimes_{A} w_{i}&\mapsto &f_{i}(w)w_{i}
\end{eqnarray*}
for $i=1, \dots, m$
such that 
$$F(w)=\sum_{i=1}^{m}\tr_{M_{i}}f_{i}(w)$$
for $w\in M$.
\end{lemma}
\pf
This is a result in the classical theory of semisimple
associative algebras. We just give the idea of the proof. It is easy 
for the reader to fill in all the details using the theory presented 
in, for example, \cite{J} or \cite{FD}. One can also find a proof in 
\cite{M1}.

Since $A$ is semisimple, it must be a direct sum of simple ideals.
Because of this fact, we can reduce our lemma to the case that $A$ is 
simple. But when $A$ is simple, we know that $A$ is isomorphic to a 
matrix algebra. The lemma now can be verified directly for simple matrix
algebras.
\epfv

\renewcommand{\theequation}{\thesection.\arabic{equation}}
\renewcommand{\thethm}{\thesection.\arabic{thm}}
\setcounter{equation}{0}
\setcounter{thm}{0}

\section{Modular invariance}

In this section, we prove the modular invariance of the space 
of chiral genus-one correlation functions using the results we 
have obtained in the preceding sections. 

We first discuss the modular invariance of the system 
(\ref{eqn1})--(\ref{eqn2}). We need the following
modular transformation formulas (see, for example, 
\cite{Ko}): For any 
$$\left(\begin{array}{cc}
\alpha&\beta\\
\gamma&\delta
\end{array}\right)\in SL(2, \mathbb{Z}),$$
we have 
\begin{eqnarray}
G_{2}\left(\frac{\alpha\tau+\beta}{\gamma\tau +\delta}\right)
&=&(\gamma\tau+\delta)^{2}G_{2}(\tau)
-2\pi i\gamma(\gamma\tau+\delta),\label{mod-g2}\\
G_{2k}\left(\frac{\alpha\tau+\beta}{\gamma\tau+\delta}\right)&=&
(\gamma\tau+\delta)^{2k}G_{2k}(\tau),\label{mod-g2k}\\
\wp_{m}\left(\frac{z}{\gamma\tau+\delta}; 
\frac{\alpha\tau+\beta}{\gamma\tau+\delta}\right)&=&
(\gamma\tau+\delta)^{m}\wp_{m}(z, \tau),\label{mod-wpm}
\end{eqnarray}
for $k\ge 2$ and $m\ge 1$. 

We have:

\begin{prop}\label{eqn-mod-inv}
Let $\varphi(z_{1}, \dots, z_{n}; \tau)$ be a solution of 
the system (\ref{eqn1})--(\ref{eqn2}) with $q=q_{\tau}$. 
Then for any 
$$\left(\begin{array}{cc}
\alpha&\beta\\
\gamma&\delta
\end{array}\right)\in SL(2, \mathbb{Z}),$$
$$\left(\frac{1}{\gamma\tau +\delta}\right)^{\swt w_{1}+\cdots \swt w_{n}}
\varphi\left(\frac{z_{1}}{c\tau +d}, \dots, \frac{z_{n}}{c\tau +d}; 
\frac{\alpha\tau+\beta}{\gamma\tau +\delta}\right)$$
is also a solution of the system (\ref{eqn1})--(\ref{eqn2})
with $q=q_{\tau}$.
\end{prop}
\pf
The proof is a straightforward calculation using 
(\ref{mod-g2})--(\ref{mod-wpm}).
\epfv

From this result, we know that the space of solutions of the system 
(\ref{eqn1})--(\ref{eqn2})  with $q=q_{\tau}$ 
is invariant under the action of the 
modular group $SL(2, \mathbb{Z})$. But this is not enough for 
the modular invariance we would like to prove because we want to 
prove that the space of those solutions obtained from the $q_{\tau}$-traces
of products of geometrically modified intertwining operators 
are invariant under this action of $SL(2, \mathbb{Z})$.

We need the following:

\begin{thm}\label{summary}
Let $V$ be a vertex operator algebra satisfying the following conditions:

\begin{enumerate}

\item For $n<0$, $V_{(n)}=0$ and $V_{(0)}=\mathbb{C}\mathbf{1}$.

\item Every $\mathbb{N}$-gradable weak $V$-module is completely reducible.

\item $V$ is $C_{2}$-cofinite.

\end{enumerate}

Then all the conclusions of the results in Sections 1--6 hold. 
\end{thm}
\pf
By Corollary \ref{finitely-many}, there are only finitely many
inequivalent irreducible $V$-modules.  By a result of Anderson-Moore
\cite{AM} and Dong-Li-Mason \cite{DLM}, every irreducible $V$-module is
in fact $\mathbb{Q}$-graded.  Also in this case, it is clear that every
finitely-generated lower-truncated generalized $V$-module is a
$V$-module.  In \cite{H11}, the author proved that for such a vertex
operator algebra, the direct sum of a complete set of 
inequivalent irreducible
$V$-modules has a natural structure of intertwining operator algebra.
By a result of Abe, Buhl and Dong \cite{ABD}, we also know that for such
a vertex operator algebra $V$, every $V$-module is $C_{2}$-cofinite.
Thus the conditions for $V$ needed in Sections 1--6 are all satisfied.
\epfv 

Let $V$ be a vertex operator algebra satisfying the conditions in
Theorem \ref{summary}. 
Let $W_{i}$ be $V$-modules and $w_{i}\in W_{i}$ for $i=1, \dots, n$.
For  any $V$-modules
$\tilde{W}_{i}$ and any intertwining operators
$\mathcal{Y}_{i}$, $i=1, \dots, n$,  of 
types $\binom{\tilde{W}_{i-1}}{W_{i}\tilde{W}_{i}}$, respectively,
we have a genus-one correlation function
$$\overline{F}_{\mathcal{Y}_{1}, \dots, \mathcal{Y}_{n}}(w_{1}, \dots, w_{n};
z_{1}, \dots, z_{n}; \tau).$$
Note that these multivalued functions actually have preferred
branches in the region $1>|q_{z_{1}}|>\cdots >|q_{z_{n}}|>|q_{\tau}|>0$
given by the intertwining operators $\mathcal{Y}_{1}, 
\dots, \mathcal{Y}_{n}$. Thus linear 
combinations of these functions make sense. 
For fixed $V$-modules $W_{i}$ and $w_{i}\in W_{i}$
for $i=1, \dots, n$, we now denote the vector space
spanned by all such functions by  $\mathcal{F}_{w_{1}, \dots, w_{n}}$. 
The following theorem is the main result of this section:

\begin{thm}\label{mod-inv}
Let $V$ be a vertex operator algebra satisfying the
conditions in Theorem \ref{summary}.
Then for any $V$-modules $\tilde{W}_{i}$
and any
intertwining operators $\mathcal{Y}_{i}$ ($i=1, \dots, n$)  of 
types $\binom{\tilde{W}_{i-1}}{W_{i}\tilde{W}_{i}}$, respectively,
and any 
$$\left(\begin{array}{cc}
\alpha&\beta\\
\gamma&\delta
\end{array}\right)\in SL(2, \mathbb{Z}),$$
\begin{eqnarray*}
\lefteqn{\overline{F}_{\mathcal{Y}_{1}, \dots, \mathcal{Y}_{n}}
\Biggl(\left(\frac{1}{\gamma\tau+\delta}\right)^{L(0)}w_{1}, \dots,
\left(\frac{1}{\gamma\tau+\delta}\right)^{L(0)}w_{n};}\nn
&&\quad\quad\quad\quad\quad\quad\quad\quad\quad\quad\quad\quad
\quad\quad\quad
\frac{z_{1}}{\gamma\tau+\delta}, \dots, \frac{z_{n}}{\gamma\tau+\delta}; 
\frac{\alpha\tau+\beta}{\gamma\tau+\delta}\Biggr)
\end{eqnarray*}
is in $\mathcal{F}_{w_{1}, \dots, w_{n}}$. 
\end{thm}
\pf
By Theorem \ref{summary}, all the results in Sections 1--6 can be used.

The case of $n=1$, as we have mentioned in the introduction,
was proved by Miyamoto in \cite{M1} using the method 
of Zhu in \cite{Z}. Since the differential equations we obtained
in this paper are explicitly modular invariant, 
we have a simpler proof than the one given in \cite{M1}. 
Here we give this proof.
In this case, the identity 
(\ref{mod-inv-der}) becomes
\begin{eqnarray}\label{mod-inv-1}
\lefteqn{\left(2\pi i\frac{\partial}{\partial \tau}
+G_{2}(\tau)\wt w_{1}
+G_{2}(\tau)z_{1}\frac{\partial}{\partial z_{1}}\right)
F_{\mathcal{Y}_{1}}(w_{1};
z_{1}; q_{\tau})}\nn
&&=F_{\mathcal{Y}_{1}}(L(-2)w_{1};
z_{1}; q_{\tau}) -\sum_{k\in \mathbb{Z}_{+}}G_{2k+2}(\tau)
F_{\mathcal{Y}_{1}}(L(2k)w_{1};
z_{1}; q_{\tau}).
\end{eqnarray}
But in this case, 
\begin{eqnarray*}
\lefteqn{\frac{\partial}{\partial z_{1}}
F_{\mathcal{Y}_{1}}(w_{1}; z_{1}; q_{\tau})}\nn
&&=\tr_{\tilde{W}_{1}}\mathcal{Y}_{1}((2\pi iL(0)+
2\pi iq_{z_{1}}L(-1))
\mathcal{U}(q_{z_{1}})w_{1}, 
q_{z_{1}})q^{L(0)-\frac{c}{24}}\nn
&&=2\pi i\tr_{\tilde{W}_{1}}[L(0), \mathcal{Y}_{1}(
\mathcal{U}(q_{z_{1}})w_{1}, 
q_{z_{1}})]q^{L(0)-\frac{c}{24}}\nn
&&=0.
\end{eqnarray*}
In other words, $F_{\mathcal{Y}_{1}}(w_{1};
z_{1}; q_{\tau})$ is in fact independent of $z_{1}$. 
So from (\ref{mod-inv-1}), we obtain
\begin{eqnarray}\label{mod-inv-2}
\lefteqn{\left(2\pi i\frac{\partial}{\partial \tau}
+G_{2}(\tau)\wt w_{1}\right)
F_{\mathcal{Y}_{1}}(w_{1};
z_{1}; q_{\tau})}\nn
&&=F_{\mathcal{Y}_{1}}\left(L(-2)w_{1}-\sum_{k\in \mathbb{Z}_{+}}G_{2k+2}(\tau)
L(2k)w_{1};
z_{1}; q_{\tau}\right).
\end{eqnarray}

We use $\tau'$ to denote $\frac{\alpha\tau+\beta}{\gamma\tau+\delta}$ and 
$z'_{1}=\frac{z_{1}}{\gamma\tau+\delta}$. 
Then (\ref{mod-inv-2}) also holds with $\tau$ replaced by $\tau'$
and $z_{1}$ replaced by $z'_{1}$.
Using the modular transformation property (\ref{mod-g2})
and (\ref{mod-g2k}) for $G_{2k}(\tau)$ for $k\in \mathbb{Z}_{+}$
and the fact that $F_{\mathcal{Y}_{1}}(w_{1};
z_{1}; q_{\tau})$ is independent of $z_{1}$, we see that by 
a straightforward calculation, 
(\ref{mod-inv-2}) with $\tau$, $z_{1}$ and $w_{1}$ replaced by $\tau'$,  
$z'_{1}$ and $(\gamma\tau+\delta)^{-L(0)}w_{1}$, respectively,
is equivalent to 
\begin{eqnarray*}
\lefteqn{\left(2\pi i\frac{\partial}{\partial \tau}
+G_{2}(\tau)\wt w_{1}\right)
F_{\mathcal{Y}_{1}}((\gamma\tau+\delta)^{-L(0)}w_{1};
z'_{1}; q_{\tau'})}\nn
&&=F_{\mathcal{Y}_{1}}((\gamma\tau+\delta)^{-L(0)}L(-2)w_{1};
z'_{1}; q_{\tau'})\nn
&&\quad -\sum_{k\in \mathbb{Z}_{+}}G_{2k+2}(\tau)
F_{\mathcal{Y}_{1}}((\gamma\tau+\delta)^{-L(0)}L(2k)w_{1};
z'_{1}; q_{\tau'})
\end{eqnarray*}
or equivalently
\begin{eqnarray}\label{mod-inv-3}
\lefteqn{2\pi i\frac{\partial}{\partial \tau}
F_{\mathcal{Y}_{1}}((\gamma\tau+\delta)^{-L(0)}w_{1};
z'_{1}; q_{\tau'})}\nn
&&=F_{\mathcal{Y}_{1}}((\gamma\tau+\delta)^{-L(0)}L(-2)w_{1};
z'_{1}; q_{\tau'})\nn
&&\quad  -\sum_{k\in \mathbb{N}}G_{2k+2}(\tau)
F_{\mathcal{Y}_{1}}((\gamma\tau+\delta)^{-L(0)}L(2k)w_{1};
z'_{1}; q_{\tau'})
\end{eqnarray}
The $n=1$ cases of the identities (\ref{identity0.5}) and (\ref{identity2})
give
\begin{equation}\label{mod-inv-4}
F_{\mathcal{Y}_{1}}((\gamma\tau+\delta)^{-L(0)}u_{0}w_{1};
z'_{1}; q_{\tau'})=0
\end{equation}
and 
\begin{eqnarray}\label{mod-inv-5}
\lefteqn{F_{\mathcal{Y}_{1}}
\Biggl((\gamma\tau+\delta)^{-L(0)}u_{-2}w_{1}}\nn
&&\quad\quad\quad\quad\quad\quad
 +\sum_{k\in \mathbb{Z}_{+}}(2k+1)G_{2k+2}(\tau)(\gamma\tau+\delta)^{-L(0)}
u_{2k}w_{1};
z'_{1}; q_{\tau'}\Biggr)\nn
&&=0,
\end{eqnarray}
respectively, where in (\ref{mod-inv-5}), we have used (\ref{mod-g2k}).

Using (\ref{wp0}), we have 
\begin{eqnarray}\label{mod-inv-6}
\lefteqn{L(-2)w_{1}-\sum_{k\in \mathbb{N}}G_{2k+2}(\tau)
L(2k)w_{1}}\nn
&&=\res_{x}(\wp_{1}(x; \tau)-G_{2}(\tau)x)Y(\omega, x)w_{1}
\end{eqnarray}
and 
\begin{eqnarray}\label{mod-inv-7}
&{\displaystyle u_{-2}w_{1}+\sum_{k\in \mathbb{Z}_{+}}(2k+1)G_{2k+2}(\tau)
u_{2k}w_{1}}&\nn
&{\displaystyle =\res_{x}\wp_{2}(x; \tau)Y(u, x)w_{1}.}
\end{eqnarray}
Using (\ref{wp.5}) and (\ref{chg-var0}), we see that the constant terms 
of the expansions of (\ref{mod-inv-6}) and (\ref{mod-inv-7})
as power series in $q_{\tau}$ 
are 
\begin{eqnarray}\label{mod-inv-8}
\lefteqn{\res_{x}\pi i\frac{e^{2\pi ix}+1}{e^{2\pi ix}-1}
Y(omega, x)w_{1}}\nn
&&=\res_{x}\pi i\frac{2\pi ie^{2\pi ix}}{e^{2\pi ix}-1}
Y(\omega, x)w_{1}
-\pi i \res_{x} Y(u, x)w_{1}\nn
&&=\omega\bullet w_{1}-\pi i \omega_{0}w_{1}
\end{eqnarray}
and 
\begin{eqnarray}\label{mod-inv-9}
\lefteqn{\res_{x}\wp_{2}(x; \tau)Y(u, x)w_{1}}\nn
&&=\res_{x}\left(-\frac{\pi^{2}}{3}
+\frac{(2\pi i)^{2}e^{2\pi ix}}{(e^{2\pi ix}-1)^{2}}\right)
Y(u, x)w_{1}\nn
&&=-\frac{\pi^{2}}{3}u_{0}w_{1}+2\pi i
\res_{x}\frac{2\pi ie^{2\pi ix}}{(e^{2\pi ix}-1)^{2}}
Y(u, x)w_{1},\nn
&&
\end{eqnarray}
respectively.

By Proposition \ref{eqn-mod-inv},  
we know that $F_{\mathcal{Y}_{1}}((\gamma\tau+\delta)^{-L(0)}w_{1};
z'_{1}; q_{\tau'})$ satisfies the same equation of regular singular points
as the one for $F_{\mathcal{Y}_{1}}(w_{1};
z_{1}; q_{\tau})$. Thus we have
$$F_{\mathcal{Y}_{1}}((\gamma\tau+\delta)^{-L(0)}w_{1};
z'_{1}; q_{\tau'})=\sum_{k=0}^{K}\sum_{l=1}^{N}\sum_{m\in \mathbb{N}}
C_{k, l, m}(w_{1})\tau^{k}q_{\tau}^{r_{l}+m}$$
where $r_{l}$ for $l=1, \dots, N$ are real numbers such that 
$r_{l_{1}}-r_{l_{2}}\not \in \mathbb{Z}$ when $l_{1}\ne l_{2}$. From 
(\ref{mod-inv-3})--(\ref{mod-inv-9}), we obtain
\begin{eqnarray}
C_{K, l, 0}(u\bullet w_{1})&=&C_{K, l, 0}(w_{1}\bullet u),\label{mod-inv-10}\\
C_{K, l, 0}(\tilde{O}(W_{1}))&=&0,\label{mod-inv-11}\\
C_{K, l, 0}\left(\left(\omega-\frac{c}{24}-r_{l}\right)\bullet w_{1}\right)
&=&0.\label{mod-inv-12}
\end{eqnarray}
Thus we see that $C_{K, l, 0}$ gives a linear functional on 
the $\tilde{A}(V)$-module $\tilde{A}(W_{1})$ satisfying the conditions in
Lemma \ref{tr}. 
By Lemma \ref{tr}, we can find  irreducible left
$\tilde{A}(V)$-modules $M_{i}$
and left  $\tilde{A}(V)$-module maps 
$$f_{i}: 
\tilde{A}(W_{1})\otimes_{\tilde{A}(V)}M_{i}\to M_{i}$$ 
such that
for $i=1, \dots, m$ such that 
$$C_{K, l, 0}(w_{1})=\sum_{i=1}^{m}\tr_{M_{i}}f_{i}(w_{1})$$
for $w_{1}\in W_{1}$. By Propositions \ref{functor-f}, \ref{reductive}
and \ref{rho-isom}, 
there exist irreducible $V$-modules 
$W^{(1)}_{i}$ and intertwining operators 
$\mathcal{Y}^{(1)}_{i}$
of types $\binom{W^{(1)}_{i}}{W_{1}W^{(1)}_{i}}$ 
for $i=1, \dots, p$ such that 
$M_{i}=T(W^{(1)}_{i})$ and $\mathcal{Y}^{(1)}_{i}$ correspond
to $f_{i}$ for $i=1, \dots, p$. 

It is clear that 
$$\sum_{l=1}^{N}\sum_{m\in \mathbb{N}}
C_{K, l, m}(w_{1})q_{\tau}^{r_{l}+m}
-\sum_{i=1}^{p}F_{\mathcal{Y}^{(1)}_{i}}(w_{1};
z_{1}; q_{\tau})$$
must be of the form 
$$\sum_{l=1}^{N^{(1)}}\sum_{m\in \mathbb{N}}
C^{(1)}_{K, l, m}(w_{1})q_{\tau}^{r^{(1)}_{l}+m}$$
where for $l=1, \dots, N^{(1)}$, there must be $k$ satisfying $1\le k\le N$
such that  $r^{(1)}_{l}-r_{k}$ is a positive integer. In addition, this series 
also satisfies (\ref{mod-inv-3})--(\ref{mod-inv-5}) and thus 
$C^{(1)}_{K, l, 0}(\cdot)$ also satisfy (\ref{mod-inv-10})--(\ref{mod-inv-12}).
Repeating the argument above again and again, say $s$ times
and noticing that there are only finitely many inequivalent irreducible 
$V$-modules, 
finally we can 
find $V$-modules $W^{(s)}_{j}$ and intertwining operators 
$\mathcal{Y}^{(s)}_{j}$
of type $\binom{W^{(s)}_{j}}{W_{1}W^{(s)}_{j}}$ for $j=1, \dots, p^{(s)}$ 
such that 
$$\sum_{l=1}^{N}\sum_{m\in \mathbb{N}}
C_{K, l, m}(w_{1})q_{\tau}^{r_{l}+m}
-\sum_{j=1}^{p^{(s)}}F_{\mathcal{Y}^{(s)}_{j}}(w_{1};
z_{1}; q_{\tau})$$
is of the form 
$$\sum_{1=1}^{N^{(s)}}\sum_{m\in \mathbb{N}}
C^{(s)}_{K, l, m}(w_{1})q_{\tau}^{r^{(s)}_{l}+m},$$
where for $l=1, \dots, N^{(s)}$, $r^{(s)}_{l}$ are larger than the
lowest weights of all irreducible $V$-modules. 
But this happens only when this series is $0$. So we obtain
$$\sum_{l=1}^{N}\sum_{m\in \mathbb{N}}
C_{K, l, m}(w_{1})q_{\tau}^{r_{l}+m}
=\sum_{j=1}^{p^{(s)}}F_{\mathcal{Y}^{(s)}_{j}}(w_{1};
z_{1}; q_{\tau}).$$

We still need to show that $K=0$. Assume $K\ne 0$. Then $K-1\ge 0$. Let 
$$S_{k}(w_{i}; \tau)=\sum_{l=1}^{N}\sum_{m\in \mathbb{N}}
C_{k, l, m}(w_{1})q_{\tau}^{r_{l}+m}$$
for $k=1, \dots, K$. Then 
$$F_{\mathcal{Y}_{1}}((\gamma\tau+\delta)^{-L(0)}w_{1};
z'_{1}; q_{\tau'})=\sum_{k=1}^{K}S_{k}(w_{i}; \tau)\tau^{k}.$$
Using (\ref{mod-inv-3})--(\ref{mod-inv-5}), 
we obtain
\begin{eqnarray}
&S_{K-1}(u_{0}w_{1}; \tau)=0,&\label{mod-inv-13}\\
&{\displaystyle S_{K-1}\left(u_{-2}w_{1}+\sum_{k\in \mathbb{Z}_{+}}
(2k+1)G_{2k+2}(\tau)
u_{2k}w_{1}; \tau\right)=0,}& \label{mod-inv-14}\\
\lefteqn{(2\pi i)^{2}K S_{K}(w_{i}; \tau)+2\pi 
i\frac{d}{d\tau}S_{K-1}(w_{1}; \tau)}\nn
&{\displaystyle  =S_{K-1}\left(\left(L(-2)-\sum_{k\in \mathbb{N}}G_{2k+2}(\tau)
L(2k)\right)w_{1}; \tau\right),}&\label{mod-inv-15}\\
&{\displaystyle 2\pi i\frac{d}{d\tau}S_{K}(w_{1}; \tau)
=S_{K}\left(\left(L(-2)-\sum_{k\in \mathbb{N}}G_{2k+2}(\tau)
L(2k)\right)w_{1}; \tau\right).}&\label{mod-inv-16}
\end{eqnarray}
From  (\ref{mod-inv-15}) and (\ref{mod-inv-16}), we obtain
\begin{eqnarray}\label{mod-inv-17}
\lefteqn{4\pi i S_{K-1}(w_{i}; \tau)+(2\pi i)^{2}\frac{d^{2}}{d\tau^{2}}
S_{K-1}(w_{1}; \tau)}\nn
&&=S_{K-1}\left(\left(L(-2)-\sum_{k\in \mathbb{N}}G_{2k+2}(\tau)
L(2k)\right)^{2}w_{1}; \tau\right).
\end{eqnarray}
From (\ref{mod-inv-13}), (\ref{mod-inv-14}) and (\ref{mod-inv-17}),
we obtain
\begin{eqnarray*}
C_{K-1, l, 0}(u\bullet w_{1})&=&C_{K, l, 0}(w_{1}\bullet u),\\
C_{K-1, l, 0}(\tilde{O}(W_{1}))&=&0,\\
C_{K-1, l, 0}\left(\left(\omega-\frac{c}{24}-r_{l}\right)^{2}
\bullet w_{1}\right)
&=&0.
\end{eqnarray*}
The same argument as we have used above shows that 
$$S_{K-1}(w_{1}; \tau)=\sum_{j=1}^{p}F_{\hat{\mathcal{Y}}_{j}}(w_{1};
z_{1}; q_{\tau})$$
for some irreducible modules $\hat{W}_{j}$ and intertwining operators of type
$\binom{\hat{W}_{j}}{W_{1}\hat{W}_{j}}$ for $j=1, \dots, p$. 
Since $F_{\hat{\mathcal{Y}}_{j}}(w_{1};
z_{1}; q_{\tau})$, $j=1, \dots, p$, satisfy
(\ref{mod-inv-3}), we must have 
$$2\pi i\frac{d}{d\tau}S_{K-1}(w_{1}; \tau)
=S_{K-1}\left(\left(L(-2)-\sum_{k\in \mathbb{N}}G_{2k+2}(\tau)
L(2k)\right)w_{1}; \tau\right)$$
which together with (\ref{mod-inv-15}) gives
$$(2\pi i)^{2}KS_{K}(w_{1}; \tau)=0,$$
a contradiction.  

Now we  prove the case for any $n\ge 1$. This is the case where 
the method of Zhu \cite{Z}, further developed by Dong-Li-Mason
\cite{DLM} and Miyamoto \cite{M1}, cannot be used since there 
is no recurrence formula. 
We use induction. When $n=1$, the theorem is just proved. 
Assume that when $n=k$, the theorem is proved. 
We now prove the case $n=k+1$ using the genus-one associativity
(Theorem \ref{g1-asso}). By Theorem \ref{g1-asso}, we have 
\begin{eqnarray}\label{mod-inv-18}
\lefteqn{\overline{F}_{\mathcal{Y}_{1}, \dots, 
\mathcal{Y}_{n+1}}(w_{1}, \dots, w_{n+1};
z_{1}, \dots, z_{n+1}; \tau)}\nn
&&=\sum_{r\in \mathbb{R}}
\overline{F}_{\mathcal{Y}_{1}, \dots, \mathcal{Y}_{n-1},
\hat{\mathcal{Y}}_{n+1}}(w_{1}, 
\dots, w_{n-1}, P_{r}(\hat{\mathcal{Y}}_{n}(w_{n}, z_{n}-z_{n+1})
w_{n+1}); \nn
&&\quad\quad\quad\quad\quad\quad\quad\quad\quad\quad
\quad\quad\quad\quad\quad\quad\quad
z_{1}, \dots, z_{n-1}, z_{n+1}; \tau).
\end{eqnarray}
Using the induction assumption, 
we know that for any $r\in \mathbb{R}$, 
\begin{eqnarray*}
\lefteqn{\overline{F}_{\mathcal{Y}_{1}, \dots, \mathcal{Y}_{n-1},
\hat{\mathcal{Y}}_{n+1}}((\gamma\tau+\delta)^{-L(0)}w_{1}, 
\dots, (\gamma\tau+\delta)^{-L(0)}w_{n-1}, }\nn
&&\quad\quad\quad
(\gamma\tau+\delta)^{-L(0)}
P_{r}(\hat{\mathcal{Y}}(w_{n}, z_{n}-z_{n+1})
w_{n+1}); 
z'_{1}, \dots, z'_{n-1}, z'_{n+1}; \tau')\nn
&&
\end{eqnarray*}
is in 
$$\mathcal{F}_{w_{1}, \dots, w_{n-1}, 
P_{r}(\hat{\mathcal{Y}}_{n}(w_{n}, z_{n}-z_{n+1})
w_{n+1})}.$$
Thus by Theorem \ref{g1-asso} again,
\begin{eqnarray}\label{mod-inv-19}
\lefteqn{\sum_{r\in \mathbb{R}}
\overline{F}_{\mathcal{Y}_{1}, \dots, \mathcal{Y}_{n-1},
\hat{\mathcal{Y}}_{n+1}}((\gamma\tau+\delta)^{-L(0)}w_{1}, 
\dots, (\gamma\tau+\delta)^{-L(0)}w_{n-1},} \nn
&&\quad\quad\quad
(\gamma\tau+\delta)^{-L(0)}
P_{r}(\hat{\mathcal{Y}}(w_{n}, z_{n}-z_{n+1})
w_{n+1}); z'_{1}, \dots, z'_{n-1}, z'_{n+1}; \tau').\nn
&&
\end{eqnarray}
is absolutely convergent and is a linear combination
of absolutely convergent series of the form
\begin{eqnarray}\label{mod-inv-20}
\lefteqn{\sum_{r\in \mathbb{R}}
\overline{F}_{\breve{\mathcal{Y}}_{1}, \dots, \breve{\mathcal{Y}}_{n-1},
\bar{\mathcal{Y}}_{n+1}}(w_{1}, 
\dots, w_{n-1}, P_{r}(\bar{\mathcal{Y}}_{n}(w_{n}, z_{n}-z_{n+1})
w_{n+1}); }\nn
&&\quad\quad\quad\quad\quad\quad\quad\quad\quad
\quad\quad\quad\quad\quad\quad\quad\quad
z_{1}, \dots, z_{n-1}, z_{n+1}; \tau)
\end{eqnarray}
for suitable intertwining operators $\breve{\mathcal{Y}}_{1},
\dots, \breve{\mathcal{Y}}_{n-1}$, $\bar{\mathcal{Y}}_{n}$ and 
$\bar{\mathcal{Y}}_{n+1}$.
Moreover,  there exist suitable intertwining 
operators $\breve{\mathcal{Y}}_{n}$ and $\breve{\mathcal{Y}}_{n+1}$
such that (\ref{mod-inv-20}) is equal to 
\begin{equation}\label{mod-inv-21}
\overline{F}_{\breve{\mathcal{Y}}_{1}, \dots, \breve{\mathcal{Y}}_{n+1}}
(w_{1}, \dots, w_{n+1}; z_{1}, \dots, z_{n+1}; \tau)
\end{equation}
which is in $\mathcal{F}_{w_{1}, \dots,
w_{n+1}}$. Thus, (\ref{mod-inv-19}) as a linear combination of 
elements of the form (\ref{mod-inv-21}) is also in 
$\mathcal{F}_{w_{1}, \dots,
w_{n+1}}$. By (\ref{mod-inv-18}), we see that 
$$\overline{F}_{\mathcal{Y}_{1}, \dots, 
\mathcal{Y}_{n+1}}(w_{1}, \dots, w_{n+1};
z_{1}, \dots, z_{n+1}; \tau)$$
is also in $\mathcal{F}_{w_{1}, \dots,
w_{n+1}}$.
\epfv

\begin{rema}
{\rm Let $W_{1}=\cdots=W_{n}=V$. Assume that 
$\tilde{W}_{i}$ for $i=1, \dots, n$ are irreducible. 
Then any intertwining operator of type 
$\binom{\tilde{W}_{i}}{V\tilde{W}_{i+1}}$
is $0$ when $\tilde{W}_{i}$ is not isomorphic to
$\tilde{W}_{i+1}$ and is a multiple of the vertex operator 
$Y_{\tilde{W}_{i}}$ for the $V$-module $\tilde{W}_{i}$
when $\tilde{W}_{i}$ is isomorphic (and then is identified with)
$\tilde{W}_{i+1}$. In this case,  if we take
$w_{i}=(\mathcal{U}(1))^{-1}v_{i}\in V$
for $i=1, \dots, n$, Theorem \ref{mod-inv} 
gives the modular invariance result proved by Zhu in \cite{Z}. 
Similarly, if we take all but one of the $V$-modules 
$W_{1}, \dots, W_{n}$ to be $V$, Theorem \ref{mod-inv} 
gives the  generalization of Zhu's result 
by Miyamoto  in \cite{M1}.}
\end{rema}

\begin{rema}
{\rm If we replace  Conditions 1 and 3 by the conditions that 
for $n<0$, $V_{(n)}=0$ and every $V$-module is $C_{2}$-cofinite,
then the conclusions of Theorems \ref{summary} and \ref{mod-inv} 
are also true. }
\end{rema}

\begin{rema}
{\rm Geometrically, Theorems \ref{eqn-mod-inv} and \ref{mod-inv} and the
fact that the coefficients of the system (\ref{eqn1}) is doubly-periodic 
in $1$ and $\tau$
actually say that the space of the solutions of the system 
(\ref{eqn1})--(\ref{eqn2})  and the 
space $\mathcal{F}_{w_{1}, \dots,
w_{n+1}}$ are the spaces of (multivalued) global sections 
of some vector bundles with flat connections 
over the moduli space of genus-one Riemann surfaces with $n$ punctures
and standard local coordinates vanishing at these punctures.}
\end{rema}

\noindent {\small \sc Department of Mathematics, Rutgers University,
110 Frelinghuysen Rd., Piscataway, NJ 08854-8019 }

\noindent {\em E-mail address}: yzhuang@math.rutgers.edu

\end{document}